%% file: Shortest_Minkowski_billiard_trajectories_on_convex_bodies.tex
\documentclass[12pt]{amsart}

\usepackage[german,english]{babel}

\usepackage{dirtytalk}

\usepackage{amsmath}

\usepackage[mathx,matha]{mathabx}

\usepackage{amstext}

\usepackage{amssymb}

\usepackage{array}

\usepackage{amsthm}

\usepackage{dsfont}

\usepackage[latin1]{inputenc}

\usepackage[T1]{fontenc}

\usepackage{mathtools}

\usepackage{appendix}

\usepackage{xcolor}

\usepackage{ulem}

\usepackage[arrow, matrix, curve]{xy}

\usepackage{verbatim}

\usepackage{graphicx}
\graphicspath{{images/}}

\usepackage{subfigure}

\usepackage{hyperref}

\usepackage{cases}

\usepackage[update,prepend]{epstopdf}

\usepackage{bbold}

\usepackage[a]{esvect}

\usepackage{tikz}

\newtheorem{theorem}{Theorem}[section]

\newtheorem*{theorem*}{Theorem}

\theoremstyle{plain}

\newtheorem*{conjecture*}{Conjecture}

\newtheorem{assumption}[theorem]{Assumption}

\newtheorem{corollary}[theorem]{Corollary}

\newtheorem{proposition}[theorem]{Proposition}

\newtheorem{lemma}[theorem]{Lemma}

\theoremstyle{remark}

\newtheorem{condition}{Condition}

\theoremstyle{definition}

\newtheorem{definition}[theorem]{Definition}

\newtheorem{remark}[theorem]{Remark}

\def\eps{\varepsilon}

\def\bi{\begin{itemize}}

\def\ei{\end{itemize}}

\newcommand{\R}{\mathbb{R}}

\newcommand{\N}{\mathbb{N}}

\newcommand{\Z}{\mathbb{Z}}

\renewcommand{\dim}{{\rm dim}\,}

\renewcommand{\phi}{\varphi}

\DeclareMathOperator{\conv}{conv}

\DeclareOldFontCommand{\it}{\normalfont\itshape}{\mathit}

\newcommand{\bspm}{\left(\begin{smallmatrix}}\newcommand{\espm}{\end{smallmatrix}\right)}

\newcommand{\bpm}{\begin{pmatrix}}\newcommand{\epm}{\end{pmatrix}}

\def\bs{\begin{satz}}\def\es{\end{satz}}

\def\blem{\begin{lemma}}\def\elem{\end{lemma}}

\def\bthm{\begin{theorem}}\def\ethm{\end{theorem}}

\def\bcor{\begin{corollary}}\def\ecor{\end{corollary}}

\def\beq{\begin{equation}}\def\eeq{\end{equation}}

\def\beqq{\begin{equation*}}\def\eeqq{\end{equation*}}

\def\bal{\begin{align}}\def\eal{\end{align}}

\def\bpf{\begin{proof}}\def\epf{\end{proof}}

\def\bex{\begin{example}}\def\eex{\end{example}}

\def\brem{\begin{remark}}\def\erem{\end{remark}}

\def\bass{\begin{assumption}}\def\eass{\end{assumption}}

\def\bprop{\begin{proposition}}\def\eprop{\end{proposition}}

\def\bdefi{\begin{definition}}\def\edefi{\end{definition}}

\def\bcond{\begin{condition}}\def\econd{\end{condition}}

\def\bconj{\begin{conjecture*}}\def\econj{\end{conjecture*}}

\DeclareFontEncoding{FMS}{}{}

\DeclareFontSubstitution{FMS}{futm}{m}{n}

\DeclareFontEncoding{FMX}{}{}

\DeclareFontSubstitution{FMX}{futm}{m}{n}

\DeclareSymbolFont{fouriersymbols}{FMS}{futm}{m}{n}

\DeclareSymbolFont{fourierlargesymbols}{FMX}{futm}{m}{n}

\DeclareMathDelimiter{\VERT}{\mathord}{fouriersymbols}{152}{fourierlargesymbols}{147}

\def\bi{\begin{itemize}}

\def\ei{\end{itemize}}

\def\ben{\begin{enumerate}}

\def\een{\end{enumerate}}

\usepackage{enumitem}

\usepackage{tcolorbox}

\newtcolorbox{implementation}[2][]{colframe=blue!75!black,colbacktitle=green!10!white,colback=green!10!white,coltitle=green!75!black,title={#2},fonttitle=\bfseries,#1}

\begin{document}

\title[Shortest Minkowski billiard trajectories on convex bodies]{Shortest Minkowski billiard trajectories on convex bodies}

\author{Stefan Krupp and Daniel Rudolf}


\date{\today}

\maketitle

\begin{abstract}
We rigorously investigate closed Minkowski/Finsler billiard trajectories on $n$-dimensional convex bodies. We outline the central properties in comparison and differentiation from the Euclidean special case and establish two main results for length-minimizing closed Minkowski/Finsler billiard trajectories: one is a regularity result, the other is of geometric nature. Building on these results, we develop an algorithm for computing length-minimizing closed Minkowski/Finsler billiard trajectories in the plane.
\end{abstract}

\section{Introduction and main results}\label{Sec:Introminkbill}

Minkowski/Finsler billiards are the natural extensions of Euclidean billiards to the Minkowski/Finsler setting.

Euclidean billiards are associated to the local Euclidean billiard reflection rule: The angle of reflection equals the angle of incidence (here, we assume that the relevant normal vector as well as the incident and the reflected ray lie in the same two-dimensional affine flat). This local Euclidean billiard reflection rule follows from the global least action principle, colloquially expressed: A billiard trajectory between a pair of points of the billiard table minimizes the Euclidean length among all paths connecting these two points via a reflection at the billiard table boundary.

In Finsler geometry, the notion of length of vectors in $\R^n$, $n\geq 1$, is given by a convex body $T\subset\R^n$, i.e., a compact convex set in $\R^n$ which has the origin in its interior (in $\R^n$). The Minkowski functional
\beqq \mu_{T}(x)=\min\{t\geq 0 : x\in tT\},\; x\in\R^n,\eeqq
determines the distance function, where we recover the Euclidean setting when $T$ is the $n$-dimensional Euclidean unit ball. Then, heuristically, billiard trajectories are defined via the global least action principle with respect to $\mu_T$--we specify this in a moment--, because in Finsler geometry, there is no useful notion of angles.

There is generally much interest in the study of billiards: Problems in almost every mathematical field can be related to problems in mathematical billiards, see for example \cite{Gutkin2012}, \cite{Katok2005} and \cite{Tabachnikov2005} for comprehensive surveys. Euclidean billiards in the plane have been investigated very intensively. Nevertheless, so far not much is known about Euclidean billiards on higher dimensional billiard tables. But even much less is known for Minkowski/Finsler billiards. Although the applications are numerous and important in many different areas, to the authors' knowledge, \cite{AkopBal2015, ElApp2015, AkopKar2019, AkopSchTab2020, AlkoumiSchlenk2014, ArtOst2012, ArtFlorOstRos2017, BlagHarTabZieg2017, GutkinTabachnikov2002, Radnovic2003} are the only publications concerning different aspects of Minkowski/Finsler billiards so far.

After Minkowski/Finsler billiards were introduced in \cite{GutkinTabachnikov2002}, their study was intensified when in \cite{ArtOst2012} the relationship to the EHZ-capacity of convex Lagrangian products in $\R^{2n}$ was proven. When studied isoperimetric Minkowski/Finsler billiard inequalities, this opened the connection to Viterbo's conjecture (cf.\;\cite{Viterbo2000}) within symplectic geometry. With \cite{ArtKarOst2013}, this consequently also allowed to analyze the famous Mahler conjecture (cf.\;\cite{Mahler1939}) from convex geometry.

In this paper, we prove fundamental properties of Minkowski/Finsler billiards. Particular attention is paid to length-minimizing closed Minkowski/Finsler billiard trajectories. As part of the investigation of the latter, we state two main results: one mainly is a regularity result, the other is geometric in nature. Together they can be seen as the generalization of Theorem 1.2 in \cite{KruppRudolf2020} to the Minkowski/Finsler setting. Based on these results, we derive an algorithm for computing length-minimizing closed Minkowski/Finsler billiard trajectories in the plane.

Before we state the main results of this paper, let us precisely define Minkowski/Finsler billiards, while we suggest from now on to call them just \textit{Minkowski billiards} following \cite{ArtOst2012} and \cite{ArtFlorOstRos2017}.

In the further course of this paper, we will see that it makes sense to differentiate between weak and strong Minkoswski billiard trajectories. We begin by introducing weak Minkowski billiard trajectories.

\begin{figure}[h!]
\centering
\def\svgwidth{300pt}
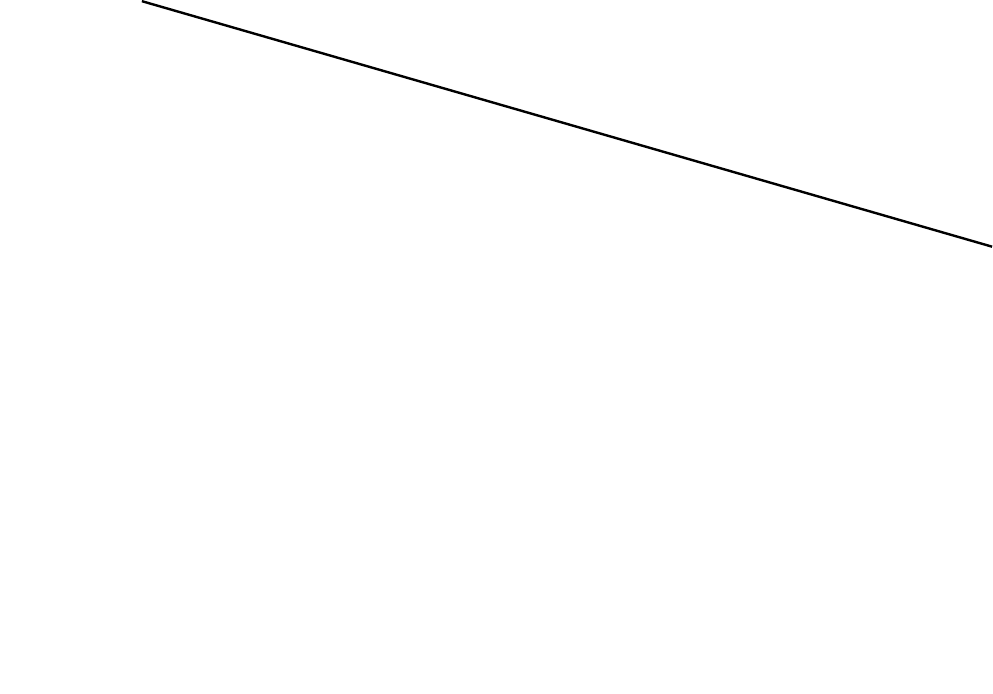
\caption[The weak Minkowski billiard relection rule]{The weak Minkowski billiard reflection rule: $q_j$ minimizes \eqref{eq:Minkowskipolar} over all $\widebar{q}_j\in H_j$, where $H_j$ is a $K$-supporting hyperplane through $q_j$.}
\label{img:Stossregelminkbill}
\end{figure}

\bdefi[Weak Minkowski billiards]\label{Def:weakMinkowskiBilliards}
Let $K\subset\R^n$ be a convex body which from now on we call the \textit{billiard table}. Let $T\subset\R^n$ be another convex body and $T^\circ\subset\R^n$ its polar body. We say that a closed polygonal curve\footnote{For the sake of simplicity, whenever we talk of the vertices $q_1,...,q_m$ of a closed polygonal curve, we assume that they satisfy $q_j\neq q_{j+1}$ and $q_j$ is not contained in the line segment connecting $q_{j-1}$ and $q_{j+1}$ for all $j\in\{1,...,m\}$. Furthermore, whenever we settle indices $1,...,m$, then the indices in $\Z$ will be considered as indices modulo $m$.\label{foot:polygonalline}} with vertices $q_1,...,q_m$, $m\in \N_{\geq 2}$, on the boundary of $K$ (denoted by $\partial K$) is a \textit{closed weak $(K,T)$-Minkowski billiard trajectory} if for every $j\in \{1,...,m\}$, there is a $K$-supporting hyperplane $H_j$ through $q_j$ such that $q_j$ minimizes
\beq \mu_{T^\circ}(\widebar{q}_j-q_{j-1})+\mu_{T^\circ}(q_{j+1}-\widebar{q}_j),\label{eq:Minkowskipolar}\eeq
over all $\widebar{q}_j\in H_j$ (cf.\;Figure \ref{img:Stossregelminkbill}). We encode this closed weak $(K,T)$-Minkowski billiard trajectory by $(q_1,...,q_m)$ and call its vertices \textit{bouncing points}. Its \textit{$\ell_T$-length}\footnote{This length-definition can be generalized to any closed polygonal curve.} is given by
\beqq \ell_T (( q_1,...,q_m )) = \sum\limits_{j = 1}^m \mu_{T^\circ}(q_{j+1}-q_j).\eeqq
\edefi

We call a boundary point $q\in\partial K$ \textit{smooth} if there is a unique $K$-supporting hyperplane through $q$. We say that $\partial K$ is \textit{smooth} if every boundary point is smooth (we also say $K$ is smooth while we actually mean $\partial K$).

Concerning the definition of weak Minkowski billiards, we remark that, in general, the $K$-supporting hyperplanes $H_j$ are not uniquely determined. This is only the case for smooth and strictly convex $T$ (cf.\;Proposition \ref{Prop:onlyonehyperplane}). 

We note that the weak Minkowski billiard reflection rule defined in Definition \ref{Def:weakMinkowskiBilliards} does not only generalize the Euclidean billiard reflection rule to Minkowski/Finsler geometries, it also extends the classical understanding of billiard trajectories to non-smooth billiard table boundaries\footnote{Classical billiard trajectories are usually understood as trajectories with bouncing points in smooth boundary points (billiard table gangs) while they stop in non-smooth boundary points (billiard table holes).}. In the latter case, there can exist infinitely many different $K$-supporting hyperplanes through non-smooth bouncing points at the billiard table boundary, and consequently, from a constructive viewpoint, the weak Minkowski billiard reflection rule may produce different bouncing points following two already known consecutive ones.

\begin{figure}[h!]
\centering
\def\svgwidth{420pt}
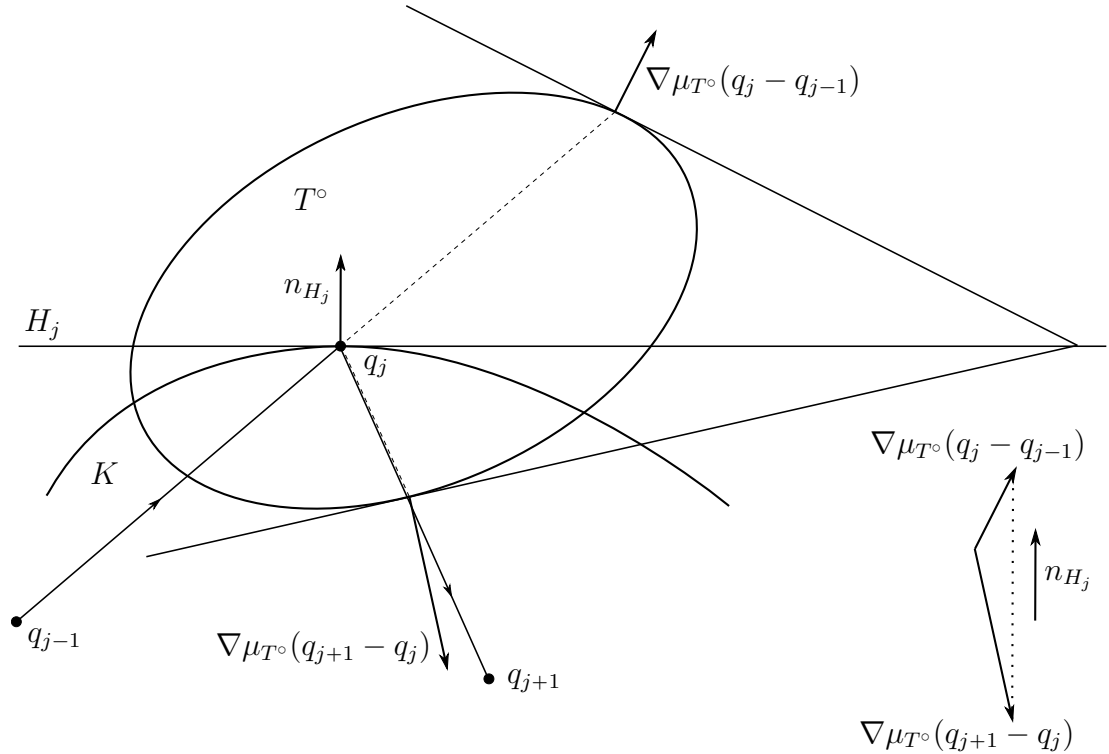
\caption[Illustration I of the reformulated Minkowski billiard reflection rule]{$T^\circ$ is a smooth and strictly convex body in $\R^2$ and its boundary plays the role of the indicatrix, i.e., the set of vectors of unit Finsler (with respect to $T^\circ$) length, which therefore is an $1$-level set of $\mu_{T^\circ}$. Note that the two $T^\circ$-supporting hyperplanes intersect on $H_j$ due to the condition $\nabla\mu_{T^\circ}(q_j-q_{j-1})-\nabla\mu_{T^\circ}(q_{j+1}-q_j)= \mu_j n_{H_j}$.}
\label{img:ReflectionRule1}
\end{figure}

In the case when $T^\circ$ is smooth and strictly convex, the definition of weak Minkowski billiards yields a geometric interpretation of the billiard reflection rule: On the basis of Lagrange's multiplier theorem, one derives the condition
\beqq \nabla_{\widebar{q}_j}\Sigma_j(\widebar{q}_j)_{\vert \widebar{q}_j=q_j}=\nabla\mu_{T^\circ}(q_j-q_{j-1})-\nabla\mu_{T^\circ}(q_{j+1}-q_j)= \mu_j n_{H_j},\eeqq
where $\mu_j>0$, since the strict convexity of $T^\circ$ implies
\beqq \nabla\mu_{T^\circ}(q_j-q_{j-1}) \neq \nabla\mu_{T^\circ}(q_{j+1}-q_j),\eeqq
and where $n_{H_j}$ is the outer unit vector normal to $H_j$. This implies that the weak Minkowski billiard reflection rule can be illustrated as within Figure \ref{img:ReflectionRule1}. For smooth, strictly convex, and centrally symmetric $T^ \circ\subset\R^n$, this interpretation is due to \cite[Lemma 3.1, Corollary 3.2 and Lemma 3.3]{GutkinTabachnikov2002} (this interpretation has also been referenced in \cite{AlkoumiSchlenk2014}). For the extension to just smooth and strictly convex $T^\circ\subset\R^n$, it is due to \cite[Lemma 2.1]{BlagHarTabZieg2017}. However, from the constructive point of view, this interpretation has its limitations.

We come to the definition of strong Minkowski billiards. 

\bdefi[Strong Minkowski billiards]\label{Def:strongMinkowskiBilliards}
Let $K,T\subset\R^n$ be convex bodies. We say that a closed polygonal curve $q$ with vertices $q_1,...,q_m$, $m\in \N_{\geq 2}$, on $\partial K$ is a \textit{closed strong $(K,T)$-Minkowski billiard trajectory} if there are points $p_1,...,p_m$ on $\partial T$ such that
\beq \begin{cases} q_{j+1}-q_j \in N_T(p_j), \\ p_{j+1}-p_j \in - N_K(q_{j+1})\end{cases}\label{eq:System}\eeq
is satisfied for all $j\in\{1,...,m\}$. We call $p=(p_1,...,p_m)$ a \textit{closed dual billiard trajectory in $T$}. We say that $q=(q_1,...,q_m)$ is a closed strong $(K,T)$-Minkowski billiard trajectory with respect to $H_1,...,H_m$, if $H_1,...,H_m$ are the $K$-supporting hyperplanes through $q_1,...,q_m$ which are normal to $n_K(q_1),...,n_K(q_m)$, where
\beq p_{j+1}-p_{j}=-\mu_{j+1} n_K(q_{j+1}),\; \mu_{j+1} \geq 0,\label{eq:weakstrongcondition}\eeq
for all $j\in\{1,...,m\}$.
\edefi

\begin{figure}[h!]
\centering
\def\svgwidth{420pt}
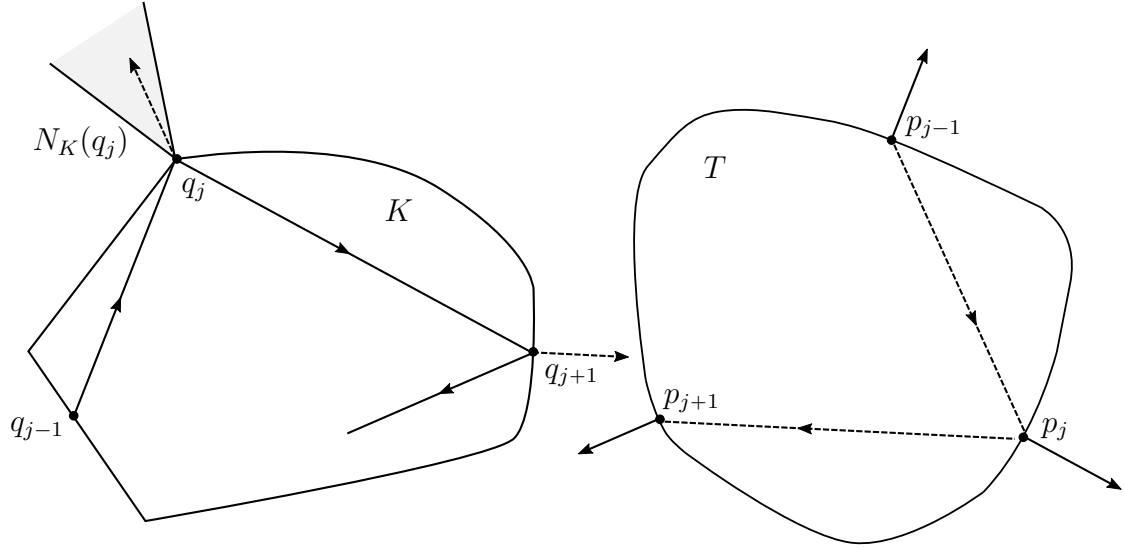
\caption[Illustration II of the reformulated Minkowski billiard reflection rule]{The pair $(q,p)$ fulfills \eqref{eq:System}, namely: $q_j-q_{j-1}\in N_T(p_{j-1})$, $q_{j+1}-q_j\in N_T(p_j)$, $p_j-p_{j-1}\in -N_K(q_j)$, and $p_{j+1}-p_j\in - N_K(q_{j+1})$.}
\label{img:ReflectionRule2}
\end{figure}

This definition appeared implicitly in \cite[Theorem 7.1]{GutkinTabachnikov2002}, then later the first time explicitly in \cite{ArtOst2012}. It yields a different interpretation of the billiard reflection rule. Without requiring a condition on $T$, the billiard reflection rule can be represented as within Figure \ref{img:ReflectionRule2}. From the constructive point of view, this interpretation is much more appropriate in comparison to the one for weak Minkowski billiards.

We remark that, in general, the closed dual billiard trajectory in $T$ is not uniquely determined. This is only the case when $T$ is strictly convex (cf.\;Proposition \ref{Prop:dualbilliardtrajunique}). Further, we remark that under the condition that $T$ is strictly convex and smooth, the closed dual billiard trajectory $p$ in Definition \ref{Def:strongMinkowskiBilliards} also is a closed Minkowski billiard trajectory--we refer to Proposition \ref{Prop:dualbilliard} for the precise statement.

The natural follow-up question concerns the relationship between weak and strong Minkowski billiards. The following theorem gives an answer to this question.

\bthm\label{Thm:weakstrongbilliards}
Let $K,T\subset\R^n$ be convex bodies. Then, every closed strong $(K,T)$-Minkowski billiard trajectory is a weak one, more precisely, if $q=(q_1,...,q_m)$ is a closed strong $(K,T)$-Minkowski billiard trajectory with respect to $H_1,...,H_m$, then $q$ fulfills the weak Minkowski billiard reflection rule with respect to $H_1,...,H_m$.

If $T$ is strictly convex, then every closed weak $(K,T)$-Minkowski billiard trajectory is a strong one, more precisely, if $q=(q_1,...,q_m)$ is a closed weak $(K,T)$-Minkowski billiard trajectory fulfilling the weak Minkowski billiard reflection rule with respect to the $K$-supporting hyperplanes $H_1,...,H_m$, then $q$ is a closed strong $(K,T)$-Minkowski billiard trajectory with respect to $H_1,...,H_m$.
\ethm

Theorem \ref{Thm:weakstrongbilliards} is sharp in the following sense: One can construct convex bodies $K,T\subset\R^n$ (where $T$ is not strictly convex) and a closed weak $(K,T)$-Minkowski billiard trajectory which is not a strong one (cf.\;Example A in Section \ref{Sec:Examples}).

In the following, if the risk of confusion is excluded, we will call strong Minkowski billiards trajectories just Minkowski billiard trajectories. Although the connection to the least action principle, which is the basis for the definition of the weak version, is lost for non-strictly convex bodies, we make the strong version to the center of our investigations, since one can use them to compute the EHZ-capacity of convex Lagrangian products (cf.\;\cite{Rudolf2022}).

In the following theorem, the main properties of closed Minkowski billiard trajectories are collected. We recall that for convex body $K\subset\R^n$, we denote by $F(K)$ the set of subsets of $\R^n$ which cannot be translated into $\mathring{K}$ and by $F^{cp}_{n+1}(K)$ the set closed polygonal curves in $\R^n$ with at most $n+1$ vertices which cannot be translated into $\mathring{K}$. Furthermore, we denote by $M_{n+1}(K,T)$ the set of closed $(K,T)$-Minkowski billiard trajectories with at most $n+1$ bouncing points.

\bthm\label{Thm:collectionproperties}
Let $K,T\subset\R^n$ be convex bodies.
\begin{itemize}
\item[(i)] [Proposition \ref{Prop:lengthdualbilliard}] Let $q$ be a closed $(K,T)$-Minkowski billiard trajectory with dual billiard trajectory $p$. Then, we have
\beqq \ell_T(q)=\ell_{-K}(p).\eeqq
\item[(ii)] [Proposition \ref{Prop:dualbilliard}] Let $T$ be strictly convex and smooth. Let $q$ be a closed $(K,T)$-Minkowski billiard trajectory with dual billiard trajectory $p$. Then, $p$ is a closed $(T,-K)$-Minkowski billiard trajectory.
\item[(iii)] [Proposition \ref{Prop:notranslation}] Let $q$ be a closed $(K,T)$-Minkowski billiard trajectory with dual billiard trajectory $p$. Then, we have
\beqq q\in F(K)\; \text{ and } \; p\in F(T).\eeqq
\item[(iv)] [Proposition \ref{Prop:GenBezdek}] Let $q=(q_1,...,q_m)$ be a closed $(K,T)$-Minkowski billiard trajectory with respect to $H_1,...,H_m$ and let $U$ be the inclusion minimal linear subspace of $\R^n$ containing the outer unit vectors $n_K(q_1),...,n_K(q_m)$ which are normal to $H_1,...,H_m$. Then, there is a selection $\{i_1,...,i_{\dim U +1}\}\subseteq \{1,...,m\}$ such that
\beqq \{q_{i_1},...,q_{i_{\dim U +1}}\}\in F(K).\eeqq
\item[(v)] [Theorem \ref{Thm:onetoone}] Let $T$ be strictly convex. Then, every $\ell_T$-minimizing closed $(K,T)$-Minkowski billiard trajectory is an $\ell_T$-minimizing element of $F_{n+1}^{cp}(K)$ and, conversely, every $\ell_T$-minimizing element of $F_{n+1}^{cp}(K)$ can be translated in order to be an $\ell_T$-minimizing closed $(K,T)$-Minkowski billiard trajectory.

Especially, one has
\beq \min_{q\in F^{cp}_{n+1}(K)}\ell_T(q) = \min_{q \in M_{n+1}(K,T)} \ell_T(q).\label{eq:onetoone0}\eeq
\end{itemize}
\ethm

We note that in \cite[Theorem 2.2]{Rudolf2022}, we actually prove that \eqref{eq:onetoone0} holds without the condition of strict convexity of $T$ (in the general case, the other part of the statement in Theorem \ref{Thm:collectionproperties}(v) holds slightly changed). However, for the purposes of this paper, this formulation is enough for us.

The main interest we pursue in the following part of this paper is the investigation of the subsequent questions (in analogy to the questions in \cite{KruppRudolf2020}). Given an $\ell_T$-minimizing closed $(K,T)$-Minkowski billiard trajectory--depending on the conditions on $K\subset\R^n$ and $T\subset\R^n$:
\begin{itemize}
\item[(a)] What can be said concerning the number of bouncing points?
\item[(b)] What can be said concerning the convex cone which is spanned by the normal vectors related to the Minkowski billiard reflection rule? Is it a subspace? What is its dimension?
\item[(c)] What can be said concerning the dimension of the inclusion minimal affine section of $K$ containing this billiard trajectory?
\item[(d)] What can be said concerning the regularity of the bouncing points? Are they smooth? Are they smooth within the inclusion minimal section of $K$ containing this billiard trajectory?
\end{itemize}
All these questions are important in order to develop an algorithm for computing $\ell_T$-length-minimizing closed $(K,T)$-Minkowski billiard trajectories.


Our first main result is the following regularity result:

\bthm\label{Thm:RegularityResult}
Let $K,T\subset\R^n$ be convex bodies, where $T$ is additionally assumed to be strictly convex and smooth. Let $q=(q_1,...,q_m)$ be an $\ell_T$-minimizing closed $(K,T)$-Minkowski billiard trajectory (which fulfills the Minkowski billiard reflection rule with respect to $H_1,...,H_m$) and let $U$ be the convex cone spanned by the outer unit vectors $n_K(q_1),...,n_K(q_m)$ which are normal to $H_1,...,H_m$. Then, $U$ is a linear subspace of $\R^n$ with
\beqq \dim U = m-1 \eeqq
and
\beq \dim (N_K(q_j)\cap U)=1\label{eq:regularityformula}\eeq
for all $j\in\{1,...,m\}$.
\ethm

Since in general $\dim U\leq n$, it follows that
\beqq m\leq n+1.\eeqq

In fact, \eqref{eq:regularityformula} is a \textit{regularity} result: If $m=n+1$, meaning that $U=\R^n$, then \eqref{eq:regularityformula} becomes
\beqq \dim(N_K(q_j))=1\eeqq
for all $j\in\{1,...,m\}$, i.e., $q$ is regular, meaning that all bouncing points of $q$ are smooth boundary points of $K$.

Some special cases of Theorem \ref{Thm:RegularityResult} were already known: For $T$ equals the Euclidean unit ball in $\R^n$, Theorem \ref{Thm:RegularityResult} coincides with \cite[Theorem 1.2]{KruppRudolf2020}, since by \cite[Proposition 2.6]{KruppRudolf2020}, in this Euclidean case, $U$ equals $V_0$ which is the underlying linear subspace of the affine subspace $V\subseteq\R^n$ which is chosen such that $K\cap V$ is the inclusion minimal affine section of $K$ containing $q$. In the situation of Theorem \ref{Thm:RegularityResult}, for $\ell_T$-minimizing closed $(K,T)$-Minkowski billiard trajectories $q=(q_1,...,q_m)$, it has been proven in \cite{ElApp2015} that $m$ is bounded by $n+1$, and in \cite[Theorem 4.8]{AkopBal2015}, that the $\ell_T$-minimizing closed $(K,T)$-Minkowski billiard trajectories with $n+1$ bouncing points are regular. In \cite{AlkoumiSchlenk2014}, it has been proven the statement of Theorem \ref{Thm:RegularityResult} for $n=2$.

As in the less general Euclidean case, Theorem \ref{Thm:RegularityResult} refutes the presumption which at first appears to be intuitively correct, that every $\ell_T$-minimizing closed $(K,T)$-Minkowski billiard trajectory with more than two bouncing points is regular within the inclusion minimal section of $K$ containing this billiard trajectory.

For the sharpness of Theorem \ref{Thm:RegularityResult}, we refer on the one hand to the sharpness of \cite[Theorem 2.1]{KruppRudolf2020}. There, for $T$ equals the Euclidean unit ball, we showed that an $\ell_T$-minimizing closed $(K,T)$-Minkowski billiard trajectory may not be regular within the inclusion minimal sections of $K$ containing the billiard trajectory (which in this case is a translate of $U$). This can even appear for the unique $\ell_T$-minimizing closed $(K,T)$-Minkowski billiard trajectory. Then, we also showed that an $\ell_T$-minimizing closed $(K,T)$-Minkowski billiard trajectory can have bouncing points in vertices as well as in more than $0$-dimensional faces of $K$.

On the other hand, we can prove that in Theorem \ref{Thm:RegularityResult}, the smoothness of $T$ is a necessary condition. For that, we refer to Example F. Furthermore, we can show that within the weaker situation of weak Minkowski billiard trajectories, in general, the statement of Theorem \ref{Thm:RegularityResult} is not true without requiring the strict convexity of $T$. For that, we refer to Example G (cf.\;Section \ref{Sec:Examples}). In this latter example, we will see that without requiring the strict convexity of $T$, even the notion of $\ell_T$-minimizing closed weak $(K,T)$-Minkowski billiard trajectories does not make sense: There are configurations for which such minimizers do not exist.

The second main result generalizes a property of length-minimizing closed Euclidean billiard trajectories on $K$ (cf.\;\cite[Theorem 2.1]{KruppRudolf2020}) to the Minkowski/Finsler setting $(K,T)$ under the condition that both $K$ and $T$ are strictly convex and smooth:

\bthm\label{Thm:maximxallyspanning}
Let $K,T\subset\R^n$ be strictly convex and smooth bodies. Let $q=(q_1,...,q_m)$ be an $\ell_T$-minimizing closed $(K,T)$-Minkowski billiard trajectory and $V$ an affine subspace of $\R^n$ such that $K\cap V$ is the inclusion minimal affine section of $K$ containing $q$. Then, we have
\beqq \dim V=m-1,\eeqq
i.e., $q$ is maximally spanning, meaning that the dimension of the convex hull of the bouncing points $q_1,...,q_m$ is $m-1$.
\ethm

From Theorem \ref{Thm:RegularityResult} and a consideration in the context of Theorem \ref{Thm:maximxallyspanning}, we immediately derive the following corollary for $n=2$:

\bcor\label{Cor:RegularityResult}
Let $K,T\subset\R^2$ be convex bodies, where $T$ is additionally assumed to be strictly convex and smooth. Then, every $\ell_T$-minimizing closed $(K,T)$-Minkowski billiard trajectory $q=(q_1,...,q_m)$ has two or three bouncing points, i.e., $m\in\{2,3\}$, where in the latter case the billiard trajectory is regular. Furthermore, when $q$ fufills the Minkowski billiard reflection rule with respect to $H_1,...,H_m$, then we have:
\begin{itemize}
\item[(i)] If $U$ is the convex cone spanned by the unit vectors $n_K(q_1),...,n_K(q_m)$ which are normal to $H_1,...,H_m$, then $U$ is a linear subspace of $\R^2$ with
\beqq \dim U=m-1.\eeqq
\item[(ii)] If $V\subseteq\R^2$ is an affine subspace such that $K\cap V$ is the inclusion minimal affine section of $K$ containing $q$, then
\beqq \dim V = m-1.\eeqq
\end{itemize}
\ecor

Here, $(i)$ follows from Theorem \ref{Thm:RegularityResult} and $(ii)$ from the fact that, by definition, both $2$- and $3$-bouncing Minkowski billiard trajectories are maximally spanning (so, for $(ii)$, the strict convexity of $K$, as required in Theorem \ref{Thm:maximxallyspanning} for general dimension, is not necessary).

Corollary \ref{Cor:RegularityResult} can be used for the construction of $\ell_T$-minimizing closed $(K,T)$-Minkowski billiard trajectories when $K\subset\R^2$ is a convex polytope. In Section \ref{Subsec:GeneralConstruction}, we describe an algorithm, while in Section \ref{Subsec:Implementation}, we provide a detailed survey concerning the efficiency and methods used within the implementation. We remark especially that our implemented algorithm is the solution to the open problem of finding the Fagnano triangle in a Minkowski triangle stated in \cite{AlkoumiSchlenk2014}.

Let us briefly give an overview of the structure of this paper: In Section \ref{Sec:PreliminariesConvexGeometry}, we enumerate properties of the polar body and the Minkowski functional that we will use repeatedly within this paper. In Section \ref{Sec:Minkowskibilliards}, we prove fundamental properties of Minkowski billiard trajectories which we utilize in Sections \ref{Sec:Proof1} and \ref{Sec:Proof2} for the proofs of Theorems \ref{Thm:RegularityResult} and \ref{Thm:maximxallyspanning}, respectively. In Section \ref{Sec:Examples}, we give various examples in order to show the sharpness of the statements in Section \ref{Sec:Minkowskibilliards} as well as the sharpness of the main results. In Section \ref{Sec:Construction}, we discuss the algorithm for computing $\ell_T$-minimizing closed $(K,T)$-Minkowski billiard trajectories for $K,T\subset\R^2$. In Section \ref{Sec:ObtuseTriangle}, we present a note on the existence of closed regular Minkowski billiard trajectories with three bouncing points in obtuse triangles.

\section{Preliminaries from convex geometry}\label{Sec:PreliminariesConvexGeometry}

In this section, we collect some useful properties of the polar body and the Minkowski functional. Since these properties are well-known, we will just state them and refer for the proofs to the usual literature on this topic.

Let $T\subset\R^n$ be a convex body. Then, the \textit{polar body} of $T$ is
\beqq T^\circ = \left\{x\in\R^n : \langle x,y\rangle \leq 1 \text{ for all }y\in T\right\}\subset\R^n.\eeqq
The polar body satisfies the following properties:

\bprop\label{Prop:polarbody}
Let $P,Q\subset\R^n$ be convex bodies. Then:
\begin{itemize}
\item[(i)] $P^\circ$ is in $\mathcal{C}(\R^n)$.
\item[(ii)] For $c\neq 0$ we have $(cP)^\circ = \frac{1}{c}P^\circ$.
\item[(iii)] If $P \subseteq Q$, then $P^\circ \supseteq Q^\circ$.
\item[(iv)] It is $(P^\circ)^\circ = P$.
\end{itemize}
\eprop

The \textit{Minkowski functional} $\mu_{T}$, defined by
\beqq \mu_{T}(x):=\min\{t\geq 0 : x\in tT\},\; x\in\R^n,\eeqq
determines a distance function, where we recover the Euclidean distance when $T$ is the $n$-dimensional Euclidean unit ball. The \textit{support function} $h_T$ of $T$ is given by
\beqq h_T(x):=\max\{\langle x,y\rangle : y\in T\}, \; x\in\R^n.\eeqq
The following lemma clarifies the connection between the Minkowski functional and the support function and will be useful throughout many proofs of the following sections:

\bprop\label{Prop:Minkowskisupport}
Let $T\subset\R^n$ be a convex body.
\begin{itemize}
\item[(i)] One has
\beqq h_T(x)=\mu_{T^\circ}(x)\eeqq
for all $x\in\R^n$.
\item[(ii)] The following equivalence holds for all $x\in\R^n$: One has
\beqq h_T(x)=\langle x,y\rangle\; \Leftrightarrow\; x\in N_T(y)\eeqq
under the constraint $y\in\partial T$.
\end{itemize}
\eprop

The next proposition collects properties of the Minkowski functional:

\bprop\label{Prop:IntroMinkowski}
Let $S,T\subset\R^n$ be convex bodies.
\begin{itemize}
\item[(i)] If $T$ is additionally assumed to be strictly convex and if
\beqq x,y\in\R^n\setminus \{0\} \; \text{ with } \; x \neq \lambda y\; \forall\lambda\in\R,\eeqq
then
\beq \mu_{T^\circ}(x+y)<\mu_{T^\circ}(x) + \mu_{T^\circ}(y).\label{eq:IntroMinkowski1}\eeq
We note that, when just requiring convexity of $T$, one has $\leq$ in \eqref{eq:IntroMinkowski1}.
\item[(ii)] With $S\subseteq T$ we have $T^\circ \subseteq S^\circ$. This implies
\beqq \mu_{S^\circ}(x)\leq \mu_{T^\circ}(x)\quad \forall x\in\R^n.\eeqq
\item[(iii)] For $c\neq 0$ we have
\beqq \mu_{T^\circ}(cx)=\mu_{(cT)^\circ}(x)=c\mu_{T^\circ}(x) \quad \forall x\in\R^n.\eeqq
\item[(iv)] The map
\beqq \mu_{T^\circ}: (\R^n,|\cdot|)\rightarrow (\R_{\geq 0},|\cdot|)\eeqq
is continuous.
\end{itemize}
\eprop

\section{Properties of closed $(K,T)$-Minkowski billiard trajectories}\label{Sec:Minkowskibilliards}

We begin with the following lemma:

\blem\label{Lem:strictconvexitynormalcones}
Let $T\subset\R^n$ be a strictly convex body and $p_i,p_j\in\partial T$. Then one has the equivalence
\beqq q_i\neq q_j \Leftrightarrow N_T(q_i)\cap N_T(q_j)=\{0\}.\eeqq
\elem

\bpf
From
\beqq N_T(q_i)\cap N_T(q_j)=\{0\},\eeqq
together with
\beqq \{0\}\subsetneqq N_T(q_i) \; \text{ and } \; \{0\}\subsetneqq N_T(q_j),\eeqq
it directly follows $q_i\neq q_j$.

Let $q_i\neq q_j$. If there is a nonzero
\beqq n\in N_T(q_i)\cap N_T(q_j),\eeqq
then it follows from the definition of strict convexity that
\beqq \langle n,z_1-q_i\rangle < 0 \quad \forall z_1\in T\eeqq
and
\beqq \langle n,z_2-q_j\rangle < 0 \quad \forall z_2\in T.\eeqq
Choosing $z_1=q_j$ and $z_2=q_i$ yields
\beqq 0 > \langle n,q_j-q_i\rangle,\eeqq
and therefore
\beqq 0 > \langle n,q_i-q_j\rangle = - \langle n,q_j-q_i\rangle >0,\eeqq
a contradiction. Therefore, it follows
\beqq N_T(q_i)\cap N_T(q_j)=\{0\}.\eeqq
\epf

The subsequent proposition clarifies the uniqueness of closed dual billiard trajectories.

\bprop\label{Prop:dualbilliardtrajunique}
Let $K,T\subset\R^n$ be convex bodies, where $T$ is additionally assumed to be strictly convex. Let $q$ be a closed strong $(K,T)$-Minkowski billiard trajectory. Then, the closed dual billiard trajectory $p$ in $T$ is uniquely determined.
\eprop

\bpf
Referring to \eqref{eq:System}, $p=(p_1,...,p_m)$ fulfills
\beqq q_{j+1}-q_j \in N_T(p_j)\eeqq
for all $j\in\{1,...,m\}$. Using Lemma \ref{Lem:strictconvexitynormalcones} ($T$ is assumed to be strictly convex), this implies that $p_1,...,p_m$ are uniquely determined.
\epf

Without requiring the strict convexity of $T$, $p$ is not necessarily uniquely determined (cf.\;Example B in Section \ref{Sec:Examples}).

We proceed with the proof of Theorem \ref{Thm:weakstrongbilliards}, which clarifies the relationship between weak and strong Minkowski billiards.

\bpf[Proof of Theorem \ref{Thm:weakstrongbilliards}]
We first prove that every closed strong Minkowski billiard trajectory is a weak one. For that, let $T$ be an arbitrary convex body in $\R^n$. Given a closed strong $(K,T)$-Minkowski billiard trajectory $q=(q_1,...,q_m)$ together with its dual billiard trajectory $p=(p_1,...,p_m)$. We let $n_K(q_1),...,n_K(q_m)$ be the unit vectors in $N_K(q_1),...,N_K(q_m)$ for which
\beqq p_{j+1}-p_{j}=-\mu_{j+1}n_K(q_{j+1}),\; \mu_{j+1}\geq 0,\eeqq
holds for all $j\in\{1,...,m\}$. Now, let $j\in\{1,...,m\}$ be arbitrarily chosen and let $H_j$ be the $K$-supporting hyperplane through $q_j$ which is normal to $n_K(q_j)$. Then, one has
\beqq \langle q_j-q_j^*,p_j-p_{j-1}\rangle =0\eeqq
for all $q_j^*\in H_j$ (since $p_j-p_{j-1}=-\mu_j n_K(q_j)$) and therefore together with
\beqq q_j-q_{j-1}\in N_T(p_{j-1})\;\text{ and }\; q_{j+1}-q_j\in N_T(p_j)\eeqq
and Proposition \ref{Prop:Minkowskisupport} that
\allowdisplaybreaks{\begin{align*}
\Sigma_j(q_j)=&\mu_{T^\circ}(q_j-q_{j-1})+\mu_{T^\circ}(q_{j+1}-q_j)\\
=&\langle q_j-q_{j-1},p_{j-1}\rangle + \langle q_{j+1}-q_j,p_j\rangle\\
=& \langle q_j-q_{j-1},p_{j-1}\rangle + \langle q_{j+1}-q_j,p_j\rangle + \langle q_j-q_j^*,p_j-p_{j-1}\rangle\\
=& \langle q_j^*-q_{j-1},p_{j-1}\rangle + \langle q_{j+1}-q_j^*,p_j\rangle\\
=& \langle q_j^*-q_{j-1},p_{j-1}^*\rangle + \langle q_j^*-q_{j-1},p_{j-1}-p_{j-1}^*\rangle \\
& +  \langle q_{j+1}-q_j^*,p_j^*\rangle + \langle q_{j+1}-q_j^*,p_j-p_j^*\rangle\\
=&\Sigma_j(q_j^*) + \langle q_j^*-q_{j-1},p_{j-1}-p_{j-1}^*\rangle + \langle q_{j+1}-q_j^*,p_j-p_j^*\rangle
\end{align*}}%
for all $q_j^*\in H_j$, where $p_j^*,p_{j-1}^*\in \partial T$ were chosen (possibly not uniquely) to fulfill
\beq q_j^*-q_{j-1}\in N_T(p_{j-1}^*)\;\text{ and }\; q_{j+1}-q_j^*\in N_T(p_j^*).\label{eq:System3}\eeq
By \eqref{eq:System3}, it follows from the definition of the normal cone and $p_j,p_{j-1}\in\partial T$ together with the convexity of $T$ that
\beqq \langle q_j^*-q_{j-1},p_{j-1}-p_{j-1}^*\rangle \leq 0\; \text{ and }\; \langle q_{j+1}-q_j^*,p_j-p_j^*\rangle \leq 0. \eeqq
This implies
\beqq \Sigma_j(q_j)\leq \Sigma_j(q_j^*)\eeqq
for all $q_j^*\in H_j$. We conclude that the polygonal curve segment $(q_{j-1},q_j,q_{j+1})$ fulfills the weak Minkowski billiard reflection rule in $q_j$ with respect to the $K$-supporting hyperplane $H_j$ through $q_j$ which is normal to $n_K(q_j)$. This eventually means that $(q_1,...,q_m)$ is a closed weak $(K,T)$-Minkowski billiard trajectory.

We proceed by proving that for strictly convex body $T\subset\R^n$, every closed weak Minkowski billiard trajectory is a strong one. So, let $q=(q_1,...,q_m)$ be a closed weak $(K,T)$-Minkowski billiard trajectory. We define a closed polygonal curve $p=(p_1,...,p_m)$ by
\beq q_{j+1}-q_j\in N_T(p_j).\label{eq:System2}\eeq
We note that due to the strict convexity of $T$ (cf.\;Lemma \ref{Lem:strictconvexitynormalcones}), $p$ is uniquely determined. Let $j\in\{1,...,m\}$ be arbitrarily chosen. Since there is a $K$-supporting hyperplane $H_j$ through $q_j$ such that $q_j$ minimizes $\Sigma_j(q_j^*)$ over all $q_j^*\in H_j$, we conclude by Lagrange's multiplier theorem that there is a $\mu_j\in\R$ with
\beqq \nabla_{q_j^*}\Sigma_j(q_j^*)_{\vert q_j^*=q_j}=\mu_j n_{H_j},\eeqq
where $n_{H_j}$ is the outer unit vector normal to $H_j$. We note that the differentiability of $\Sigma_j$ follows from the strict convexity of $T$--cf.\;the following calculation or, more basically, the duality between strict convexity of $T$ and smoothness of $T^\circ$ (cf.\;\cite[Theorem 11.13]{RockafellarWets2009}). We calculate the left side:
\allowdisplaybreaks{\begin{align*}
&\nabla_{q_j^*}\Sigma_j(q_j^*)_{\vert q_j^*=q_j}\\
\stackrel{(a)}{=}& \nabla_{q_j^*}\left(\langle q_j^*-q_{j-1},p_{j-1}(q_j^*)\rangle + \langle q_{j+1}-q_j^*,p_j(q_j^*)\rangle\right)_{\vert q_j^*=q_j}\\
\stackrel{(b)}{=}& \left(\partial_{q_j^*,i} \left[\langle q_j^*-q_{j-1},p_{j-1}(q_j^*)\rangle + \langle q_{j+1}-q_j^*,p_j(q_j^*)\rangle\right]_{\vert q_j^*=q_j}\right)_{i\in\{1,...,n\}}\\
\stackrel{(c)}{=}&\bigg(\lim_{\eps\rightarrow 0} \frac{1}{\eps}\big[\langle q_j+\eps e_i-q_{j-1},p_{j-1}(q_j+\eps e_i)\rangle - \langle q_j-q_{j-1},p_{j-1}\rangle\\
&\quad \quad \quad \quad \quad \quad+\langle q_{j+1}-q_j-\eps e_i,p_j(q_j+\eps e_i)\rangle - \langle q_{j+1}-q_j,p_j\rangle\big]\bigg)_{i\in\{1,...,n\}}\\
=&\bigg(\lim_{\eps\rightarrow 0}\frac{1}{\eps}\big[\langle q_j-q_{j-1},p_{j-1}(q_j+\eps e_i)-p_{j-1}\rangle + \langle q_{j+1}-q_j,p_j(q_j+\eps e_i)-p_j\rangle\\
& \quad \quad \quad \quad \quad \quad +\langle \eps e_i,p_{j-1}(q_j+\eps e_i)-p_j(q_j+\eps e_i)\rangle\big]\bigg)_{i\in\{1,...,n\}}\\
\stackrel{(d)}{=}& \left( \lim_{\eps\rightarrow 0}\langle e_i,p_{j-1}(q_j+\eps e_i)-p_j(q_j+\eps e_i)\rangle\right)_{i\in\{1,...,n\}}\\
\stackrel{(e)}{=}&\left(\langle e_i,p_{j-1}-p_j \rangle\right)_{i\in\{1,...,n\}}\\
=&p_{j-1}-p_j,
\end{align*}}%
where in equality ($a$), for every $q_j^*\in\R^n$ $p_{j-1}(q_j^*)$, $p_j(q_j^*)$ (here $p_{j-1}$ and $p_j$ are acting as functions) are the boundary points of $T$ (uniquely determined for $q_j^*\in\R^n$ near $q_j$ due to the strict convexity of $T$ and Lemma \ref{Lem:strictconvexitynormalcones}) fulfilling
\beq q_j^*-q_{j-1}\in N_T(p_{j-1}(q_j^*))\;\text{ and }\; q_{j+1}-q_j^*\in N_T(p_j(q_j^*)),\label{eq:System4}\eeq
where we note 
\beqq p_{j-1}(q_j)=p_{j-1}\; \text{ and } \;p_j(q_j)=p_j,\eeqq
in equality ($b$), by $\partial_{q_j^*,i}$ we denote the $i$-th partial derivative with respect to $q_j^*$, and in equality ($c$), by $e_i$ we denote the $i$-th standard unit vector in $\R^n$. In equality ($d$), we used
\beq \lim_{\eps\rightarrow 0}\frac{1}{\eps}\langle q_j-q_{j-1},p_{j-1}(q_j+\eps e_i)-p_{j-1}\rangle = 0\label{eq:refoften}\eeq
for all $i\in\{1,...,n\}$.

Indeed, if
\beq \dim N_T(p_{j-1})=1,\label{eq:case1}\eeq
then
\beqq \lim_{\eps\rightarrow 0} \frac{p_{j-1}(q_j+\eps e_i)-p_{j-1}}{\eps}\eeqq
is a tangent vector at $\partial T$ in $p_j$ and therefore
\beqq \langle q_j-q_{j-1},\lim_{\eps\rightarrow 0} \frac{p_{j-1}(q_j+\eps e_i)-p_{j-1}}{\eps}\rangle =0,\eeqq
and consequently
\begin{align*}
&\lim_{\eps\rightarrow 0} \frac{1}{\eps}\langle q_j-q_{j-1}, p_{j-1}(q_j+\eps e_i)-p_{j-1} \rangle \\
=& \langle q_j-q_{j-1},\lim_{\eps\rightarrow 0} \frac{p_{j-1}(q_j+\eps e_i)-p_{j-1}}{\eps}\rangle \\
=& 0.
\end{align*}
If
\beq \dim N_T(p_{j-1})>1\; \text{ and }\; q_j-q_{j-1}\in \mathring{N}_T(p_{j-1}),\label{eq:case2}\eeq
then it follows
\beqq p_{j-1}(q_j+\eps e_i)-p_{j-1}=p_{j-1}-p_{j-1}=0\eeqq
for $|\eps|$ small and therefore \eqref{eq:refoften}. If
\beqq \dim N_T(p_{j-1})>1\; \text{ and }\; q_j-q_{j-1}\in \partial N_T(p_{j-1}),\eeqq
then for $\eps >0$ \eqref{eq:refoften} follows from the argument either made for the case \eqref{eq:case1} or for the case \eqref{eq:case2}. Similarly, for $\eps <0$ \eqref{eq:refoften} follows from the argument either made for the case \eqref{eq:case1} or for the case \eqref{eq:case2}.

By similar reasoning, we derive
\beqq \lim_{\eps\rightarrow 0}\frac{1}{\eps}\langle q_{j+1}-q_j,p_j(q_j+\eps e_i)-p_j\rangle = 0\; \text{ for all }i\in\{1,...,n\}.\eeqq
In equality ($e$), we applied the continuity (which holds due to the strict convexity of $T$) of the functions $p_j$ and $p_{j-1}$ defined by \eqref{eq:System4}. Therefore, we conclude
\beq p_{j}-p_{j-1} = -\mu_j n_{H_j}.\label{eq:System5}\eeq
It remains to show $\mu_j \geq 0$. For that, scalar multiplication of \eqref{eq:System5} by $n_{H_j}$ implies
\beqq \langle p_j-p_{j-1},n_{H_j}\rangle = -\mu_j.\eeqq
From
\beqq \langle q_j-q_{j-1},n_{H_j}\rangle \geq 0\;\text{ and }\; \langle q_{j+1}-q_j,n_{H_j}\rangle \leq 0\eeqq
together with
\beqq q_j-q_{j-1}\in N_T(p_{j-1}) \; \text{ and } \; q_{j+1}-q_j\in N_T(p_j), \eeqq
it follows from the convexity of $T$ that
\beqq \langle p_j-p_{j-1},n_{H_j}\rangle \leq 0,\eeqq
and therefore $\mu_j \geq 0$.

\begin{figure}[h!]
\centering
\def\svgwidth{360pt}
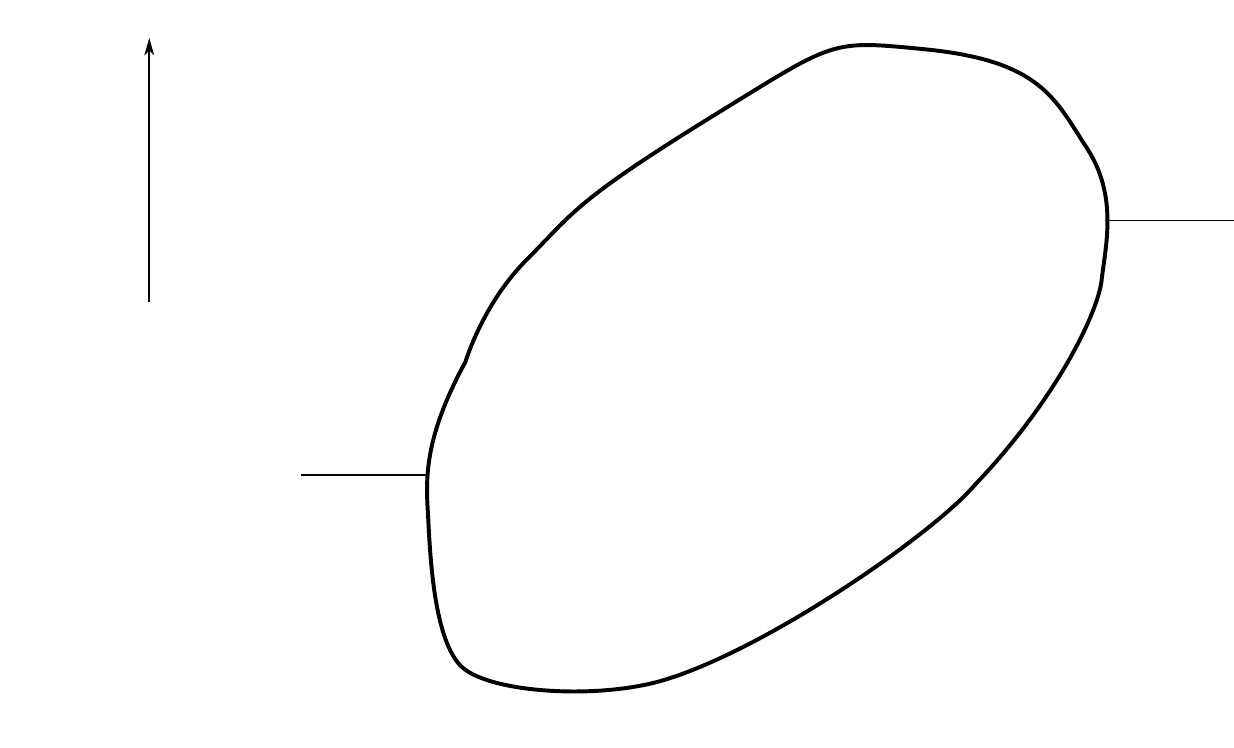
\caption[Illustration of an argument within the proof of Theorem \ref{Thm:weakstrongbilliards} I]{Illustration of $\partial T_{H_j^+}$ and $\partial T_{H_j^-}$. By $H_j^-$ and $H_j^+$ we denote the closed half-spaces of $\R^n$ bounded by $H_j$.}
\label{img:upperlowerboundary}
\end{figure}

Indeed,
\beqq \langle q_j-q_{j-1},n_{H_j}\rangle \geq 0 \; \text{ together with } \; q_j-q_{j-1}\in N_T(p_{j-1})\eeqq
implies that
\beq p_{j-1}\in \{p'\in \partial T : \langle n,n_{H_j} \rangle \geq 0 \; \forall n\in N_T(p')\}=:\partial T_{H_j^+},\label{eq:upperboundary}\eeq
and
\beqq \langle q_{j+1}-q_j,n_{H_j}\rangle \leq 0 \; \text{ together with } \; q_{j+1}-q_j\in N_T(p_j)\eeqq
implies
\beq p_j\in \{p'\in \partial T : \langle n,n_{H_j}\rangle \leq 0 \; \forall n\in N_T(p')\}=:\partial T_{H_j^-}.\label{eq:lowerboundary}\eeq
If $\mu_j <0$, i.e., $p_j-p_{j-1}$ is a positive multiple of $n_{H_j}$, then it follows from the strict convexity of $T$ (cf.\;Figure \ref{img:upperlowerboundary}) that
\beqq p_j\in \partial T_{H_j^+}\;\text{ and }\; p_{j-1}\in \partial T_{H_j^-},\eeqq
a contradiction to \eqref{eq:upperboundary} and \eqref{eq:lowerboundary}. Therefore it follows $\mu_j \geq 0 $.

\begin{figure}[h!]
\centering
\def\svgwidth{380pt}
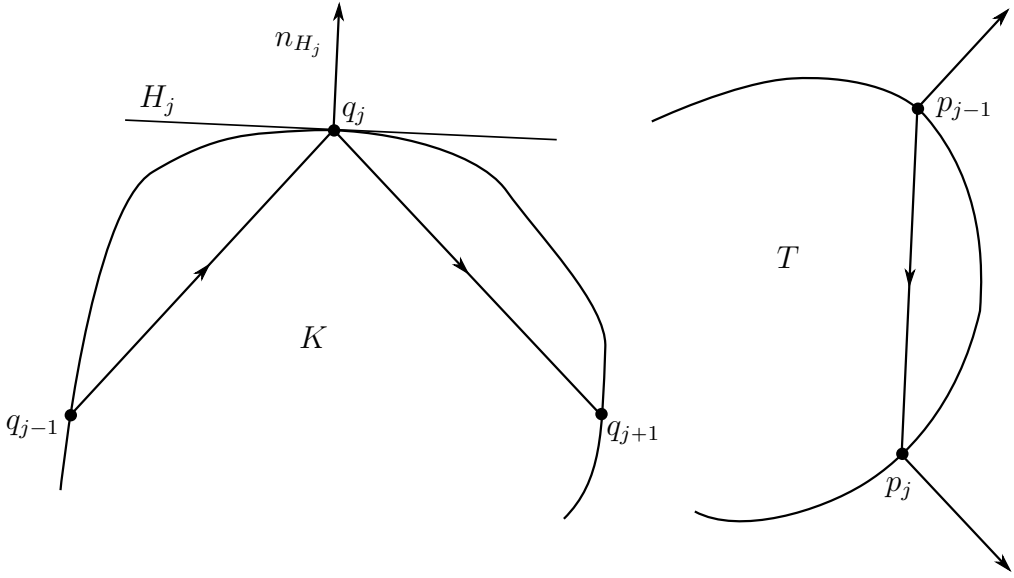
\caption[Illustration of an argument within the proof of Theorem \ref{Thm:weakstrongbilliards} II]{We have $q_j-q_{j-1}\in N_T(p_{j-1})$, $q_{j+1}-q_j\in N_T(p_j)$ and $p_j-p_{j-1}=-\mu_j n_{H_j}$, $\mu_j \geq 0$.}
\label{img:bouncingruleinqj}
\end{figure}

If we define
\beqq n_K(q_j):=n_{H_j}\quad  \forall j\in\{1,...,m\},\eeqq
then the pair $(q,p)$ fulfills \eqref{eq:System} by referring to \eqref{eq:System2} and \eqref{eq:System5}.
\epf


The following rather obvious proposition is needed for the proof of Theorem \ref{Thm:RegularityResult} and \ref{Thm:maximxallyspanning}. It follows immediately from within the proof of Theorem \ref{Thm:weakstrongbilliards} when $T$ is additionally required to be smooth.

\bprop\label{Prop:onlyonehyperplane}
Let $K,T\subset\R^n$ be convex bodies, where $T$ is additionally assumed to be strictly convex and smooth. Let $q=(q_1,...,q_m)$ be a closed $(K,T)$-Minkowski billiard trajectory. Then, for every $j\in\{1,...,m\}$, there is only one $K$-supporting hyperplane through $q_j$ for which the weak Minkowski billiard reflection rule in $q_j$ is satisfied.
\eprop

Before we prove this proposition, it is important to note that the combination of Proposition \ref{Prop:dualbilliardtrajunique} and Theorem \ref{Thm:weakstrongbilliards} does not imply this statement. That is because the vertices of the closed dual billiard trajectory can coincide, and as consequence, condition \eqref{eq:weakstrongcondition} may not determine a unique hyperplane (cf.\;Example C in Section \ref{Sec:Examples}). We will see that the smoothness of $T$ guarantees that the closed dual billiard trajectory is a closed polygonal curve (cf.\;Footnote \ref{foot:polygonalline}).

\bpf[Proof of Proposition \ref{Prop:onlyonehyperplane}]
Since $q$ is a closed $(K,T)$-Minkowski billiard trajectory, due to the proof of Theorem \ref{Thm:weakstrongbilliards}, there is a hyperplane $H_j$ through $q_j$ and an outer unit vector $n_{H_j}$ normal to $H_j$ such that
\beq \nabla_{q_j^*}\Sigma_j(q_j^*)|_{q_j^*=q_j}=p_{j-1}-p_j=\mu_j n_{H_j},\label{eq:hintergrund}\eeq
where $\mu_j \geq 0$ and $p_{j-1}$ and $p_j$ are uniquely determined by
\beqq q_{j+1}-q_j\in N_T(p_j)\; \text{ and }\; q_j-q_{j-1}\in N_T(p_{j-1}).\eeqq

We aim to show that $H_j$, respectively $n_{H_j}$, is uniquely determined. Against the background of \eqref{eq:hintergrund}, it would be enough to show that $\mu_j >0$.

If $\mu_j=0$, then $p_{j-1}=p_j$. But this implies that $q_{j+1}-q_j$ and $q_j-q_{j-1}$ are in the same normal cone, while they are not parallel (since $q_j$ is not on the line connecting $q_{j-1}$ and $q_{j+1}$). Therefore, there are two linearly independent nonzero vectors in $N_T(p_j)=N_T(p_{j-1})$, a contradiction to the smoothness of $T$.

Consequently, it follows $\mu_j >0$ and therefore the uniquess of $n_{H_j}$ and $H_j$.
\epf

The statement of Proposition \ref{Prop:onlyonehyperplane} is not true without requiring $T$ to be both strictly convex and smooth. For the necessity of the strict convexity, we refer to Example B, for the necessity of the smoothness, we refer to Example C (cf.\;Section \ref{Sec:Examples}).

Its notion suggests that dual billiard trajectories are in fact billiard trajectories. Indeed, the following Propositions \ref{Prop:lengthdualbilliard} and \ref{Prop:dualbilliard} show that for strictly convex and smooth body $T\subset\R^n$ the closed dual billiard trajectory of a closed $(K,T)$-Minkowski billiard trajectory $q$ is a closed $(T,-K)$-Minkowski billiard trajectory $p$ with
\beqq \ell_{-K}(p)=\ell_T(q).\eeqq

For the case $T$ is a strictly convex body in $\R^n$, this also implies that the $\ell_T$-length of $\ell_T$-minimizing closed $(K,T)$-Minkowski billiard trajectories equals the $\ell_{-K}$-length of $\ell_{-K}$-minimizing closed $(T,-K)$-Minkowski billiard trajectories\footnote{The existence of these minimums will be proved in Theorem \ref{Thm:onetoone} under the condition of strict convexity of $T$ and in \cite[Theorem 2.2]{Rudolf2022} for the general case.}. Later, we will use this fact for the proof of Theorem \ref{Thm:maximxallyspanning}.

\bprop\label{Prop:lengthdualbilliard}
Let $K,T\subset\R^n$ be convex bodies. Let $q=(q_1,...,q_m)$ be a closed $(K,T)$-Minkowski billiard trajectory with closed dual billiard trajectory $p=(p_1,...,p_m)$ in $T$. Then, we have
\beqq \ell_T(q)=\ell_{-K}(p).\eeqq
\eprop

\bpf
By definition of the Minkowski billiard reflection rule, we have
\beqq q_{j+1}-q_j\in N_T(p_j)\;\text{ and }\; p_{j+1}-p_j\in -N_K(q_{j+1}).\eeqq
for all $j\in\{1,...,m\}$. By recalling Proposition \ref{Prop:Minkowskisupport}, this implies
\beqq \mu_{T^\circ}(q_{j+1}-q_j)=h_T(q_{j+1}-q_j)=\langle q_{j+1}-q_j,p_j\rangle\eeqq
and
\beqq \mu_{K^\circ}(p_j-p_{j+1})=h_K(p_j-p_{j+1})=\langle p_j-p_{j+1},q_{j+1}\rangle\eeqq
for all $j\in\{1,...,m\}$. Then, we compute
\allowdisplaybreaks{\begin{align*}
\ell_{T}(q)=&\sum_{j=1}^m \mu_{T^\circ}(q_{j+1}-q_{j}) =\sum_{j=1}^m \langle q_{j+1}-q_{j},p_{j}\rangle \\
=&(\langle q_2,p_1\rangle - \langle q_1,p_1\rangle) + (\langle q_3,p_2\rangle - \langle q_2,p_2\rangle)+...\\
& ... + (\langle q_{m},p_{m-1}\rangle -\langle q_{m-1},p_{m-1} \rangle) + (\langle q_{m+1},p_{m}\rangle -\langle q_{m},p_{m} \rangle)\\
=& -\langle q_1,p_1\rangle +( \langle q_2,p_1\rangle -\langle q_2,p_2\rangle) +( \langle q_3,p_2\rangle - \langle q_3,p_3 \rangle) + ... \\
&... +(\langle q_{m-1},p_{m-2}\rangle - \langle q_{m-1},p_{m-1}\rangle) +(\langle q_{m},p_{m-1} \rangle - \langle q_{m},p_{m}\rangle ) +\langle q_{m+1},p_{m} \rangle\\
=& \sum_{j=1}^m \langle p_{j}-p_{j+1},q_{j+1}\rangle =\sum_{j=1}^m \mu_{K^\circ}(p_j-p_{j+1}) = \sum_{j=1}^m \mu_{-K^\circ }(p_{j+1}-p_{j})\\
=& \sum_{j=1}^m \mu_{(-K)^\circ}(p_{j+1}-p_j)\\
=&\ell_{-K}(p),
\end{align*}}%
where we used $q_{m+1}=q_1$, $p_{m+1}=p_1$ and the property
\beqq \mu_{K^\circ}(-\cdot)=\mu_{-K^\circ}(\cdot)\eeqq
of the Minkowski functional utilizing Proposition \ref{Prop:IntroMinkowski}(iii) and
\beqq (-K)^\circ = -K^\circ\eeqq
(cf.\;Proposition \ref{Prop:polarbody}).
\epf

\bprop\label{Prop:dualbilliard}
Let $K,T\subset\R^n$ be convex bodies and $T$ is additionally assumed to be strictly convex and smooth. Let $q=(q_1,...,q_m)$ be a closed $(K,T)$-Minkowski billiard trajectory with its closed dual billiard trajectory $p=(p_1,...,p_m)$ in $T$. Then, $p$ is a closed $(T,-K)$-Minkowski billiard trajectory with
\beqq -q^{+1}:=(-q_2,...,-q_m,-q_1)\eeqq
as closed dual billiard trajectory on $-K$.
\eprop

\bpf
Since the pair $(q,p)$ fulfills \eqref{eq:System}, $p_1,...,p_m$ are uniquely determined (cf.\;Lemma \ref{Lem:strictconvexitynormalcones}) by the condition
\beq q_{j+1}-q_j\in N_T(p_j)\quad \forall j\in\{1,...,m\}.\label{eq:dualbilliard1}\eeq
Since $q_1,...,q_m$ satisfy $q_j\neq q_{j+1}$ for all $j\in\{1,...,m\}$ and $q_j$ is not contained in the line segment connecting $q_{j-1}$ and $q_{j+1}$ for all $j\in\{1,...,m\}$ (cf.\;Footnote \ref{foot:polygonalline}),
\beq q_2-q_1,...,q_m-q_{m-1},q_1-q_m\label{eq:dualbilliard2}\eeq
are nonzero and satisfy
\beqq q_j-q_{j-1} \nparallel q_{j+1}-q_j\quad \forall j\in\{1,...,m\}.\eeqq
Then, \eqref{eq:dualbilliard1} together with the strict convexity and smoothness of $T$ implies that $p_1,...,p_m$ satisfy $p_j\neq p_{j+1}$ for all $j\in\{1,...,m\}$ and $p_j$ is not contained in the line segment connecting $p_{j-1}$ and $p_{j+1}$ for all $j\in\{1,...,m\}$.

This implies
\beq \begin{cases} p_{j+1}-p_j=-\mu_{j+1}n_K(q_{j+1})\in -N_K(q_{j+1})=N_{-K}(-q_{j+1})=N_{-K}(-q_j^{+1})\\(-q_{j+1}^{+1})-(-q_j^{+1})=(-q_{j+2})-(-q_{j+1})\in - N_T(p_{j+1})\end{cases}\label{eq:dualbilliard3}\eeq
for all $j\in\{1,...,m\}$, where we used
\beqq N_{-K}(-q_{j+1})=-N_K(q_{j+1})\quad \forall j\in\{1,...,m\}.\eeqq
From \eqref{eq:dualbilliard3}, we conclude that the pair $(p,-q^{+1})$ fulfills \eqref{eq:System} for the configuration $(T,-K)$. Therefore, $p$ is a closed $(T,-K)$-Minkowski billiard trajectory with $-q^{+1}$ as its closed dual billiard trajectory on $-K$.
\epf

In order to show the necessity of the smoothness of $T$ in Proposition \ref{Prop:dualbilliard}, we refer to Example C (cf.\;Section \ref{Sec:Examples})--there, one can construct a closed polygonal curve $p$, but which does not satisfy $p_j\neq p_{j+1}$ for all $j$.




We have the following proposition as generalization of \cite[Proposition 2.3]{KruppRudolf2020} to the Minkowski/Finsler setting:

\bprop\label{Prop:SectionInvariance}
Let $K,T\subset\R^n$ be convex bodies. Let $q=(q_1,...,q_m)$ be a closed weak $(K,T)$-Minkowski billiard trajectory and $V\subseteq\R^n$ an affine subspace such that $K\cap V$ is an affine section of $K$ containing $q$. Then, $q$ is a closed weak $(K\cap V,T)$-Minkowski billiard trajectory.\footnote{We notice that the dimension of $K\cap V$ is possibly smaller than the dimension of $T$. For these cases, we consider that Definiton \ref{Def:weakMinkowskiBilliards} can be easily extended to covex bodies $T$ which are allowed to have dimension greater $n$.}
\eprop

\bpf
Since $q$ is a closed weak $(K,T)$-Minkowski billiard trajectory, there are $K$-supporting hyperplanes $H_1,...,H_m$ through $q_1,...,q_m$ such that $q_j$ minimizes
\beq \mu_{T^\circ}(\widebar{q}_j-q_{j-1})+\mu_{T^\circ}(q_{j+1}-\widebar{q}_j)\label{eq:minimizationprop}\eeq
over all $\widebar{q}_j\in H_j$ for all $j\in\{1,...,m\}$. Since $K\cap V$ contains $q$ it follows that $q_j$ minimizes \eqref{eq:minimizationprop} over all $\widebar{q}_j\in H_j\cap V$ for all $j\in\{1,...,m\}$. This implies that $q$ is a closed weak $(K\cap V,T)$-Minkowski billiard trajectory.
\epf

Clearly, in general, the converse is not true: We can imagine an affine section $K\cap V$ of $K$ that can be translated into $\mathring{K}$. Then, every closed weak $(K\cap V,T)$-Minkowski billiard trajectory $q$ can be translated into $\mathring{K}$. But in Proposition \ref{Prop:notranslation}, we will prove that $q\in F(K)$, a contradiction.

In \cite[Examples A, B, C, and D]{KruppRudolf2020}, we have seen (for $T$ equals the Euclidean unit ball) that in general the length minimality of a closed weak $(K,T)$-Minkowski billiard trajectory is not invariant under going to (inclusion minimal) affine sections of $K$ containing the billiard trajectory. This billiard trajectory may not even locally minimize the length of closed polygonal curves in $F(K\cap V)$. We note that these examples can be easily generalized to settings when the weak Minkowski billiard reflection rule is not necessarily governed by the Euclidean unit ball.

The next two propositions make a statement concerning the positional relationship of the hyperplanes which determine the weak Minkowski billiard reflection rule.

\bprop\label{Prop:normalvectorsspanning}
Let $K,T\subset\R^n$ be convex bodies. Let $q=(q_1,...,q_m)$ be a closed $(K,T)$-Minkowski billiard trajectory with respect to $H_1,...,H_m$. Then, we have
\beqq 0\in \conv\{n_K(q_1),...,n_K(q_m)\},\eeqq
where $n_K(q_1),...,n_K(q_m)$ are the outer unit vectors normal to $H_1,...,H_m$.

If $T$ is assumed to be smooth and if we denote by $\widetilde{U}$ the convex cone spanned by 
\beqq n_K(q_1),...,n_K(q_m),\eeqq
then $\widetilde{U}$ is a linear subspace of $\R^n$ with dimension less or equal than $m-1$.
\eprop

\bpf
Let $p=(p_1,...,p_m)$ be a closed dual billiard trajectory of $q$. Then, there are $\mu_1,...,\mu_{m}\geq 0$ with
\beq p_{j+1}-p_j=-\mu_{j+1}n_K(q_{j+1})\quad \forall j\in\{1,...,m\}.\label{eq:normalvectorspanning2}\eeq
We first consider the case
when $\mu_j>0$ for all $j\in\{1,...,m\}$. Then, we define
\beqq s_j:=\frac{\mu_{j+1}}{\mu_1+...+\mu_m} \quad \forall j\in\{1,...,m\}\eeqq
and conclude
\beq \sum_{j=1}^m s_j\frac{1}{\mu_{j+1}}(p_{j+1}-p_j)=\frac{1}{\mu_1+...+\mu_m}\sum_{j=1}^m (p_{j+1}-p_j)=0 \label{eq:normalvectorspanning1}\eeq
while
\beq \sum_{j=1}^m s_j = \sum_{j=1}^m \frac{\mu_{j+1}}{\mu_1+...+\mu_m}=1.\label{eq:normalvectorspanning11}\eeq
This implies by the definition of the convex hull that
\begin{align}
0&\in \left\{\sum_{j=1}^m \widetilde{s}_j\frac{1}{\mu_{j+1}} (p_{j+1}-p_j):\sum_{j=1}^m \widetilde{s}_j =1, \;\widetilde{s}_j \geq 0\right\}\notag\\
&=\conv\left\{\frac{1}{\mu_2}(p_2-p_1),...,\frac{1}{\mu_m}(p_m-p_{m-1}),\frac{1}{\mu_1}(p_1-p_m)\right\}\label{eq:normalvectorspanning3}
\end{align}
and therefore, together with \eqref{eq:normalvectorspanning2},
\beqq 0\in\conv\{-n_K(q_1),...,-n_K(q_m)\}\eeqq
and consequently
\beq 0\in\conv\{n_K(q_1),...,n_K(q_m)\}.\label{eq:normalvectorspanning4}\eeq
If $\mu_j=0$ for some $j\in\{1,...,m\}$ (for all is impossible), then, by \eqref{eq:normalvectorspanning2}, also
\beqq p_{j}-p_{j-1}=0\eeqq
for all these $j$. But then, the vector corresponding to $p_{j}-p_{j-1}$ can be removed from within the set of vectors building the convex hull in \eqref{eq:normalvectorspanning3} without influencing \eqref{eq:normalvectorspanning1} and \eqref{eq:normalvectorspanning11}. Therefore, $0$ is in the convex hull of the nonzero
\beqq \frac{1}{\mu_{j+1}}(p_{j+1}-p_j),\eeqq
what implies that $0$ is in the convex hull of the associated unit normal vectors $n_K(q_{j+1})$. But the latter convex hull is subset of
\beqq \conv\{n_K(q_1),...,n_K(q_m)\}.\eeqq
Therefore, we derive \eqref{eq:normalvectorspanning4}.

Let us assume $T$ is smooth. Then, this implies
\beq p_{j+1}-p_j\neq 0\quad \forall j\in\{1,...,m\}\label{eq:normalvectorspanning5}\eeq
(cf.\;Proposition \ref{Prop:onlyonehyperplane}). It is
\beqq (p_2-p_1) + ... + (p_m-p_{m-1}) + (p_1-p_m) = 0.\eeqq
For
\beqq s_j:=\frac{1}{m} \quad \forall j\in\{1,...,m\}, \eeqq
we also have
\beqq s_1(p_2-p_1) + ... + s_{m-1} (p_m-p_{m-1}) + s_m (p_1-p_m) = 0.\eeqq
Since
\beqq \sum_{j=1}^m s_j =1 \; \text{ and } \;  s_j\neq 0 \; \; \forall j\in\{1,...,m\},\eeqq
it follows that $0$ lies in the relative interior of
\beqq \conv\{p_2-p_1,...,p_m-p_{m-1}, p_1-p_m\}.\eeqq
But this implies that the convex cone spanned by
\beq p_2-p_1,...,p_m-p_{m-1},p_1-p_{m}\label{eq:normalvectorspanning6}\eeq
and therefore, by \eqref{eq:normalvectorspanning2} and \eqref{eq:normalvectorspanning5}, also $\widetilde{U}$ is a linear subspace of $\R^n$. Obviously, then $\widetilde{U}$ is the inclusion minimal linear subspace containing the vectors \eqref{eq:normalvectorspanning6} and consequently has dimension less or equal than $m-1$.
\epf

The necessity of the smoothness of $T$ for the second statement follows by referring to Example C. In this case, the convex cone spanned by $n_K(q_1),n_K(q_2),n_K(q_3)$ is not a linear subspace of $\R^2$. Furthermore, for the weaker situation of closed weak $(K,T)$-Minkowski billiard trajectory, one can show--having in mind Theorem \ref{Thm:weakstrongbilliards}--the necessity of the strict convexity of $T$ for the second statement by referring to Example D--and also for the first statement by referring to Example E (cf.\;Section \ref{Sec:Examples} for the three examples).

\bprop\label{Prop:noteinproof}
Let $K,T\subset\R^n$ be convex bodies. Let $q=(q_1,...,q_m)$ be a closed $(K,T)$-Minkowski billiard trajectory with respect to $H_1,...,H_m$ and let $U$ be the inclusion minimal linear subspace of $\R^n$ containing the outer unit vectors $n_K(q_1),...,n_K(q_m)$ which are normal to $H_1,...,H_m$. We denote by $H_1^+,...,H_m^+$ the closed half-spaces of $\R^n$ which are bounded by $H_1,...,H_m$ and contain $K$. Further, let $W$ be the orthogonal complement to $U$ in $\R^n$. Then, we can write
\beqq H_j=(H_j\cap U)\oplus W\; \text{ and }\; H_j^+=(H_j^+\cap U)\oplus W\eeqq
for all $j\in\{1,...,m\}$ and have that
\beqq \bigcap_{j=1}^m \left(H_j^+\cap U\right)\text{ is nearly bounded in }U,\quad \bigcap_{j=1}^m H_j^+\text{ is nearly bounded in }\R^n.\eeqq
If $T$ is assumed to be smooth, then $U$ coincides with the convex cone spanned by $n_K(q_1),...,n_K(q_m)$ and we have that
\beqq \bigcap_{j=1}^m \left(H_j^+\cap U\right)\text{ is bounded in }U,\quad \bigcap_{j=1}^m H_j^+\text{ is nearly bounded in }\R^n.\eeqq
\eprop

\bpf
Since $U$ is a linear subspace of $\R^n$ containing $n_K(q_1),...,n_K(q_m)$, we can write
\beq H_j=(H_j\cap U)\oplus W\;\text{ and }\; H_j^+=\left(H_j^+\cap U\right)\oplus W.\label{eq:noteinproof0}\eeq
Let $\widetilde{U}$ be the convex cone spanned by $n_K(q_1),...,n_K(q_m)$. By Proposition \ref{Prop:normalvectorsspanning}, we have that
\beq 0\in\conv\{n_K(q_1),...,n_K(q_m)\}.\label{eq:noteinproof00}\eeq
Now, $0$ either is an interior point or a boundary point (both with respect to $U$) of the convex hull in \eqref{eq:noteinproof00}.

If $0$ is a boundary point (with respect to $U$) of the convex hull in \eqref{eq:noteinproof00}, then $\widetilde{U}$ is subset of a $\widetilde{U}$-supporting closed half-space $H_{\widetilde{U}}$ of $U$ while
\beqq \partial H_{\widetilde{U}}\cap \widetilde{U}\eeqq
contains a selection of unit vectors
\beqq \{n_K(q_{i_1}),...,n_K(q_{i_k})\}\subseteq \{n_K(q_1),...,n_K(q_m)\},\quad k\leq m,\eeqq
with $0$ in the relative interior of
\beqq \conv\{n_K(q_{i_1}),...,n_K(q_{i_k})\}\subseteq \partial H_{\widetilde{U}} \cap \widetilde{U}.\eeqq
We denote by $H_{i_1},...,H_{i_k}$ the associated $K$-supporting hyperplanes through $q_{i_1},...,q_{i_k}$ which are normal to $n_K(q_{i_1}),...,n_K(q_{i_k})$. It follows that
\beq \bigcap_{l=1}^k \left(H_{i_l}^+\cap U\right)\label{eq:ibtw1}\eeq
intersected with the convex cone spanned by the vectors $n_K(q_{i_1}),...,n_K(q_{i_{k}})$ is bounded in $\partial H_{\widetilde{U}}$. We denote this intersection by $I$. Then, we can write
\beqq \bigcap_{l=1}^k \left(H_{i_l}^+\cap U\right) = I \oplus I^{\perp_U},\eeqq
where by $I^{\perp_U}$ we denote the orthogonal complement in $U$ to the inclusion minimal linear subspace of $U$ containing $I$. Clearly, the boundedness of $I$ in $\partial H_{\widetilde{U}}$ implies the nearly boundedness of $I\oplus I^{\perp_U}$ in $U$: Because of the boundedness of $I$ in $\partial H_{\widetilde{U}}$, there are two parallel hyperplanes $G$ and $G+c$, $c\in \partial H_{\widetilde{U}}$, in $\partial H_{\widetilde{U}}$ such that $I$ lies in-between. Then, $I\oplus I^{\perp_U}$ lies between the two hyperplanes
\beq G \oplus I^{\perp_U} \; \text{ and } \; (G + c) \oplus I^{\perp_U}\; \text{ in }\; U,\label{eq:inbetween0}\eeq
i.e., it is nearly bounded in $U$. Then, using \eqref{eq:noteinproof0},
\beqq \bigcap_{l=1}^k H_{i_l}^+= \bigcap_{l=1}^k \left(\left(H_{i_l}^+\cap U\right)\oplus W\right)=\left( \bigcap_{l=1}^k \left(H_{i_l}^+\cap U\right)\right)\oplus W\eeqq
lies between the two parallel hyperplanes
\beq \left(G \oplus I^{\perp_U}\right)\oplus W \; \text{ and } \; \left((G+c) \oplus I^{\perp_U}\right)\oplus W \; \text{ in }\; \R^n,\label{eq:inbetween}\eeq
i.e., it is nearly bounded in $\R^n$. This implies that
\beqq \bigcap_{j=1}^m \left(H_j^+\cap U\right) \subseteq \bigcap_{l=1}^k \left(H_{i_l}^+ \cap U\right) \eeqq
lies between the two hyperplanes in \eqref{eq:inbetween0} and 
\begin{align*}
\bigcap_{j=1}^m H_j^+ &=\bigcap_{j=1}^m \left(\left(H_j^+\cap U\right)\oplus W\right)\\
& = \left(\bigcap_{j=1}^m \left(H_j^+\cap U\right)\right)\oplus W \\
&  \subseteq \bigcap_{l=1}^k H_{i_l}^+ \\
&= \left(\bigcap_{l=1}^k \left(H_{i_l}^+ \cap U\right)\right) \oplus W
\end{align*}
between the two hyperplanes in \eqref{eq:inbetween}, i.e., they are nearly bounded in $U$ and $\R^n$, respectively.

If $0$ is an interior point of the convex hull in \eqref{eq:noteinproof00}, i.e., when $\widetilde{U}$ coincides with $U$--and by Proposition \ref{Prop:normalvectorsspanning} this is also the case when $T$ is assumed to be smooth--, then this directly implies that
\beq \bigcap_{j=1}^m \left(H_j^+\cap U\right)\label{eq:noteinproof1}\eeq
is bounded in $U$. From this, we conclude that there are parallel hyperplanes $H$ and $H+d$, $d\in U$, in $U$ such that \eqref{eq:noteinproof1} lies in-between. With \eqref{eq:noteinproof0}, this implies that
\beqq \bigcap_{j=1}^m H_j^+=\bigcap_{j=1}^m \left(\left(H_j^+\cap U\right)\oplus W\right) = \left(\bigcap_{j=1}^m \left(H_j^+\cap U\right)\right)\oplus W\eeqq
lies between the parallel hyperplanes
\beqq H\oplus W\; \text{ and } \; (H+d)\oplus W \; \text{ in }\; \R^n\eeqq
and therefore is nearly bounded in $\R^n$.
\epf

\bprop\label{Prop:notranslation}
Let $K,T\subset\R^n$ be convex bodies. Let $q=(q_1,...,q_m)$ be a closed $(K,T)$-Minkowski billiard trajectory with closed dual billiard trajectory $p=(p_1,...,p_m)$. Then, we have
\beqq q\in F(K)\; \text{ and } \; p\in F(T).\eeqq
\eprop

\bpf
Let $H_1,...,H_m$ be the $K$-supporting hyperplanes through $q_1,...,q_m$ which are associated to the Minkowski billiard reflection rule and let $H_1^+,...,H_m^+$ the closed half-spaces of $\R^n$ containing $K$ and which are bounded by $H_1,...,H_m$. By Proposition \ref{Prop:noteinproof}, we conclude that
\beqq H_1^+\cap ... \cap H_m^+\eeqq
is nearly bounded in $\R^n$. Then, from \cite[Lemma 2.1(ii)]{KruppRudolf2020} it follows that
\beqq \{q_1,...,q_m\}\in F(K),\eeqq
i.e., $q\in F(K)$.

By the definition of the Minkowski billiard reflection rule, there are factors $\lambda_1,...,\lambda_m >0$ (which are $>0$ due to Footnote \ref{foot:polygonalline}) and unit vectors
\beqq n_T(p_1),...,n_T(p_m)\;\text{ in }\; N_T(p_1),...,N_T(p_m)\eeqq
such that
\beqq q_{j+1}-q_j=\lambda_j n_T(p_j)\quad \forall j\in\{1,...,m\}.\eeqq
Since $q$ is closed, we justify
\beq 0\in \conv\{n_T(p_1),...,n_T(p_m)\}\label{eq:notranslation1}\eeq
in a similar way to the proof of Proposition \ref{Prop:normalvectorsspanning}. Let $U'$ be the inclusion minimal linear subspace of $\R^n$ containing $n_T(p_1),...,n_T(p_m)$. Then, as in the proof of Proposition \ref{Prop:noteinproof}, \eqref{eq:notranslation1} implies that
\beq H_1'^+\cap ... \cap H_m'^+\label{eq:notranslation2}\eeq
is nearly bounded in $\R^n$, where $H_1'^+,...,H_m'^+$ are the closed half-spaces of $\R^n$ containing $T$ and which are bounded by $H_1',...,H_m'$ which are the $T$-supporting hyperplanes of $\R^n$ through $p_1,...,p_m$ normal to $n_T(p_1),...,n_T(p_m)$. By \cite[Lemma 2.1(ii)]{KruppRudolf2020}, it follows from the nearly boundedness of \eqref{eq:notranslation2} that
\beqq \{p_1,...,p_m\}\in F(T),\eeqq
i.e., $p\in F(T)$.
\epf

The first statement of Proposition \ref{Prop:notranslation}, i.e., $q\in F(K)$, in general, is not true when $q$ is just assumed to be a closed weak $(K,T)$-Minkowski billiard trajectory and $T$ is not required to be strictly convex. To see this, we consider Example E (cf.\;Section \ref{Sec:Examples}).

\bprop\label{Prop:GenBezdek}
Let $K,T\subset\R^n$ be convex bodies. Let $q=(q_1,...,q_m)$ be a closed $(K,T)$-Minkowski billiard trajectory with respect to $H_1,...,H_m$ and let $U$ be the inclusion minimal linear subspace of $\R^n$ containing the outer unit vectors $n_K(q_1),...,n_K(q_m)$ which are normal to $H_1,...,H_m$. Then, there is a selection
\beq \{i_1,...,i_{\dim U +1}\}\subseteq \{1,...,m\}\label{eq:selection}\eeq
such that
\beqq \{q_{i_1},...,q_{i_{\dim U +1}}\}\in F(K).\eeqq
\eprop

\bpf
For $m=\dim U+1$, we can just apply Proposition \ref{Prop:notranslation} and nothing more is to prove. If $\dim U=n$, i.e., $U=\R^n$, then the claim follows immediately from Proposition \ref{Prop:notranslation} and \cite[Lemma 2.1(i)]{KruppRudolf2020} (cf.\;the equivalent expression below this Lemma).

\begin{figure}[h!]
\centering
\def\svgwidth{320pt}
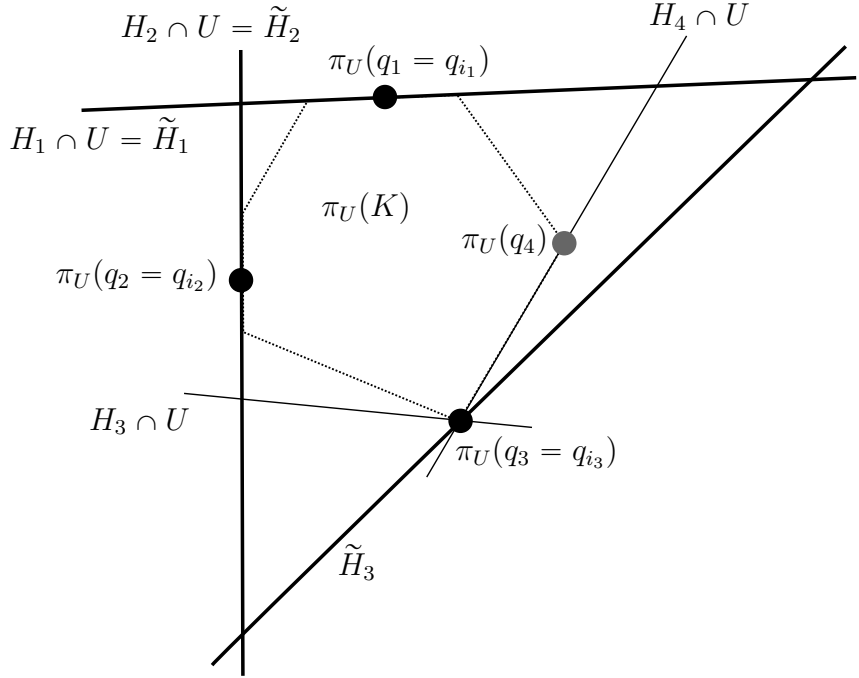
\caption[The choice of hyperplanes $\widetilde{H}_j$]{Illustration of the selection of $\{\pi_U(q_{i_1}),\pi_U(q_{i_2}),\pi_U(q_{i_3})\}$ out of $\{\pi_U(q_1),\pi_U(q_2),\pi_U(q_3),\pi_U(q_4)\}$ and the choice of $\bigcap_{j=1}^m\left(H_j^+\cap U\right)$-supporting hyperplanes $\widetilde{H}_1,\widetilde{H}_2,\widetilde{H}_3$ in $U$ such that $\bigcap_{j=1}^3\widetilde{H}^+_j$ is nearly bounded in $U$.}
\label{img:hyperplanechoice}
\end{figure}

Let
\beqq \dim U\leq \min\{n-1,m-2\}.\eeqq
Proposition \ref{Prop:noteinproof} implies, on the one hand, that we can write
\beqq H_j=(H_j\cap U)\oplus W\;\text{ and }\; H_j^+=(H_j^+\cap U)\oplus W\eeqq
for all $j\in\{1,...,m\}$, where $W$ is the orthogonal complement to $U$ in $\R^n$ and $H_1^+,...,H_m^+$ are the closed half-spaces of $\R^n$ containing $K$ and which are bounded by $H_1,...,H_m$, and, on the other hand, that
\beqq \bigcap_{j=1}^m \left(H_j^+\cap U\right)\eeqq
is nearly bounded in $U$. This implies by\cite[Lemma 2.1(ii)]{KruppRudolf2020} that
\beqq \pi_U(q)\in F\left(\bigcap_{j=1}^m \left(H_j^+ \cap U\right)\right),\eeqq
where we denote by $\pi_U$ the orthogonal projection onto $U$. Then, by \cite[Lemma 2.1(i)]{KruppRudolf2020}, there is a selection
\beqq \{i_1,...,i_{\dim U+1}\}\subset \{1,...,m\}\eeqq
such that
\beqq \left\{\pi_U(q_{i_1}),...,\pi_U(q_{i_{\dim U +1}})\right\} \in F\left(\bigcap_{j=1}^m\left(H_j^+ \cap U\right)\right).\eeqq
Referring again to \cite[Lemma 2.1(ii)]{KruppRudolf2020}, there are $\bigcap_{j=1}^m\left(H_j^+ \cap U\right)$-supporting hyperplanes\footnote{Not necessarily with $\widetilde{H}_j=H_{i_j}\cap U$.}
\beqq \widetilde{H}_1,...,\widetilde{H}_{\dim U +1}\eeqq
in $U$ through
\beqq \pi_U(q_{i_1}),...,\pi_U(q_{i_{\dim U +1}})\eeqq
such that
\beqq \bigcap_{j=1}^{\dim U +1}\widetilde{H}_j^+\eeqq
is nearly bounded in $U$ with
\beqq \bigcap_{j=1}^m \left(H_j^+\cap U\right) \subseteq \bigcap_{j=1}^{\dim U +1}\widetilde{H}_j^+,\eeqq
where $\widetilde{H}_j^+$ is the half-space bounded by $\widetilde{H}_j$ containing $\pi_U(K)$ for all $j\in\{1,...,\dim U +1\}$ (cf.\;Figure \ref{img:hyperplanechoice}).
Then, this implies that
\beqq \bigcap_{j=1}^{\dim U +1}\left(\widetilde{H}^+_j \oplus W\right)\eeqq
is nearly bounded in $\R^n$ with
\beqq K\subseteq \bigcap_{j=1}^m H_j^+ \subseteq \bigcap_{j=1}^{\dim U +1}\left(\widetilde{H}^+_j \oplus W\right)\eeqq
and
\beqq q_{i_j}\in \widetilde{H}_j \oplus W\quad \forall j\in\{1,...,\dim U +1\}.\eeqq
By using \cite[Lemma 2.1(ii)]{KruppRudolf2020}, this yields
\beqq \{q_{i_1},...,q_{i_{\dim U +1}}\} \in F(K).\eeqq
\epf

We remark that the statement of Proposition \ref{Prop:GenBezdek}, in general, is not true when requiring $q$ just to be a closed polygonal curve in $F(K)$ (and not a closed $(K,T)$-Minkowski billiard trajectory). In \cite{KruppRudolf2020}, we gave a counterexample for $T$ equals the Euclidean unit ball in $\R^n$.

In order to state/prove the upcoming Theorem \ref{Thm:onetoone}, we recall that for a convex body $K\subset\R^n$, the set $F_{n+1}^{cp}(K)$ is the set of all closed polygonal curves
\beqq q=(q_1,...,q_m)\in F(K)\eeqq
with $m\leq n+1$.

Let $(M,d)$ be a metric space and $P(M)$ the set of all nonempty compact subsets of $M$. We recall that $P(M)$ is a metric space together with the \textit{Hausdorff metric} $d_H$ which for $X,Y\in P(M)$ is defined by
\beqq d_H(X,Y)=\max\left\{\sup_{x\in X}\inf_{y\in Y} d(x,y),\sup_{y\in Y} \inf_{x\in X} d(x,y)\right\}.\eeqq

For the following Proposition \ref{Prop:F&l}, we denote by $\mathcal{C}(\R^n)$ the set of all convex bodies and by $ cp(\R^n)$ the set of all closed polygonal curves in $\R^n$. For $K\in\mathcal{C}(\R^n)$, we will consider $(F(K),d_H)$ and $(F_{n+1}^{cp}(K),d_H)$ as well as $( cp(\R^n),d_H)$ as metric subspaces of the complete metric space $(P(\R^n),d_H)$ which is induced by the Euclidean space $(\R^n,|\cdot|)$.

\bprop\label{Prop:F&l}
\begin{itemize}
\item[(i)] If $K,L\in\mathcal{C}(\R^n)$ with $K\subseteq L$, then
\beqq F(L)\subseteq F(K) \; \text{ and } \; F_{n+1}^{cp}(L) \subseteq F_{n+1}^{cp}(K).\eeqq
\item[(ii)] If $K\in\mathcal{C}(\R^n)$ and $c>0$, then
\beqq F(cK)=cF(K)  \; \text{ and } \; F_{n+1}^{cp}(cK)=cF_{n+1}^{cp}(K).\eeqq
\item[(iii)] If $S,T\in\mathcal{C}(\R^n)$ with $S\subseteq T$ and $q\in  cp(\R^n)$, then
\beqq \ell_S(q)\leq\ell_T(q).\eeqq
\item[(iv)] If $T\in\mathcal{C}(\R^n)$, $q\in cp(\R^n)$ and $c>0$, then
\beqq \ell_T(cq)=\ell_{cT}(q)=c\ell_T(q).\eeqq
\item[(v)] If $T\in\mathcal{C}(\R^n)$, then
\beqq \ell_T: ( cp(\R^n),d_H) \rightarrow (\R_{\geq 0},|\cdot|)\eeqq
is continuous.
\item[(vi)] If $q\in cp(\R^n)$, then
\beqq \Gamma_q : (\mathcal{C}(\R^n),d_H)\rightarrow (\R_{\geq 0},|\cdot|),\; \; \Gamma_q(C)= \ell_C(q),\eeqq
is continuous.
\end{itemize}
\eprop

\bpf
\begin{itemize}
\item[(i)] If $F\in F(L)$, then $F$ cannot be translated into $\mathring{L}$. With $K\subseteq L$, $F$ also cannot be translated into $\mathring{K}$. Therefore, $F\in F(K)$. This means
\beqq F(L)\subseteq F(K).\eeqq
Analogously, one argues
\beqq F_{n+1}^{cp}(L)\subseteq F_{n+1}^{cp}(K).\eeqq

\item[(ii)] If $F\in F(cK)$, then $F$ cannot be translated into $c\mathring{K}$. Scaling $F$ and $K$ by the factor $\frac{1}{c}$, we conclude that $\frac{1}{c}F$ cannot be translated into $\mathring{K}$. Therefore,
\beqq \frac{1}{c}F\in F(K),\eeqq
and consequently $F\in cF(K)$. Analogously, we conclude that $F\in cF(K)$ implies $F\in F(cK)$. This finally implies
\beqq F(cK)=cF(K).\eeqq
Analogously one argues
\beqq F_{n+1}^{cp}(cK)=cF_{n+1}^{cp}(K).\eeqq

\item[(iii)] With $S\subseteq T$, we have $T^\circ \subseteq S^\circ$. Using Propositon \ref{Prop:IntroMinkowski}(ii), this implies
\beqq \mu_{S^\circ}(x)\leq \mu_{T^\circ}(x)\quad \forall x\in\R^n.\eeqq

This directly implies
\beqq \ell_S(q)\leq \ell_T(q).\eeqq

\item[(iv)] From Proposition \ref{Prop:IntroMinkowski}(iii), it follows
\beqq \mu_{T^\circ}(cx)=\mu_{(cT)^\circ}(x)=c\mu_{T^\circ}(x) \quad \forall x\in\R^n.\eeqq

This directly implies
\beqq \ell_T(cq)=\ell_{cT}(q)=c\ell_T(q).\eeqq

\item[(v)] It is enough to prove that
\beqq \mu_{T^\circ}: (\R^n,|\cdot|)\rightarrow (\R_{\geq 0},|\cdot|)\eeqq
is continuous. But this follows from Proposition \ref{Prop:IntroMinkowski}(iv).

\item[(vi)] Let $(T_j)_{j\in\N}$ be a sequence in $\mathcal{C}(\R^n)$ $d_H$-converging to $T\in \mathcal{C}(\R^n)$. This means for all $\eps>0$ that there is $j_0=j_0(\eps)\in\N$ with
\beq (1-\eps)T \subseteq T_j \subseteq (1+\eps)T\label{eq:Lemcontinuity2}\eeq
for all $j\geq j_0$.

We consider the sequence
\beqq (\Gamma_q(T_j))_{j\in\N} = (\ell_{T_j}(q))_{j\in\N}.\eeqq
Because of $\eqref{eq:Lemcontinuity2}$ and (iv)\&(v), we have for $\eps>0$ and big enough $j_0\in\N$
\beqq (1-\eps)\ell_T(q)=\ell_{(1-\eps)T}(q)\leq \ell_{T_j}(q) \leq \ell_{(1+\eps)T}(q)=(1+\eps)\ell_T(q)\eeqq
for all $j\geq j_0$. For $\eps\rightarrow 0$, this implies that
\beqq (\ell_{T_j}(q))_{j\in\N} \; \text{ converges to } \; \Gamma_q(T)=\ell_T(q).\eeqq

Therefore, $\Gamma_q$ is continuous.
\end{itemize}
\epf


For the next theorem, we denote by $M_{n+1}(K,T)$ the set of closed $(K,T)$-Minkowski billiard trajectories with at most $n+1$ bouncing points.

\bthm\label{Thm:onetoone}
Let $K,T\subset\R^n$ be convex bodies, where $T$ is additionally assumed to be strictly convex. Then, every $\ell_T$-minimizing closed $(K,T)$-Minkowski billiard trajectory is an $\ell_T$-minimizing element of $F_{n+1}^{cp}(K)$, and, conversely, every $\ell_T$-minimizing element of $F_{n+1}^{cp}(K)$ can be translated in order to be an $\ell_T$-minimizing closed $(K,T)$-Minkowski billiard trajectory.

Especially, one has
\beq \min_{q\in F^{cp}_{n+1}(K)}\ell_T(q) = \min_{q \in M_{n+1}(K,T)} \ell_T(q).\label{eq:onetooneesp}\eeq
\ethm

We remark that Theorem \ref{Thm:onetoone} is an existence result: In general, it guarantees the existence of $\ell_T$-minimizing elements of $F_{n+1}^{cp}(K)$, and furthermore, under the condition of strict convexity of $T$, it guarantees the existence of $\ell_T$-minimizing closed $(K,T)$-Minkowski billiard trajectories.

We note that in \cite[Theorem 2.2]{Rudolf2022}, we actually prove that \eqref{eq:onetooneesp} holds without the condition of strict convexity of $T$.

\bpf[Proof of Theorem \ref{Thm:onetoone}]
It is sufficient to prove the following two points:
\begin{itemize}
\item[(i)] Every closed $(K,T)$-Minkowski billiard trajectory is either in $F_{n+1}^{cp}(K)$ or there is a strictly $\ell_T$-shorter closed polygonal curve in $F_{n+1}^{cp}(K)$.
\item[(ii)] Every $\ell_{T}$-minimizing element of $F_{n+1}^{cp}(K)$ can be translated in order to be a closed $(K,T)$-Minkowski billiard tracjectory.
\end{itemize}

Ad (i): Let $q=(q_1,...,q_m)$ be a closed $(K,T)$-Minkowski billiard trajectory. From Proposition \ref{Prop:notranslation}, we conclude $q\in F(K)$. For $m\leq n+1$, we then have $q \in F_{n+1}^{cp}(K)$. If $m> n+1$, then, by \cite[Lemma 2.1(i)]{KruppRudolf2020}, there is a selection
\beqq \{i_1,...,i_{n+1}\}\subset \{1,...,m\} \; \text{ with } \; i_1 < ... < i_{n+1}\eeqq
such that the closed polygonal curve
\beqq (q_{i_1},...,q_{i_{n+1}})\eeqq
is in $F_{n+1}^{cp}(K)$. Because of Proposition \ref{Prop:IntroMinkowski}(i), we have
\beqq \ell_T((q_{i_1},...,q_{i_{n+1}})) < \ell_T(q).\eeqq

Ad (ii): By looking only at those members of $F_{n+1}^{cp}(K)$ that lie in an $n$-dimensional ball $B_R^n(0)$ in $\R^n$ of sufficiently large radius $R>0$ and which contains $K$, we get via the $d_H$-continuity of $\ell_T$ (cf.\;Proposition \ref{Prop:F&l}(v)) and a standard compactness argument, considering the identification between
\beqq \left(F_{n+1}^{cp}(K),d_H\right) \; \text{ and } \; \left(\left\{Q\in (\R^n)^{n+1} : Q \text{ cannot be translated into } \mathring{K}\right\},||\cdot||_2\right)\eeqq
by identifying
\beqq (q_1,...,q_m)\in F_{n+1}^{cp}(K)\; \text{ with } \; (q_1,...,q_m,...,q_m) \in (\R^n)^{n+1}\eeqq
and
\beqq (q_1,...,q_{n+1})\in (\R^n)^{n+1} \; \text{ with } \; (q_1,...,q_{n+1})\in F_{n+1}^{cp}(K)\eeqq
and the fact that
\beqq \left(F_{n+1}^{cp,*_R}(K):=\left\{q\in F_{n+1}^{cp}(K): q\subset B_R^n(0)\right\},d_H\right)\eeqq
therefore can be proven to be a compact metric subspace of the complete metric space $(P(\R^n),d_H)$, that
\beqq F_{n+1}^{cp}(K)/\{\text{translations}\}\eeqq
possesses at least one element of minimal $\ell_{T}$-length, say $\Delta_{\min}$.

We show that there is a translate $\Delta_{\min}'$ of $\Delta_{\min}$ which is a closed $(K,T)$-Minkowski billiard tracjecory.

Indeed, $\Delta_{\min}$ as element of $F_{n+1}^{cp}(K)$ cannot be translated into $\mathring{K}$. Thus, with \cite[Lemma 2.1(ii)]{KruppRudolf2020} and the $\ell_T$-minimality of $\Delta_{\min}$, there is a translate $\Delta_{\min}'$ of $\Delta_{\min}$ given by vertices $q_1',...,q_m'\in\partial K$, $m\leq n+1$, and there are $K$-supporting hyperplanes $H_1,...,H_m$ through $q_1',...,q_m'$ such that
\beqq H_1^+ \cap ... \cap H_m^+\eeqq
is nearly bounded in $\R^n$, where $H_1^+,...,H_m^+$ are the closed half-spaces of $\R^n$ containing $K$ and which are bounded by $H_1,...,H_m$. Let $(q_{j-1}',q_j',q_{j+1}')$ be an arbitrary polygonal curve segment of $\Delta_{\min}'$. If this polygonal curve segment is not fulfilling the weak Minkowski billiard reflection rule with respect to $H_j$, meaning that $q_j'$ is not minimizing
\beqq \Sigma_j(q_j'^*)=\mu_{T^\circ}(q_j'^*-q_{j-1}')+\mu_{T^\circ}(q_{j+1}'-q_j'^*)\eeqq
over all $q_j'^*\in H_j$, then we find a $\widetilde{q}_j\in H_j$, $\widetilde{q}_j\neq q_j'$, such that the $\ell_T$-length of the polygonal curve segment $(q_{j-1}',\widetilde{q}_j,q_{j+1}')$ is less than the $\ell_T$-length of the polygonal curve segment $(q_{j-1}',q_j',q_{j+1}')$. We replace the polygonal curve segment $(q_{j-1}',q_j',q_{j+1}')$ within $\Delta_{\min}'$ by $(q_{j-1}',\widetilde{q}_j,q_{j+1}')$. By \cite[Lemma 2.1(ii)]{KruppRudolf2020}, the new closed polygonal curve
\beqq \widetilde{\Delta}_{\min}'=(q_1',...,q_{j-1}',\widetilde{q}_j,q_{j+1}',...,q_m')\eeqq
still cannot be translated into $\mathring{K}$, while
\beqq \ell_T(\widetilde{\Delta}_{\min}')<\ell_T(\Delta_{\min}').\eeqq
But this is a contradiction to the $\ell_T$-minimality of $\Delta_{\min}'$. Therefore, every polygonal curve segment $(q_{j-1}',q_j',q_{j+1}')$ of $\Delta_{\min}'$ is fulfilling the weak Minkowski billiard reflection rule. Consequently, referring to Theorem \ref{Thm:weakstrongbilliards} ($T$ is strictly convex), $\Delta_{\min}'$ is an $\ell_T$-minimizing closed $(K,T)$-Minkowski billiard trajectory.
\epf

\section{Proof of Theorem \ref{Thm:RegularityResult}}\label{Sec:Proof1}



For the proof of Theorem \ref{Thm:RegularityResult}, we need the following generalization of \cite[Lemma 3.2]{KruppRudolf2020}:

\blem\label{Lem:normalconerelation}
Let $K\subset\R^n$ be a convex body and $\{q_1,...,q_m\}$ a set of boundary points of $K$. Further, let $U$ be the convex cone spanned by outer unit normal vecors
\beqq n_K(q_1),...,n_K(q_m) \; \text{ in } \; N_K(q_1),...,N_K(q_m).\eeqq
Then, we have
\beqq N_K(q_j)\cap U = N_K(q_j) \cap N_{K\cap (U+q_j)}(q_j)\eeqq
for all $j\in\{1,...,m\}$.
\elem

\bpf
From
\beqq N_{K\cap (U+q_j)}(q_j) \subseteq U\eeqq
for all $j\in\{1,...,m\}$ follows
\beqq N_K(q_j)\cap U \supseteq N_K(q_j) \cap N_{K\cap (U+q_j)}(q_j)\eeqq
for all $j\in\{1,...,m\}$.

Let $j\in\{1,...,m\}$ be arbitrarily chosen. Let $n$ be a nonzero vector in
\beqq N_K(q_j)\cap U.\eeqq
Then
\beqq n\in N_K(q_j),\text{ i.e., }\; \langle n,x-q_j\rangle \leq 0 \;\; \forall x\in K,\eeqq
and $n\in U$. Because of
\beqq K\cap (U+q_j)\subseteq K,\eeqq
this implies
\beqq \langle n,x-q_j \rangle \leq 0 \;\; \forall x\in K\cap (U+q_j),\; n\in U,\; n\in N_K(q_j).\eeqq
From that, we conclude
\beqq n\in N_{K\cap (U+q_j)}(q_j) \; \text{ and } \; n\in N_K(q_j),\eeqq
and therefore
\beqq n\in N_{K\cap (U+q_j)}(q_j)\cap N_K(q_j).\eeqq
Consequently,
\beqq N_K(q_j)\cap U \subseteq N_K(q_j) \cap N_{K\cap (U+q_j)}(q_j)\eeqq
for all $j\in\{1,...,m\}$.
\epf

We come to the proof of Theorem \ref{Thm:RegularityResult}:

\bpf[Proof of Theorem \ref{Thm:RegularityResult}]
By Proposition \ref{Prop:normalvectorsspanning}, $U$ is a linear subspace of $\R^n$ with
\beqq \dim U \leq m-1 \leq n,\eeqq
where the last inequality follows from Theorem \ref{Thm:onetoone}. By Proposition \ref{Prop:GenBezdek}, there is a selection
\beqq \{i_1,...,i_{\dim U +1}\}\subseteq \{1,...,m\}\eeqq
with
\beqq \left\{q_{i_1},...,q_{i_{\dim U +1}}\right\}\in F(K).\eeqq
Without loss of generality, we can assume
\beqq i_1<...<i_{\dim U +1}\eeqq
and define the closed polygonal curve
\beqq \widetilde{q}=\left(q_{i_1},...,q_{i_{\dim U +1}}\right).\eeqq
For
\beqq \dim U +1<m,\eeqq
it follows by Proposition \ref{Prop:IntroMinkowski}(i) (requires strict convexity of $T$) that
\beqq \ell_T(\widetilde{q})<\ell_T(q).\eeqq
But with Theorem \ref{Thm:onetoone}, this is a contradiction to the $\ell_T$-minimality of $q$. Therefore, we conclude
\beqq \dim U =m-1.\eeqq

Let us denote by $H_1^+,...,H_m^+$ the closed half-spaces of $\R^n$ containing $K$ and which are bounded by $H_1,...,H_m$. By Proposition \ref{Prop:noteinproof}, we conclude that we can write
\beq H_j=\left(H_j\cap U\right)\oplus W\; \text{ and }\;H_j^+=\left(H_j^+\cap U\right)\oplus W\label{eq:regularsplitting}\eeq
for all $j\in\{1,...,m\}$, where $W$ is the orthogonal complement to $U$ in $\R^n$, and that
\beqq \bigcap_{j=1}^m \left(H_j^+\cap U\right) \text{ is bounded in }U, \; \bigcap_{j=1}^m H_j^+ \text{ is nearly bounded in }\R^n .\eeqq

By the definition of $U$, we have
\beqq n_K(q_j)\in N_K(q_j)\cap U\quad \forall j\in\{1,...,m\}\eeqq
and therefore
\beqq \dim (N_K(q_j)\cap U)\geq 1 \quad \forall j\in\{1,...,m\}.\eeqq
Let us assume there is an $i\in\{1,...,m\}$ such that
\beqq \dim (N_K(q_i)\cap U)>1.\eeqq
Then, using Lemma \ref{Lem:normalconerelation}, i.e.,
\beqq N_K(q_i)\cap U = N_K(q_i) \cap N_{K\cap (U+q_i)}(q_i),\eeqq
it follows
\beqq \dim \left(N_K(q_i)\cap N_{K\cap (U+q_i)}(q_i)\right)>1,\eeqq
and because of \cite[Lemma 3.1]{KruppRudolf2020} (for $d=m-1$ and $k=m$), we can find a unit vector
\beq n_i^{pert}\in N_{K}(q_i)\cap N_{K\cap (U+q_i)}(q_i) \; \text{ with }\; n_i^{pert}\neq n_i\label{eq:choicepiperturbed}\eeq
such that 
\beq H_{i,U}^{pert,+}\cap \left(\bigcap_{j=1,j\neq i}^m \left(H_{j}^+\cap U\right)\right)\label{eq:regularityProp2}\eeq
remains bounded in $U$, where we denote by $H_{i,U}^{pert,+}$ the closed half-space of $U$ that contains $\pi_U(K)$, where $\pi_U$ is the orthogonal projection (projection along $W$) onto $U$, and which is bounded by $H_{i,U}^{pert}$ which is the hyperplane in $U$ through $\pi_U(q_i)$ that is normal to $n_i^{pert}$. Since by Proposition \ref{Prop:onlyonehyperplane}, the weak Minkowski billiard reflection rule in $q_i$ (cf.\;Theorem \ref{Thm:weakstrongbilliards}) is no longer satisfied with respect to the perturbed hyperplane
\beqq H_i^{pert}:=H_{i,U}^{pert}\oplus W,\eeqq
the bouncing point $q_i$ can be moved along
\beqq H_{i,U}^{pert}+\left(q_i-\pi_U(q_i)\right)\subset H_i^{pert},\eeqq
say to $q_i^*$, in order to reduce the length of the polygonal curve segment
\beqq (q_{i-1},q_i,q_{i+1}).\eeqq
We define the closed polygonal curve
\beqq q^*:=(q_1,...,q_{i_1},q_i^*,q_{i+1},...,q_m)\eeqq
and argue that $q^*\in F(K)$: With the boundedness of \eqref{eq:regularityProp2} in $U$, it follows with
\beq  H_i^{pert,+}:=H_{i,U}^{pert,+}\oplus W\label{eq:regularityProp21}\eeq
and \eqref{eq:regularsplitting} the nearly boundedness of
\beq H_i^{pert,+}\cap \left(\bigcap_{j=1,j\neq i}^m H_j^+\right)\label{eq:regularityProp3}\eeq
in $\R^n$.

Indeed, when the intersection in \eqref{eq:regularityProp2} is bounded in $U$, then there is a hyperplane $H$ in $U$ such that the intersection lies between $H$ and $H+d$ for an appropriate $d\in U$. Then it follows with \eqref{eq:regularsplitting} and \eqref{eq:regularityProp21} that
\begin{align*}
&H_i^{pert,+}\cap \left(\bigcap_{j=1,j\neq i}^m H_j^+\right)\\
=&\left(H_{i,U}^{pert,+}\oplus W\right)\cap \left(\bigcap_{j=1,j\neq i}^m \left(\left(H_j^+\cap U\right)\oplus W\right)\right)\\
=&\left(H_{i,U}^{pert,+}\cap \left(\bigcap_{j=1,j\neq i}^m \left(H_j^+\cap U\right)\right)\right)\oplus W
\end{align*}
lies between the hyperplanes
\beqq H\oplus W \; \text{ and } \; (H+d)\oplus W.\eeqq

Since $H_i^{pert}$ is a $K$-supporting hyperplane through $q_i$ (what follows from the fact that by \eqref{eq:choicepiperturbed} its outer unit normal vector $n_i^{pert}$ is an element of $N_K(q_i)$), we conclude that $K$ is a subset of the intersection in \eqref{eq:regularityProp3}. Then, it follows from the nearly boundedness (in $\R^n$) of the intersection in \eqref{eq:regularityProp3} together with \cite[Lemma 2.1(ii)]{KruppRudolf2020} that
\beqq q^*\in F(K).\eeqq
By referring to Theorem \ref{Thm:onetoone}, from
\beqq \ell_T(q^*)<\ell_T(q),\eeqq
we derive a contradiction to the $\ell_T$-minimality of $q$.

Therefore:
\beqq \dim (N_K(q_i)\cap U)=1.\eeqq
\epf

We remark that for the proof of Theorem \ref{Thm:RegularityResult}, the smoothness of $T$ is a necessary condition. It guarantees the application of Proposition \ref{Prop:onlyonehyperplane} and the boundedness of
\beq \bigcap_{j=1}^m (H_j^+\cap U).\label{eq:boundedness}\eeq
Without the smoothness of $T$, from Proposition \ref{Prop:noteinproof}, we know of \eqref{eq:boundedness}'s nearly boundedness, but that is not enough in order to utilize Lemma \cite[Lemma 3.1]{KruppRudolf2020}.

\section{Proof of Theorem \ref{Thm:maximxallyspanning}}\label{Sec:Proof2}

In the proof of Theorem \ref{Thm:maximxallyspanning} we will use the following proposition:

\bprop\label{Prop:inversion}
Let $K,T\subset\R^n$ be convex bodies. Let $q=(q_1,...,q_m)$ be a closed weak $(K,T)$-Minkowski billiard trajectory. Then, $-q$ is a closed weak $(-K,-T)$-Minkowski billiard trajectory with
\beqq \ell_T(q)=\ell_{-T}(-q).\eeqq
\eprop

\bpf
Since $q$ is a closed weak $(K,T)$-Minkowski billiard trajectory, there are $K$-supporting hyperplanes $H_1,...,H_m$ through $q_1,...,q_m$ such that $q_j$ minimizes $\Sigma_j(q_j^*)$ over all $q_j^*\in H_j$ for all $j\in\{1,...,m\}$. More precisely, we have the following for every $j\in\{1,...,m\}$: $q_j$ minimizes
\beqq \Sigma(q_j^*)=\mu_{T^\circ}(q_j^*-q_{j-1}) + \mu_{T^\circ}(q_{j+1}-q_j^*)=\langle q_j^*-q_{j-1},p_{j-1}^*\rangle + \langle q_{j+1}-q_j^*,p_j^*\rangle\eeqq
over all $q_j^*\in H_j$, where $p_{j-1}^*,p_j^*\in\partial T$ (possibly not uniquely determined) fulfill
\beqq q_j^*-q_{j-1}\in N_T(p_{j-1}^*)\;\text{ and }\; q_{j+1}-q_j^*\in N_T(p_j^*).\eeqq
Because of
\beqq \langle q_j^*-q_{j-1},p_{j-1}^*\rangle = \langle -q_j^*-(-q_{j-1}),-p_{j-1}^*\rangle\eeqq
and
\beqq \langle q_{j+1}-q_j^*,p_j^*\rangle = \langle -q_{j+1}-(-q_j^*),-p_j^*\rangle\eeqq
as well as
\beqq q_j^*-q_{j-1}\in N_T(p_{j-1}^*) \Leftrightarrow -q_j^*-(-q_{j-1})\in N_{-T}(-p_{j-1}^*)\eeqq
and
\beqq q_{j+1}-q_j^*\in N_T(p_j^*) \Leftrightarrow -q_{j+1}-(-q_j^*)\in N_{-T}(-p_j^*),\eeqq
we conclude that $-q_j$ minimizes
\beqq \langle -q_j^*-(-q_{j-1}),-p_{j-1}^*\rangle + \langle -q_{j+1}-(-q_j^*),-p_j^*\rangle\eeqq
and therefore
\beqq  \Sigma_j(-q_j^*)=\mu_{(-T)^\circ}(-q_j^*-(-q_{j-1}))+\mu_{(-T)^\circ}(-q_{j+1}-(-q_j^*))\eeqq
over all $-q_j^*\in -H_j$. It follows that $-q$ is a closed weak $(-K,-T)$-Minkowski billiard trajectory.

We finally argue that
\beqq \ell_{T}(q)=\ell_{-T}(-q):\eeqq
We have
\beqq \ell_T(q)=\sum_{j=1}^m \mu_{T^\circ}(q_{j+1}-q_j)=\langle q_{j+1}-q_j,p_j\rangle ,\eeqq
where $p_j\in\partial T$ (possibly not uniquely determined) fulfills
\beqq q_{j+1}-q_j\in N_T(p_j).\eeqq
Using
\beqq q_{j+1}-q_j\in N_T(p_j) \Leftrightarrow -q_{j+1}-(-q_j)\in N_{-T}(-p_j),\eeqq
we conclude
\begin{align*}
\ell_T(q)=\sum_{j=1}^m \langle q_{j+1}-q_j,p_j\rangle &= \sum_{j=1}^m\langle -q_{j+1}-(-q_j),-p_j\rangle\\
&=\sum_{j=1}^m\mu_{(-T)^\circ}(-q_{j+1}-(-q_j))\\
&=\ell_{-T}(-q).
\end{align*}
\epf

\bpf[Proof of Theorem \ref{Thm:maximxallyspanning}]
We have
\beqq \dim V \leq m-1,\eeqq
since in general $m$ points can maximally span an $(m-1)$-dimensional cone/space.

Let us assume
\beqq \dim V<m-1.\eeqq
Since $K\cap V$ is the inclusion minimal affine section of $K$ containing $q$, we conclude that the convex cone spanned by
\beq q_2-q_1,...,q_m-q_{m-1},q_1-q_m\label{eq:vectordifferences}\eeq
is $V_0$, where $V_0$ is the linear subspace of $\R^n$ underlying $V$ ($\dim V_0 < m-1$).

Indeed, we argue similarly to within the proof of Proposition \ref{Prop:normalvectorsspanning}: we show that the convex cone spanned by \eqref{eq:vectordifferences} is in fact a linear subspace of $\R^n$. For that, we notice that
\beqq (q_2-q_1) + ... + (q_m-q_{m-1}) + (q_1 - q_m) = 0\eeqq
and therefore that
\beqq s_1 (q_2-q_1) + ... + s_{m-1}(q_m-q_{m-1}) + s_m (q_1 - q_m) = 0,\eeqq
where we defined
\beqq s_j:= \frac{1}{m}\quad \forall j\in\{1,...,m\}.\eeqq
Since
\beqq \sum_{j=1}^m s_j = 1 \; \text{ and } \; n_j\neq 0\;\; \forall j\in\{1,...,m\},\eeqq
we have that $0$ is within the relative interior of the convex cone spanned by the vectors \eqref{eq:vectordifferences}. Consequently, the convex cone in fact is the linear subspace $V_0$ of $\R^n$ which underlies $V$.

Let $p=(p_1,...,p_m)$ be the uniquely determined closed dual billiard trajectory of $q$ in $T$ (cf.\;Proposition \ref{Prop:dualbilliardtrajunique}). Then, because of Proposition \ref{Prop:onlyonehyperplane} (and therefore, there is a $\mu_j>0$ such that $p_j-p_{j-1}=-\mu_j n_K(q_j)$, where $n_K(q_j)\in N_K(q_j)$, for all $j\in\{1,...,m\}$), the pair $(q,p)$ fulfills
\beq \begin{cases} q_{j+1}-q_j=\lambda_j n_T(p_j)\in N_T(p_j),\; \lambda_j >0,\\ p_{j+1}-p_j=-\mu_{j+1}n_K(q_{j+1}),\; \mu_{j+1}> 0,\end{cases}\label{eq:maximallyspanning1}\eeq
for all $j\in\{1,...,m\}$. By Proposition \ref{Prop:dualbilliard}, $p$ is a closed $(T,-K)$-Minkowski billiard trajectory (which requires the strict convexity and smoothness of $K$). By Proposition \ref{Prop:lengthdualbilliard}, we have
\beqq \ell_T(q)=\ell_{-K}(p).\eeqq

From \eqref{eq:maximallyspanning1}, we conclude that the convex cone spanned by $n_T(p_1),...,n_T(p_m)$ is $V_0$. Then, by Proposition \ref{Prop:GenBezdek}, there is a selection
\beqq \{i_1,...,i_{\dim V_0 +1}\}\subset \{1,...,m\} \; \text{ with } \; i_1<...<i_{\dim V_0 +1}\eeqq
such that the closed polygonal curve
\beqq \widetilde{p}=(p_{i_1},...,p_{i_{\dim V_0 +1}})\eeqq
is in $F(T)$. Because of the strict convexity of $K$ and Proposition \ref{Prop:IntroMinkowski}(i), we have
\beqq \ell_{-K}(\widetilde{p})<\ell_{-K}(p).\eeqq
Applying Theorem \ref{Thm:onetoone}, there has to be an $\ell_{-K}$-minimizing closed $(T,-K)$-Minkowski billiard trajectory $p^*$ with
\beqq \ell_{-K}(p^*)\leq \ell_{-K}(\widetilde{p}).\eeqq
Let $q^*$ be its dual billiard trajectory on $-K$ which by Propositions \ref{Prop:lengthdualbilliard} and \ref{Prop:dualbilliard} is a closed $(-K,-T)$-Minkowski billiard trajectory with
\beqq \ell_{-T}(q^*)=\ell_{-K}(p^*).\eeqq
Then, it follows by Proposition \ref{Prop:inversion} that $-q^*$ is a closed $(K,T)$-Minkowski billiard trajectory with
\beqq \ell_T(-q^*)=\ell_{-T}(q^*)=\ell_{-K}(p^*)\leq \ell_{-K}(\widetilde{p})<\ell_{-K}(p)=\ell_T(q).\eeqq
This is a contradiction to the $\ell_T$-minimality of $q$.

Therefore,
\beqq \dim V = m-1.\eeqq
\epf

\section{Examples}\label{Sec:Examples}

\underline{Example A}: We consider the following example: Let $K\subset\R^2$ be the triangle with vertices
\beqq (1,0),(0,1),(-1,0)\eeqq
and $T\subset\R^2$ the square with vertices
\beqq (1,1),(-1,1),(-1,-1),(1,-1).\eeqq

\begin{figure}[h!]
\centering
\def\svgwidth{390pt}
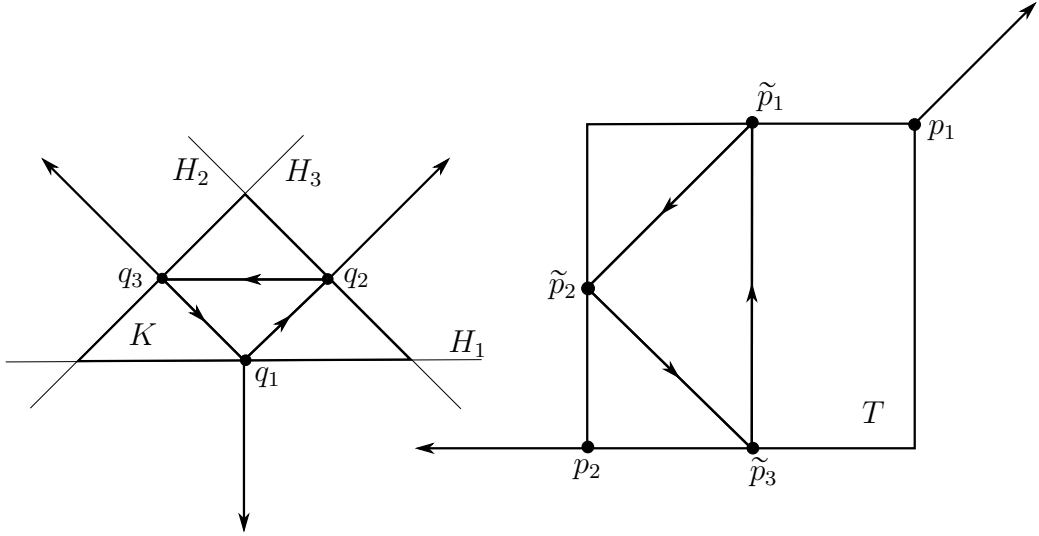
\caption[Example A]{Example A: $q=(q_1,q_2,q_3)$ is a closed weak $(K,T)$-Minkowski billiard trajectory which fulfills the weak Minkowski billiard reflection rule with respect to the $K$-supprting hyperplanes $H_1,H_2,H_3$, but it is not a strong one, i.e., there is no closed polygonal curve $p$ such that the pair $(q,p)$ fulfills \eqref{eq:System}.}
\label{img:Counterexample}
\end{figure}

Let $q=(q_1,q_2,q_3)$ be a closed polygonal curve with
\beqq q_1=(0,0),\; q_2=\left(\frac{1}{2},\frac{1}{2}\right),\; q_3=\left(-\frac{1}{2},\frac{1}{2}\right).\eeqq
We denote by $H_1,H_2,H_3$ the $K$-supporting hyperplanes through $q_1,q_2,q_3$. We claim that $(q_1,q_2,q_3)$ is a closed weak $(K,T)$-Minkowski billiard trajectory fulfilling the weak Minkowski billiard reflection rule with respect to the hyperplanes $H_1,H_2,H_3$. Exemplary, we show that the weak Minkowski billiard reflection rule is satisfied for the polygonal curve segment from $q_1$ over $q_2$ to $q_3$. For that, we show that $q_2$ minimizes
\beqq \Sigma_2(q_2^*)=\mu_{T^\circ}(q_2^*-q_1)+\mu_{T^\circ}(q_3-q_2^*)\eeqq
over all $q_2^*\in H_2$: We have
\beqq \Sigma_2(q_2)=\langle q_2-q_1,p_1\rangle + \langle q_3-q_2,p_2\rangle\eeqq
with $p_1=(1,1)$ and $p_2=(-1,-1)$ (the idea behind this example is that $T$ allows to choose $p_1$ and $p_2$ on $\partial T$ such that their connecting line is othogonal to $H_2$). Since
\beqq \langle q_2-q_2^*,p_2-p_1\rangle = 0\; \text{ for all }q_2^*\in H_2,\eeqq
we conclude for any $q_2^*\in H_2$
\allowdisplaybreaks{\begin{align*}
\Sigma_2(q_2)&= \langle q_2-q_1,p_1\rangle + \langle q_3-q_2,p_2\rangle + \langle q_2-q_2^*,p_2-p_1\rangle\\
&=\langle q_2^*-q_1,p_1\rangle + \langle q_3-q_2^*,p_2\rangle\\
&=  \langle q_2^*-q_1,p_1^*\rangle + \langle q_3-q_2^*,p_2^*\rangle + \langle q_2^*-q_1,p_1-p_1^*\rangle + \langle q_3-q_2^*,p_2-p_2^*\rangle,
\end{align*}}%
where $p_1^*,p_2^*\in\partial T$ are chosen to fulfill
\beqq q_2^*-q_1\in N_T(p_1^*)\;\text{ and }\; q_3-q_2^*\in N_T(p_2^*)\eeqq
(by this condition possibly not uniquely determined). From the convexity of $T$, it follows
\beqq  \langle q_2^*-q_1,p_1-p_1^*\rangle \leq 0\;\text{ and }\;  \langle q_3-q_2^*,p_2-p_2^*\rangle \leq 0\eeqq
and therefore
\beqq \Sigma_2(q_2)\leq \langle q_2^*-q_1,p_1^*\rangle + \langle q_3-q_2^*,p_2^*\rangle = \mu_{T^\circ}(q_2^*-q_1) + \mu_{T^0}(q_3-q_2^*) = \Sigma_2(q_2^*),\eeqq
where we used Proposition \ref{Prop:Minkowskisupport}. Consequently, $q_2$ minimizes $\Sigma_2(q_2^*)$ over all $q_2^*\in H_2$.

Similarly, one could prove that the polygonal curve segment from $q_2$ over $q_3$ to $q_1$--by choosing $p_2=(-1,1)$ and $p_3=(1,-1)$--as well as the one from $q_3$ over $q_1$ to $q_2$--by choosing $p_3=(1,-1)$ and $p_1=(1,1)$--fulfills the weak Minkowski billiard reflection rule. This gives us the idea behind choosing this example: For every polygonal curve segment of $q$ consisting of three consecutive bouncing points, the weak Minkowski billiard reflection rule is satisfied, but it is not possible to find $p_1$, $p_2$ and $p_3$ uniformly in order to construct a dual billiard trajectory.

In fact, we claim that there is no closed polygonal curve $p=(\widetilde{p}_1,\widetilde{p}_2,\widetilde{p}_3)$ with vertices on $\partial T$ such that the pair $(q,p)$ fulfills \eqref{eq:System}. If this would be the case, then there would be $\widetilde{p}_1,\widetilde{p}_2,\widetilde{p}_3\in\partial T$ with
\beq (\widetilde{p}_2-\widetilde{p}_1)+(\widetilde{p}_3-\widetilde{p}_2)+(\widetilde{p}_1-\widetilde{p}_3)=0\label{eq:ExampleA1}\eeq
and, additionally, there would be outer unit vectors $n_K(q_1),n_K(q_2),n_K(q_3)$ at $K$ normal to $H_1,H_2,H_3$ such that
\begin{align}
\label{eq:2columns}
&q_2-q_1\in N_T(\widetilde{p}_1)& \widetilde{p}_2-\widetilde{p}_1=-\mu_2 n_K(q_2)\notag\\
&q_3-q_2\in N_T(\widetilde{p}_2) & \widetilde{p}_3-\widetilde{p}_2=-\mu_3 n_K(q_3)\\
&q_1-q_3\in N_T(\widetilde{p}_3)  & \widetilde{p}_1-\widetilde{p}_3=-\mu_1 n_K(q_1)\notag
\end{align}
where $\mu_1,\mu_2,\mu_3\geq 0$. From $\widetilde{p}_1,\widetilde{p}_2,\widetilde{p}_3\in\partial T$ together with \eqref{eq:ExampleA1} and the second column in \eqref{eq:2columns}, we then could conclude
\beqq \widetilde{p}_1=(0,1),\;\widetilde{p}_2=(-1,0),\;\widetilde{p}_3=(0,-1).\eeqq
But this would imply that the conditions within the first column in \eqref{eq:2columns} cannot be satisfied. Consequently, the pair $(q,p)$ does not fulfill \eqref{eq:System}.

Summarized, this example shows that, without requiring $T$ to be strictly convex, it can happen that a closed weak $(K,T)$-Minkowski billiard trajectory is not a closed strong $(K,T)$-Minkowski billiard trajectory.\qed

\underline{Example B}: Let $K\subset\R^2$ be the square given by the vertices
\beqq (1,0), (0,1), (-1,0), (0,-1)\eeqq
and $T\subset\R^2$ the union of the square given by the vertices
\beqq (1,1), (-1,1), (-1,-1), (1,-1)\eeqq
and the two balls
\beqq B_1^2  + (1,0) \; \text{ and } \; B_1^2  + (-1,0).\eeqq
Then, $T$ is smooth, but not strictly convex.

\begin{figure}[h!]
\centering
\def\svgwidth{390pt}
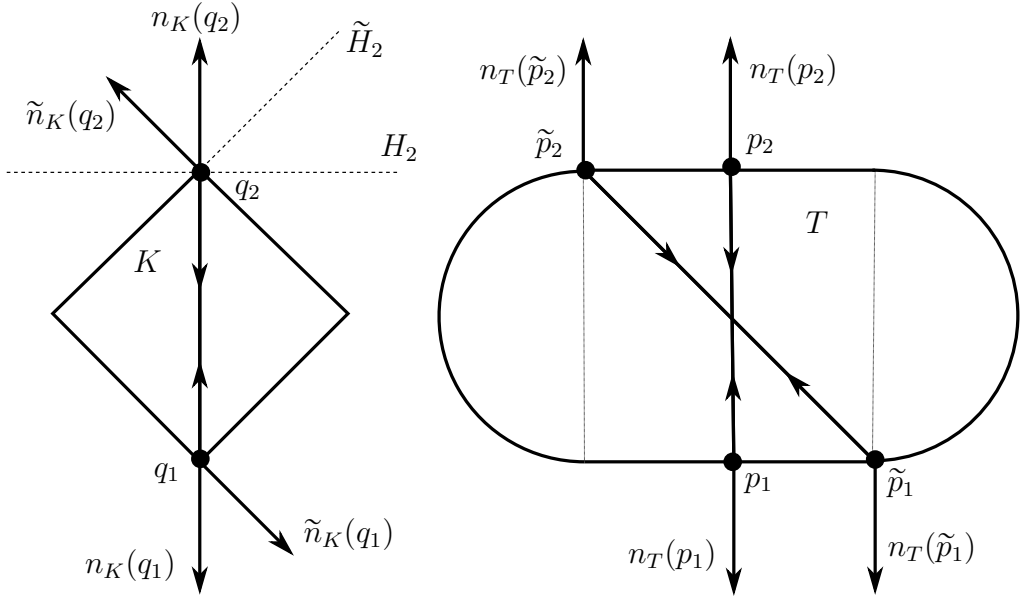
\caption[Example B]{Example B: The closed $(K,T)$-Minkowski billiard trajectory $q$ has $p=(p_1,p_2)$ as well as $\widetilde{p}=(\widetilde{p}_1,\widetilde{p}_2)$ as closed dual billiard trajectory in $T$. Furthermore, $q_2$ minimizes \eqref{eq:cs1} over all $\widebar{q}_2\in H_2$ as well as over all $\widebar{q}_2\in \widetilde{H}_2$.}
\label{img:cs1}
\end{figure}

Then, we can easily check that
\beqq q=(q_1,q_2) \; \text{ with } \; q_1=(0,-1) \; \text{ and } \; q_2=(0,1)\eeqq
is a closed $(K,T)$-Minkowski billiard trajectory: \eqref{eq:System} is satisfied for the pair $(q,p)$ for
\beqq p=(p_1,p_2) \; \text{ with } \; p_1=q_1 \; \text{ and } \; p_2=q_2\eeqq
with the corresponding outer unit normal vectors
\beqq n_K(q_1)=n_T(p_1)=(0,-1),\; n_K(q_2)=n_T(p_2)=(0,1).\eeqq

Since every closed (strong) Minkowski billiard trajectory is also a weak one (cf.\;Theorem \ref{Thm:weakstrongbilliards}), $q_2$ minimizes
\beq \mu_{T^\circ}(\widebar{q}_2-q_1)+\mu_{T^\circ}(q_1-\widebar{q}_2)\label{eq:cs1}\eeq
over all $\widebar{q}_2\in H_2$, where $H_2$ is the $K$-supporting horizontal line through $q_2$.

However, \eqref{eq:System} is also satisfied for the pair $(q,\widetilde{p})$ for
\beqq \widetilde{p}=(\widetilde{p}_1,\widetilde{p}_2) \; \text{ with } \; \widetilde{p}_1=(1,-1) \; \text{ and } \; \widetilde{p}_2=(-1,1)\eeqq
with the corresponding outer unit normal vectors
\beqq \widetilde{n}_K(q_1)=\left(\frac{1}{\sqrt{2}},-\frac{1}{\sqrt{2}}\right),\; \widetilde{n}_K(q_2)=\left(-\frac{1}{\sqrt{2}},\frac{1}{\sqrt{2}}\right),\eeqq
\beqq n_T(\widetilde{p}_1)=n_T(p_1),\; \widetilde{n}_T(p_2)=n_T(p_2).\eeqq
Then, again referring to Theorem \ref{Thm:weakstrongbilliards}, $q_2$ minimizes \eqref{eq:cs1} over all $\widebar{q}_2\in \widetilde{H}_2$, where $\widetilde{H}_2$ is the $K$-supporting line through $q_2$ with slope $1$.

Summarized, this example shows that, without requiring $T$ to be strictly convex, it can happen that the closed dual billiard trajectory of a closed $(K,T)$-Minkowski billiard trajectory is not uniquely determined. Furthermore, it shows that, without requiring $T$ to be strictly convex, it can happen that the $K$-supporting hyperplanes corresponding via the weak Minkowski billiard reflection rule to a closed weak $(K,T)$-Minkowski billiard trajectory are not uniquely determined.\qed

\underline{Example C}: Let $K\subset\R^2$ be the triangle given by the vertices
\beqq (1,0),(-1,2),(-1,-2)\eeqq
and $T$ the intersection of the two balls
\beqq B_1^2 +\left(-\frac{1}{2},0\right) \; \text{ and } \; B_1^2 +\left(\frac{1}{2},0\right).\eeqq
Then, $T$ is strictly convex, but not smooth.

\begin{figure}[h!]
\centering
\def\svgwidth{380pt}
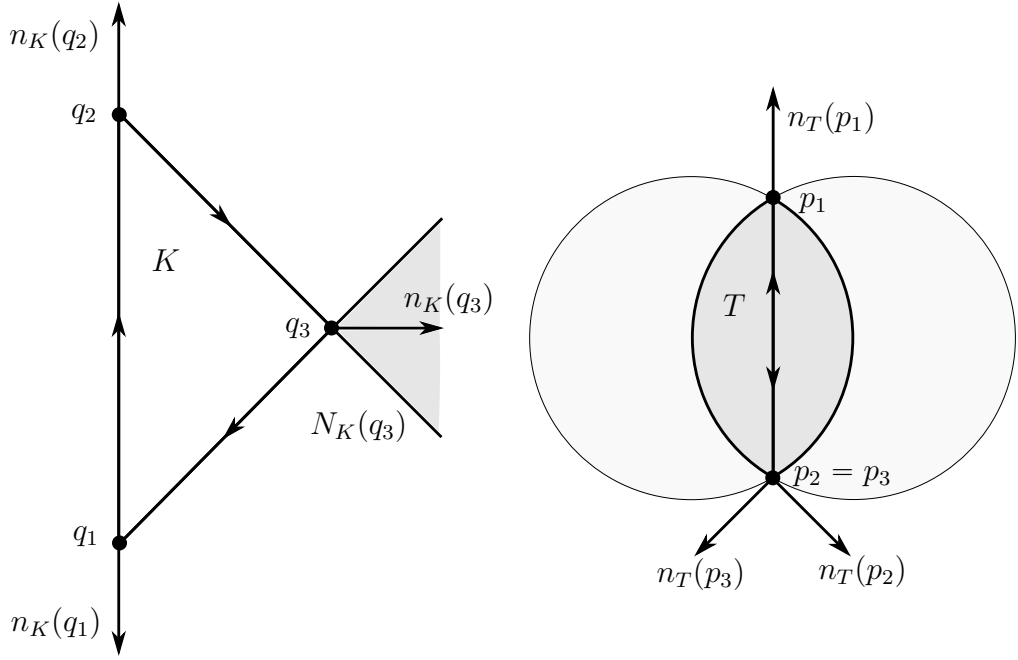
\caption[Example C]{Example C: $q=(q_1,q_2,q_3)$ is a closed $(K,T)$-Minkowski billiard trajectory, no matter which outer normal unit vector in $q_3$ is chosen in order to show that $(q,p)$ fulfills \eqref{eq:System}. The vertices $p_2,p_3$ of $q$'s closed dual billiard trajectory $p$ coincide. The convex cone spanned by $n_K(q_1),n_K(q_2),n_K(q_3)$ does not span the whole $\R^2$.}
\label{img:exampleB2}
\end{figure}

Then, we can easily check that
\beqq q=(q_1,q_2,q_3) \; \text{ with } \; q_1=(-1,-2),\; q_2=(-1,2) \; \text{and} \; q_3=(1,0)\eeqq
is a closed $(K,T)$-Minkowski billiard trajectory: \eqref{eq:System} is satisfied for the pair $(q,p)$ for
\beqq p=(p_1,p_2,p_3) \; \text{ with } \; p_1=(0,x) \; \text{and} \; p_2=p_3=(0,-x)\eeqq
with the corresponding outer unit normal vectors
\beqq n_K(q_1)=(0,-1),\; n_K(q_2)=(0,1), \; n_K(q_3)=(1,0),\eeqq
\beqq n_T(p_1)=(0,1),\; n_T(p_2)=\left(\frac{1}{\sqrt{2}},-\frac{1}{\sqrt{2}}\right),\; n_T(p_3)=\left(-\frac{1}{\sqrt{2}},-\frac{1}{\sqrt{2}}\right).\eeqq

However, we notice that \eqref{eq:System} is also satisfied when we replace $n_K(q_3)=(1,0)$ by any other outer unit normal vector within the normal cone $N_K(q_3)$.

Summarized, this example shows that, without requiring $T$ to be smooth, it can happen that the closed dual billiard trajectory corresponding to a closed $(K,T)$-Minkowski billiard trajectory is not a polygonal curve in the sense of Footnote \ref{foot:polygonalline}. This implies that, without requiring $T$ to be smooth, it can happen that the $K$-supporting hyperplanes corresponding via the weak Minkowski billiard reflection rule to a closed weak $(K,T)$-Minkowski billiard trajectory are not uniquely determined. Furthermore, this example shows that, without requiring $T$ to be smooth, it can happen that the convex cone spanned by the outer unit vectors normal to the hyperplanes which correspond via the weak Minkowski billiard reflection rule to a closed weak $(K,T)$-Minkowski billiard trajectory does not is the whole $\R^n$.




\underline{Example D}: Let $K\subset\R^2$ be the square given by the vertices
\beqq (1,1), (-1,1), (-1,-1), (1,-1)\eeqq
and $\widetilde{T}\subset\R^2$ the triangle given by the vertices
\beqq (2,1), (-2,1), (0,-1).\eeqq
By rounding off these vertices, $\widetilde{T}$ can be made smooth; after that we denote it by $T$.

\begin{figure}[h!]
\centering
\def\svgwidth{380pt}
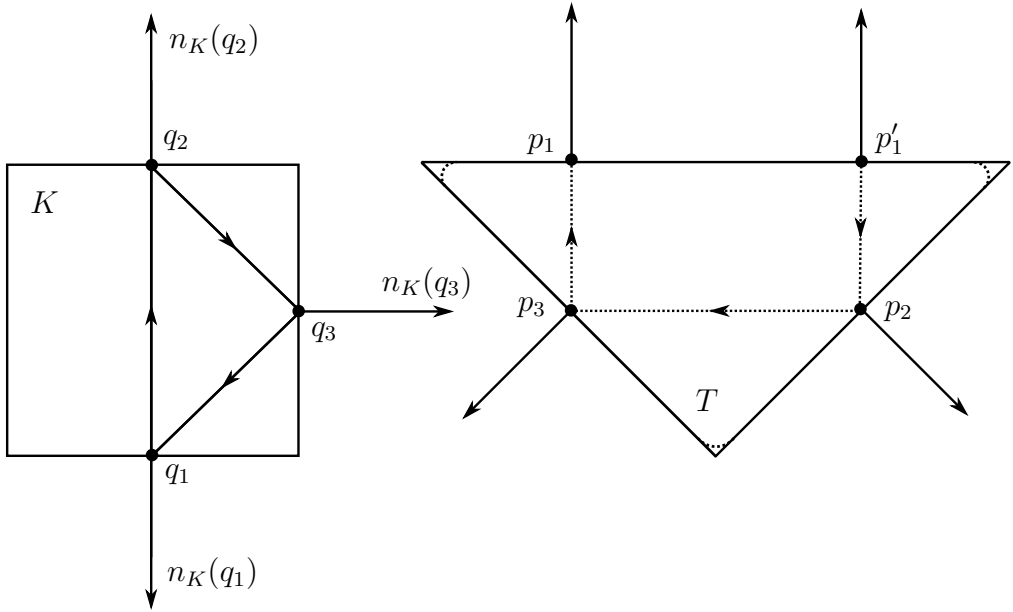
\caption[Example D]{Example D: $T\subset\R^2$ is smooth (the vertices are rounded off). $q$ is a closed weak $(K,T)$-Minkowski billiard trajectory, but the convex cone spanned by $n_K(q_1),n_K(q_2),n_K(q_3)$ is not a linear subspace of $\R^2$.}
\label{img:normalvectorspanning}
\end{figure}

Then, we can check that $q=(q_1,q_2,q_3)$ with
\beqq q_1=(0,-1),\; q_2=(0,1),\; q_3=(1,0)\eeqq
is a closed weak $(K,T)$-Minkowski billiard trajectory.

Indeed, let $H_1,H_2,H_3$ be the uniquely determined $K$-supporting hyperplanes through $q_1,q_2,q_3$. We show that $q_j$ minimizes
\beqq \Sigma_j(q_j^*)=\mu_{T^\circ}(q_j^*-q_{j-1})+\mu_{T^\circ}(q_{j+1}-q_j^*)\eeqq
over all $q_j^*\in H_j$ for all $j\in\{1,2,3\}$. We have
\beqq \Sigma_1(q_1)=\langle q_1-q_3,p_3\rangle + \langle q_2-q_1,p_1\rangle\eeqq
for
\beqq p_1=(-1,1) \; \text{ and } \; p_3=(-1,0).\eeqq
Since
\beqq \langle q_1-q_1^*,p_1-p_3\rangle = 0 \; \text{ for all }q_1^*\in H_1,\eeqq
we conclude for any $q_1^*\in H_1$
\allowdisplaybreaks{\begin{align*}
\Sigma_1(q_1)&= \langle q_1-q_3,p_3\rangle + \langle q_2-q_1,p_1\rangle + \langle q_1-q_1^*,p_1-p_3\rangle\\
&=\langle q_1^*-q_3,p_3\rangle + \langle q_2-q_1^*,p_1\rangle\\
&=\langle q_1^*-q_3,p_3^*\rangle + \langle q_2-q_1^*,p_1^*\rangle + \langle q_1^*-q_3,p_3-p_3^*\rangle + \langle q_2-q_1^*,p_1-p_1^*\rangle,
\end{align*}}%
where $p_1^*,p_3^*\in\partial T$ (possibly not uniquely determined) fulfill
\beqq q_1^*-q_3\in N_T(p_3^*) \; \text{ and }\; q_2-q_1^*\in N_T(p_1^*).\eeqq
From the convexity of $T$, it follows
\beqq \langle q_1^*-q_3,p_3-p_3^*\rangle \leq 0\; \text{ and }\; \langle q_2-q_1^*,p_1-p_1^*\rangle \leq 0\eeqq
and therefore
\beqq \Sigma_1(q_1)\leq \langle q_1^*-q_3,p_3^*\rangle + \langle q_2-q_1^*,p_1^*\rangle = \Sigma_1(q_1^*).\eeqq
Consequently, $q_1$ minimizes $\Sigma_1(q_1^*)$ over all $q_1^*\in H_1$. The same argumentation yields
\beqq \Sigma_2(q_2)=\langle q_2-q_1,p_1'\rangle + \langle q_3-q_2,p_2\rangle \leq \Sigma_2(q_2^*)\eeqq
for all $q_2^*\in H_2$, where
\beqq p_1'=(1,1) \; \text{ and } \; p_2=(1,0),\eeqq
and also
\beqq \Sigma_3(q_3)=\langle q_3-q_2,p_2\rangle + \langle q_1-q_3,p_3\rangle \leq \Sigma_3(q_3^*)\eeqq
for all $q_3^*\in H_3$.

We note that the convex cone spanned by the outer unit vectors
\beqq n_K(q_1)=(0,-1),\; n_K(q_2)=(0,1),\; n_K(q_3)=(1,0)\eeqq
which are normal to $H_1,H_2,H_3$ is not a linear subspace of $\R^2$.

Summarized, this example shows that, without requiring $T$ to be strictly convex, it can happen that the convex cone of the outer unit vectors normal to the $K$-supporting hyperplanes which correspond via the weak Minkowski billiard reflection rule to a closed weak $(K,T)$-Minkowski billiard trajectory does not is the whole $\R^n$.\qed

\underline{Example E}: Let $K\subset\R^2$ be the trapezoid given by the vertices
\beqq (1,-1), (4,2), (-4,2), (-1,-1)\eeqq
and $T\subset\R^2$ the triangle given by the vertices
\beqq \left(\frac{1}{2},2\right), \left(-\frac{1}{2},0\right), \left(\frac{1}{2},-2\right).\eeqq

\begin{figure}[h!]
\centering
\def\svgwidth{400pt}
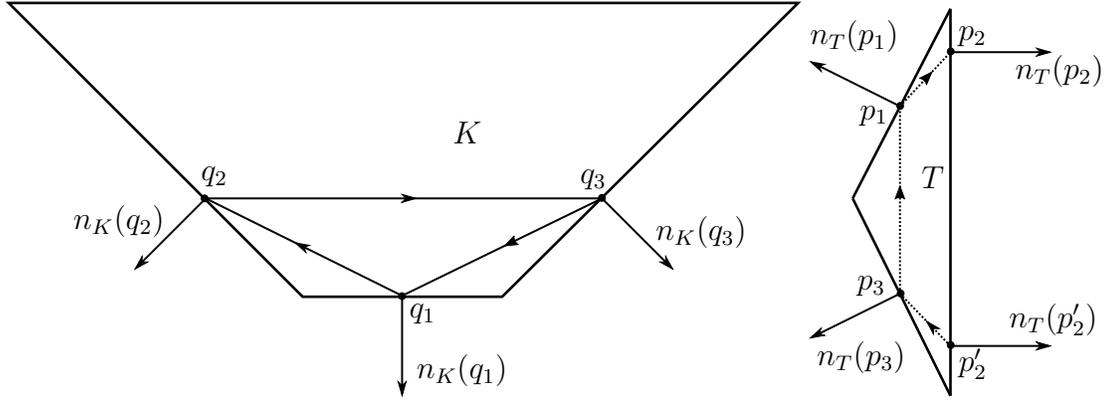
\caption[Example E]{Example E: $T$ is not strictly convex. $q=(q_1,q_2,q_3)$ is a closed $(K,T)$-Minkowski billiard trajectory which can be translated into the interior of $K$. Accordingly, the convex hull of $n_K(q_1),n_K(q_2),n_K(q_3)$ does not include the origin.}
\label{img:notranslation}
\end{figure}

Then, $q=(q_1,q_2,q_3)$ with
\beqq q_1=(0,-1), \; q_2=(-2,0),\; q_3=(2,0)\eeqq
is a closed weak $(K,T)$-Minkowski billiard trajectory.

Indeed, we denote the uniquely determined $K$-supporting hyperplanes through $q_1,q_2,q_3$ by $H_1,H_2,H_3$ and show that $q_j$ minimizes
\beqq \Sigma_j(q_j^*)=\mu_{T^\circ}(q_j^*-q_{j-1})+\mu_{T^\circ}(q_{j+1}-q_j^*)\eeqq
over all $q_j^*\in H_j$ for all $j\in\{1,2,3\}$. We have
\beqq \Sigma_2(q_2)=\langle q_2-q_1,p_1\rangle + \langle q_3-q_2,p_2\rangle\eeqq
for
\beqq p_1=(0,1) \; \text{ and } \; p_2=\left(\frac{1}{2},\frac{3}{2}\right).\eeqq
Since
\beqq \langle q_2-q_2^*,p_2-p_1\rangle =0\text{ for all }q_2^*\in H_2,\eeqq
we conclude for any $q_2^*\in H_2$ that
\allowdisplaybreaks{\begin{align*}
\Sigma_2(q_2)&=\langle q_2-q_1,p_1\rangle + \langle q_3-q_2,p_2\rangle + \langle q_2-q_2^*,p_2-p_1\rangle\\
&= \langle q_2^*-q_1,p_1\rangle + \langle q_3-q_2^*,p_2\rangle\\
&= \langle q_2^*-q_1,p_1^*\rangle + \langle q_3-q_2^*,p_2^*\rangle + \langle q_2^*-q_1,p_1 - p_1^*\rangle + \langle q_3-q_2^*,p_2 - p_2^*\rangle,
\end{align*}}%
where $p_1^*,p_2^*\in\partial T$ (possibly not uniquely determined) fulfill
\beqq q_2^*-q_1\in N_T(p_1^*)\;\text{ and }\; q_3-q_2^*\in N_T(p_2^*).\eeqq
From the convexity of $T$, it follows
\beqq  \langle q_2^*-q_1,p_1 - p_1^*\rangle \leq 0\;\text{ and }\; \langle q_3-q_2^*,p_2 - p_2^*\rangle \leq 0\eeqq
and therefore
\beqq \Sigma_2(q_2)\leq \langle q_2^*-q_1,p_1^*\rangle + \langle q_3-q_2^*,p_2^*\rangle = \Sigma_2(q_2^*).\eeqq
Consequently, $q_2$ minimizes $\Sigma_2(q_2^*)$ over all $q_2^*\in H_2$. The same argumentation yields
\beqq \Sigma_3(q_3)=\langle q_3-q_2,p_2'\rangle + \langle q_1-q_3,p_3\rangle\leq \Sigma_3(q_3^*)\eeqq
for all $q_3^*\in H_3$, where
\beqq p_2'=\left(\frac{1}{2},-\frac{3}{2}\right) \; \text{ and } \; p_3=(0,-1),\eeqq
and also
\beqq \Sigma_1(q_1)=\langle q_1-q_3,p_3\rangle + \langle q_2-q_1,p_1\rangle \leq \Sigma_1(q_1^*)\eeqq
for all $q_1^*\in H_1$.

Finally, firstly, we note that $0$ is not within the convex hull of the unit normal vectors $n_K(q_1),n_K(q_2),n_K(q_3)$, and, secondly, that $q$ can be translated into the interior of $K$; for example
\beqq q+\left(0,\frac{1}{2}\right)\in \mathring{K}.\eeqq

Summarized, this example shows that, without requiring $T$ to be strictly convex and smooth, it can happen that the convex hull of the outer unit vectors normal to the $K$-supporting hyperplanes which correspond via the weak Minkowski billiard reflection rule to a closed weak $(K,T)$-Minkowski billiard trajectory does not include the origin. Furthermore, it shows that in the situation of non-strictly convex $T$, it can happen that closed weak $(K,T)$-Minkowski billiard trajectories can be translated into the interior of $K$.\qed

\underline{Example F}: Let $K\subset\R^2$ be the convex polytope given by the vertices
\beqq \left(\frac{1}{2},0\right),(0,1),(-2,1),(-2,-1),(0,-1)\eeqq
and $\widetilde{T}\subset\R^2$ the rhombus given by the vertices
\beqq (1,0),(0,1),(-1,0),(0,-1).\eeqq
Then, we can find a strictly convex body $T\subset\R^2$ satisfying
\beqq \widetilde{T}\subseteq T\; \text{ and } \; (1,0),(0,1),(-1,0),(0,-1)\in\partial T\eeqq
such that $N_T((0,1))$ equals the convex cone spanned by the vectors $(1,2)$ and $(-1,2)$ and
\beqq N_T((0,-1))=-N_T((0,1)).\eeqq
Then, $T$ is strictly convex, but not smooth.

\begin{figure}[h!]
\centering
\def\svgwidth{360pt}
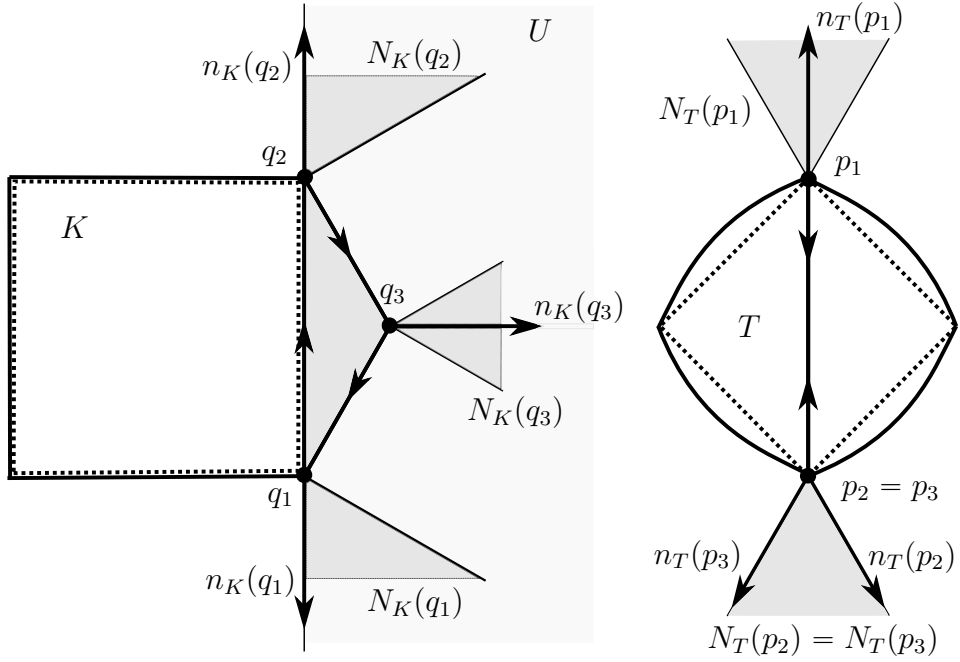
\caption[Example F]{Example F: In the figure, the dashed lines belong to $\widetilde{K}$ and $\widetilde{T}$, respectively. $q=(q_1,q_2,q_3)$ is an $\ell_T$-minimizing closed $(K,T)$-Minkowski billiard trajectory, where $p=(p_1,p_2,p_3)$ is its dual billiard trajectory in $T$. We neither have that $U$ (the convex cone spanned by $n_K(q_1),n_K(q_2),n_K(q_3)$) is a linear subspace of $\R^2$ with $\dim U = 2$, nor that that $\dim (N_K(q_j)\cap U)=1$ for all $j\in\{1,2,3\}$.}
\label{img:ExampleF}
\end{figure}

We claim that $q=(q_1,q_2,q_3)$ with
\beqq q_1=(0,-1),\;q_2=(0,1),\;q_3=\left(\frac{1}{2},0\right)\eeqq
is an $\ell_T$-minimizing closed $(K,T)$-Minkowski billiard trajectory. Then, if we denote by $U$ the convex cone spanned by the outer unit vectors
\beqq n_K(q_1),n_K(q_2),n_K(q_3)\eeqq
which are normal to $H_1,H_2,H_3$, then
\beqq U=\R_{\geq 0}\times \R\eeqq
and
\beqq \dim (N_K(q_j)\cap U)>1 \quad \forall j\in\{1,2,3\}.\eeqq

So, let us prove that $q$ is an $\ell_T$-minimizing closed $(K,T)$-Minkowski billiard trajectory. First, $q$ is a closed $(K,T)$-Minkowski billiard trajectory by checking that $p=(p_1,p_2,p_3)$ with
\beqq p_1=(0,1),\;p_2=p_3=(0,-1)\eeqq
is its dual billiard trajectory in $T$. One has
\beqq \ell_T(q)=\langle q_2-q_1,p_1\rangle + \langle q_3-q_2,p_2\rangle + \langle q_1-q_3,p_3\rangle = 2+1+1=4.\eeqq
Now, let $\widetilde{K}$ be the convex polytope defined by the vertices
\beqq (0,-1), (0,1), (-2,1), (-2,-1).\eeqq 
We have
\beqq \widetilde{K}\subseteq K \; \text{ and } \; \widetilde{T} \subseteq T.\eeqq
Therefore, we conclude that the $\ell_{\widetilde{T}}$-length of the $\ell_{\widetilde{T}}$-minimizing closed $(\widetilde{K},\widetilde{T})$-Minkowski billiard trajectories is less or equal than the $\ell_T$-length of the $\ell_T$-minimizing closed $(K,T)$-Minkowski billiard trajectories. One easily checks that the former is $4$. This implies that the latter cannot be less than $4$. Since the $\ell_T$-length of $q$ is $4$, this implies that $q$, in fact, is an $\ell_T$-minimizing closed $(K,T)$-Minkowski billiard trajectory.

Summarized, this example shows that, without requiring $T$ to be strictly convex, it can happen that there is an $\ell_T$-minimizing closed $(K,T)$-Minkowski billiard trajectory which violates all the statements made in Theorem \ref{Thm:RegularityResult}.\qed

\underline{Example G}: Let $K\subset\R^2$ be the triangle given by the vertices
\beqq (1,1),(-1,1),(0,-1)\eeqq
and $T\subset\R^2$ the rectangle given by the vertices
\beqq (1,-2),(1,2),(-1,2),(-1,-2).\eeqq

\begin{figure}[h!]
\centering
\def\svgwidth{350pt}
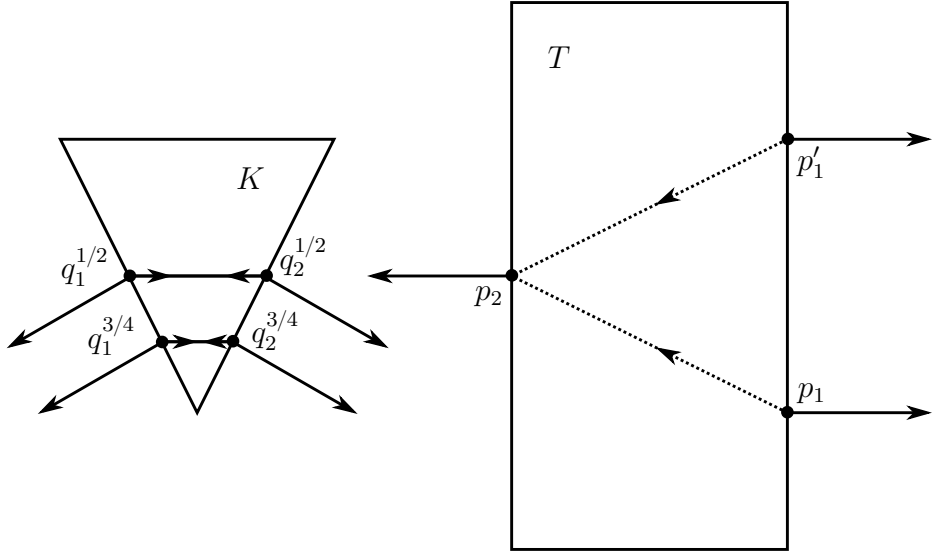
\caption[Example G]{Example G: For $a\rightarrow 1$ the closed weak $(K,T)$-Minkowski billiard trajectory $q^a=(q_1^a,q_2^a)$ $d_H$-converges to the point $(0,-1)$.}
\label{img:ExampleG}
\end{figure}

Then, we claim that for $a\in[0,1)$, $q^a=(q_1^a,q_2^a)$ with
\beqq q_1^a=(-1+a,1-2a)\;\text{ and }\; q_2^a=(1-a,1-2a)\eeqq
is a closed weak $(K,T)$-Minkowski billiard trajectory which fulfills the weak Minkowski billiard reflection rule with respect to the uniquely determined $K$-supporting hyperplanes $H_1,H_2$ through $q_1^a,q_2^a$.

Indeed, for $a\in[0,1)$ we show that $q_j^a$ minimizes
\beqq \Sigma_j(q_j^*)=\mu_{T^\circ}(q_j^*-q_{j-1}^a)+\mu_{T^\circ}(q_{j+1}^a-q_j^*)\eeqq
(note: $q_{j+1}^a=q_{j-1}^a$) over all $q_j^*\in H_j$ for $j\in\{1,2\}$. The following holds for all $a\in[0,1)$: We have
\beqq \Sigma_2(q_2^a)=\langle q_2^a-q_1^a,p_1\rangle + \langle q_1^a-q_2^a,p_2\rangle\eeqq
for
\beqq p_1=(1,-1) \; \text{ and } \; p_2=(-1,0).\eeqq
Since
\beqq \langle q_2^*-q_2^a,p_2-p_1\rangle =0 \text{ for all }q_2^*\in H_2,\eeqq
we conclude for any $q_2^*\in H_2$ that
\allowdisplaybreaks{\begin{align*}
\Sigma_2(q_2^a)&=\langle q_2^a-q_1^a,p_1\rangle + \langle q_1^a-q_2^a,p_2\rangle + \langle q_2^*-q_2^a,p_2-p_1\rangle\\
&=\langle q_2^*-q_1^a,p_1\rangle + \langle q_1^a-q_2^*,p_2\rangle\\
&=\langle q_2^*-q_1^a,p_1^*\rangle + \langle q_1^a-q_2^*,p_2^*\rangle + \langle q_2^*-q_1^a,p_1-p_1^*\rangle + \langle q_1^a-q_2^*,p_2-p_2^*\rangle,
\end{align*}}%
where $p_1^*,p_2^*\in\partial T$ (possibly not uniquely determined) fulfill
\beqq q_2^*-q_1^a\in N_T(p_1^*)\;\text{ and }\; q_1^a-q_2^*\in N_T(p_2^*).\eeqq
From the convexity of $T$, it follows
\beqq \langle q_2^*-q_1^a,p_1-p_1^*\rangle \leq 0\;\text{ and }\; \langle q_1^a-q_2^*,p_2-p_2^*\rangle \leq 0\eeqq
and therefore
\beqq \Sigma_2(q_2^a)\leq \langle q_2^*-q_1^a,p_1^*\rangle + \langle q_1^a-q_2^*,p_2^*\rangle = \Sigma_2(q_2^*).\eeqq
Consequently, $q_2^a$ minimizes $\Sigma_2(q_2^*)$ over all $q_2^*\in H_2$. The same argumentation yields
\beqq \Sigma_1(q_1^a)=\langle q_1^a-q_2^a,p_2\rangle + \langle q_2^a-q_1^a,p_1'\rangle \leq \Sigma_1(q_1^*)\eeqq
for all $q_1^*\in H_1$, where $p_1'=(1,1)$.

We have
\beqq \ell_T(q^a)=\langle q_2^a-q_1^a,p_1\rangle + \langle q_1^a-q_2^a,p_2\rangle = 2-2a + 2-2a = 4-4a\eeqq
which for $a\rightarrow 1$ goes to $0$. Therefore, there is no $\ell_T$-minimizing closed $(K,T)$-Minkowski billiard trajectory (cf.\;Footnote \ref{foot:polygonalline}).

We note that $T$ can be made smooth without loosing the above mentioned properties.

Summarized, this example shows that there are configurations $(K,T)$ for which $T$ is not strictly convex and there is no $\ell_T$-minimizing closed weak $(K,T)$-Minkowski billiard trajectory.

\section[Constructing Minkowski billiard trajectories]{Constructing shortest Minkowski billiard trajectories on convex polytopes}\label{Sec:Construction}

\subsection{General construction in two dimensions}\label{Subsec:GeneralConstruction}

In this first subsection, we describe the general construction of $\ell_T$-minimizing closed $(K,T)$-Minkowski billiard trajectories for the case of a convex polytope $K\subset\R^2$ and a strictly convex and smooth body $T\subset\R^2$. For determining the $\ell_T$-minimizing closed $(K,T)$-Minkowski billiard trajectories, we use Corollary \ref{Cor:RegularityResult}, i.e., the $\ell_T$-minimizing closed $(K,T)$-Minkowski billiard trajectories have two or three bouncing points, where in the latter case the billiard trajectories are regular.

In \cite{AlkoumiSchlenk2014}, the algorithm for finding closed $(K,T)$-Minkowski billiard trajectories with two bouncing points has already been described. For details concerning the implementation, we refer to Section \ref{Subsec:Implementation}. In \cite{AlkoumiSchlenk2014}, it was stated as open problem to find an algorithm for determining closed regular $(K,T)$-Minkowski billiard trajectories with three bouncing points. While there, for $T$ equals the Euclidean unit ball in $\R^2$, they could use the uniqueness of Fagnano triangles in acute triangles in order to find the closed regular Euclidean billiard trajectories on $K$, in the Minkowski/Finsler setting one has to find a different approach, since there are no obvious analogues of the Fagnano triangles at first. Now, this will be the task of the remainder of this subsection.

\begin{figure}[h!]
\centering
\def\svgwidth{415pt}
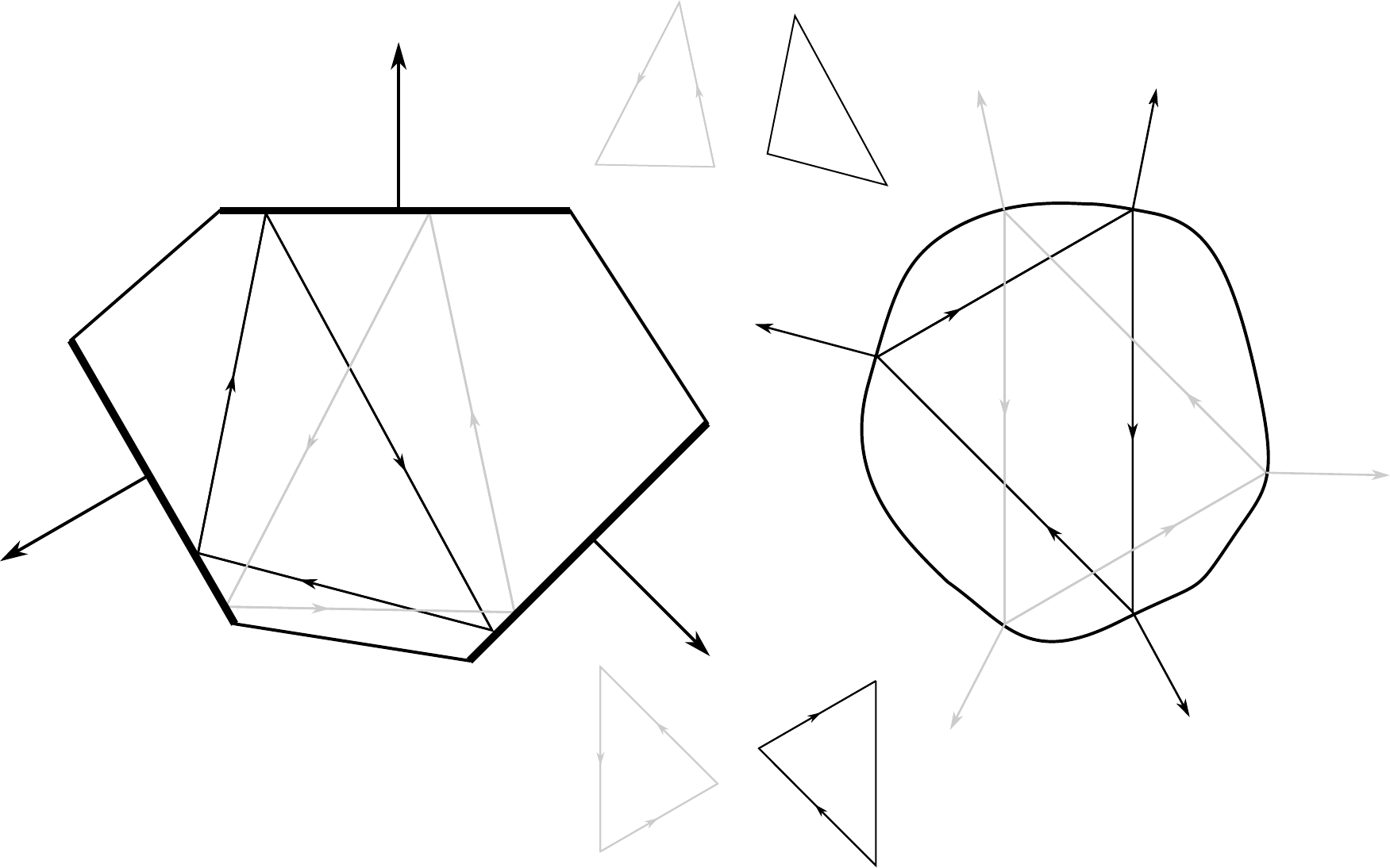
\caption{Illustration of the construction of closed regular $(K,T)$-Minkowski billiard trajectories with three bouncing points.}
\label{img:Algorithm}
\end{figure}

We do the following (cf.\;Figure \ref{img:Algorithm}):
\begin{itemize}
\item[(a)] Choose $3$ facets $F_1,F_2,F_3$ of $K$ (considering their order) such that the convex cone spanned by the associated outer normal unit vectors $n_{F_1},n_{F_2},n_{F_3}$ is $\R^2$.
\item[(b)] Construct the uniquely determined (up to scaling and translation) closed polygonal curve $(\gamma_1,\gamma_2,\gamma_3)$ with $\gamma_{i+1}-\gamma_i$ given by a negative multiple of $n_{F_i}$ for all $i\in \{1,2,3\}$.
\item[(c)] Find the uniquely determined $\lambda >0$ and $c\in \R^n$ such that
\beqq \lambda\{\gamma_1,\gamma_2,\gamma_3\}+c\subset \partial T.\eeqq
\item[(d)] Let $n_i$ be the outer normal unit vector ar $\partial T$ in the point
\beqq \lambda\gamma_i+c\eeqq
for all $i\in\{1,2,3\}$. If the convex cone spanned by $n_1,n_2,n_3$ is $\R^2$, then construct the uniquely determined (up to scaling and translation) closed polygonal curve $(\xi_1,\xi_2,\xi_3)$ with $\xi_{i+1}-\xi_i$ given by a positive multiple of $n_{i+1}$ for all $i\in \{1,2,3\}$. Otherwise: If possible: Go back to (a) and start with a choice not yet made. Otherwise: End.
\item[(e)] If possible: Find $\mu >0$ and $e\in\R^n$ such that
\beqq \mu\xi_i+e\in \mathring{F}_i\quad \forall i\in\{1,2,3\}.\eeqq
Otherwise: If possible: Go back to (a) and start with a choice not yet made. Otherwise: End.
\item[(f)] Define a closed polygonal curve
\beqq q=(q_1,q_2,q_3)\eeqq
by
\beqq q_i:=\mu\xi_i+e \quad \forall i\in\{1,2,3\}.\eeqq
By construction: $q$ is a maximally spanning, closed, regular $(K,T)$-Minkowski billiard trajectory with three bouncing points and with closed dual billiard trajectory
\beqq p=(p_1,p_2,p_3)\eeqq
given by
\beqq p_i:=\lambda \gamma_{i+1}+c \quad  \forall i\in\{1,2,3\}.\eeqq
Add $q$ to $B_3(K,T)$.
\item[(g)] If possible: Go back to (a) and start with a choice not yet made. Otherwise: End.
\end{itemize}

Finally, the set $B_3(K,T)$ contains all closed regular $(K,T)$-Minkowski billiard trajectories with three bouncing points whose $\ell_T$-length can be easily calculated:
\beqq \ell_T(q)=\sum_{j=1}^3\mu_{T^\circ}(q_{j+1}-q_j)=\sum_{j=1}^3\langle q_{j+1}-q_j,p_j\rangle.\eeqq

Let us now turn to the explanation of the individuel steps, while for the detailed justification, we refer to Section \ref{Subsec:Implementation}:

Ad (a): If there is a closed regular $(K,T)$-Minkowski billiard trajectory with three bouncing points, then we know from Proposition \ref{Prop:normalvectorsspanning} that the bouncing points lie in the interiors of three different facets of $K$ whose associated outer normal unit vectors span $\R^2$.

Ad (b): Since the convex cone spanned by $n_{F_1},n_{F_2},n_{F_3}$ is $\R^2$, solving a system of linear equations yields a uniquely determined (up to scaling and translation) $3$-tuple
\beqq (\gamma_1,\gamma_2,\gamma_3)\in(\R^2)^3\eeqq
and a uniquely determined (up to scaling--depending on the factor by which $(\gamma_1,\gamma_2,\gamma_3)$ will be scaled) $3$-tuple
\beqq (\alpha_1,\alpha_2,\alpha_3)\in (\R_{<0})^3\eeqq
fulfilling
\beqq \gamma_{i+1}-\gamma_i=\alpha_i n_{F_i} \quad i\in\{1,2,3\}.\eeqq
We understand the $3$-tuple $(\gamma_1,\gamma_2,\gamma_3)$ as a closed polygonal curve.

Ad (c): Because of the strict convexity of $T$, there is a unique combination
\beqq (\lambda,c)\in \R_{>0}\times \R^2\eeqq
such that
\beqq \lambda\{\gamma_1,\gamma_2,\gamma_3\}+c\subset \partial T.\eeqq

Ad (d): If the convex cone spanned by the unit vectors $n_1,n_2,n_3$ is $\R^2$, then, as in step (b), solving a system of linear equations yields a uniquely determined (up to scaling und translation) $3$-tuple
\beqq (\xi_1,\xi_2,\xi_3)\in(\R^2)^3\eeqq
and a uniquely determined (up to scaling--depending on the factor by which $(\xi_1,\xi_2,\xi_3)$ will be scaled) $3$-tuple
\beqq (\beta_1,\beta_2,\beta_3)\in (\R_{>0})^3\eeqq
fulfilling
\beqq \xi_{i+1}-\xi_i=\beta_{i+1} n_{i+1} \quad \forall i\in\{1,2,3\}\eeqq
We understand the $3$-tuple $(\xi_1,\xi_2,\xi_3)$ as a closed polygonal curve.

Ad (e) $\&$ (f): There is at most one combination
\beqq (\mu,e)\in \R_{>0}\times \R^2\eeqq
such that
\beqq \mu\{\xi_1,\xi_2,\xi_3\}+e\subset \partial K.\eeqq
By checking whether
\beqq \mu \xi_i +e \in \mathring{F}_i \quad \forall i\in\{1,2,3\},\eeqq
we make sure that the closed polygonal curve
\beqq q=(q_1,q_2,q_3)\eeqq
defined by
\beqq q_i:=\mu \xi_i +e\quad \forall i\in\{1,2,3\}\eeqq
has its vertices in the interiors of the facets $F_1,F_2,F_3$. $q$ is a closed $(K,T)$-Minkowski billiard trajectory, where
\beqq p=(p_1,p_2,p_3)\eeqq
with
\beqq p_i:=\lambda \gamma_{i+1}+c\quad \forall i\in\{1,2,3\}\eeqq
is its closed dual billiard trajectory on $T$. Indeed, we define
\beqq \lambda_i:=\mu\beta_i >0 \; \text{ and } \, \mu_i:=-\lambda \alpha_i >0\eeqq
and notice that the pair $(q,p)$ fulfills \eqref{eq:System}:
\beqq \begin{cases} q_{i+1}-q_i=(\mu \xi_{i+1}+e)-(\mu\xi_i +e)=\mu(\xi_{i+1}-\xi_i)=\mu\beta_{i+1} n_{i+1}=\lambda_{i+1} n_{i+1}\in N_T(p_i),\\
p_{i+1}-p_i=(\lambda\gamma_{i+2}+c)-(\lambda\gamma_{i+1}+c)=\lambda (\gamma_{i+2}-\gamma_{i+1})=-\mu_{i+1}n_{F_{i+1}}\in -N_K(q_{i+1}).\end{cases}
\eeqq

\subsection{A note concerning the general construction for higher dimensions}\label{Subsec:Generalconstructionhigherdimensions}

Let $K\subset\R^n$ be a convex polytope and $T\subset\R^n$ a strictly convex and smooth body. We know from Theorems \ref{Thm:RegularityResult} and \ref{Thm:maximxallyspanning} that there is always an $\ell_T$-minimizing closed $(K,T)$-Minkowski billiard trajectory which is maximally spanning, has at most $n+1$ bouncing points and whose corresponding outer unit normal vectors span a cone which has the same dimension as the inclusion minimal section containing this trajectory.

Instead of that the Euclidean unit ball is replaced by $T$ and one has to take into account that the linear subspaces underlying the inclusion minimal affine sections containing relevant Minkowski billiard trajectories can differ from the convex cone spanned by the corresponding outer unit normal vectors, these are the same preconditions as within the algorithms for the Euclidean setting. The necessary aspects which one has to consider for the adjustment to the Minkowski setting are indicated in Section \ref{Subsec:GeneralConstruction} for two dimensions.

We leave the detailed execution of these adjustments to further research. 

\subsection[Efficiency and used methods within the implementation of the two-dimensional case]{Efficiency and used methods within the implementation}\label{Subsec:Implementation}

We now turn our attention to the implementation of the algorithm for two dimensions which is described in Section \ref{Subsec:GeneralConstruction}.

Even though we focus on the case $n=2$, we state some of the results in this subsection for arbitrary $n$ if they hold in any dimension. In Section \ref{Subsec:GeneralConstruction}, we applied Corollary \ref{Cor:RegularityResult} and therefore required $T$ to be strictly convex and smooth. Implementing such a set can be a difficult problem because we can only make finitely many inputs. Therefore, we assume that both $K$ and $T$ are convex polytopes (in particular neither strictly convex nor smooth) in this subsection. In the following this has to be justified.

We proved in Theorem \ref{Thm:RegularityResult} for the case when $T$ is assumed to be strictly convex and smooth and when considering the closed $(K,T)$-Minkowski billiard trajectories with $n+1$ bouncing points, that, when searching for length minimizers, it is enough to just concentrate on the $\ell_T$-minimizing closed $(K,T)$-Minkowski billiard trajectories with $n+1$ bouncing points which are regular, i.e., whose normal cones in the bouncing points are one-dimensional. However, in \cite[Section 4.3.2]{Krupp2021}, it is shown that in the case when $T$ is assumed to be a convex polytope, then, the bouncing points of an $\ell_T$-minimizing closed $(K,T)$-Minkowski billiard trajectory may be in nonsmooth boundary points of $K$, but one can assume that the normal cones appearing in the system \eqref{eq:System} can be replaced by the rays which are the one-dimensional normal cones of the neighbouring facets.\footnote{This can be proved by approximating the convex polytope $T$ by a sequence of strictly convex and smooth bodies, using a line of argumentation which is similar to the one appearing in the proof of \cite[Theorem 2.1]{Rudolf2022}.} Therefore, when we look for a boundary point in some facet of $K$, we allow it to lie in the boundary of this facet, but we only consider the normal cone for some point in the relative interior of this facet.

Before we analyze the algorithm in greater detail, we argue that it is sufficient to compute finitely many pairs $(q,p)$ of closed polygonal curves fulfilling \eqref{eq:System} to find one, where $q$ is $\ell_T$-minimizing. More precisely, we show that if $(q, p)$ and $(q', p')$ are pairs fulfilling \eqref{eq:System} and if their vertices lie on the same faces of $K$ (resp. $T$), then $q$ and $q'$ have the same $\ell_T$-length.

\bprop[Theorem 4.3.6 in \cite{Krupp2021}]\label{Prop:samefacessamelength}
Let $K,T \subseteq \mathbb{R}^n$ be convex polytopes, where $F_1,...,F_m$ and $G_1,...,G_m$ are the faces of $K$ and $T$, respectively. Further, let
\beqq q=(q_1,...,q_m) \; \text{ and } \; q'=( q_1',...,q_m')\eeqq
be closed polygonal curves with vertices on $\partial K$. Assume, there are closed polygonal curves
\beqq p=(p_1,...,p_m)\; \text{ and } \; p'= (p_1',...,p_m')\eeqq
with vertices on $\partial T$ such that $(q,p)$ and $(q',p')$ fulfill \eqref{eq:System}. Further, assume for each $j \in \{1,...,m\}$ that
\beqq q_j,q'_j \in \textup{relint}(F_j) \; \text{ and } \; p_j,p'_j \in \textup{relint}(G_j)\eeqq
(unless $F_j$ is a vertex, in which case we assume $q_j,q'_j \in F_j$ instead. The same applies for $G_j$). Then
\beqq \ell_T(q) = \ell_T(q').\eeqq
\eprop

\bpf
We start the proof by stating a simple fact. If $F$ is a face of a convex polytope $P$ with $\textup{dim}(F) \geq 1$ and $y_1, y_2 \in F$, then
\begin{align}\label{Pf:samefacessamelength}
\langle y_1 - y_2 , v \rangle = 0\, , \ \forall v \in N_P(z)
\end{align}
holds for any $z \in \textup{relint}(F)$. To see this, consider the affine hull of $F$ and shift it, such that it becomes a linear space. Then the vector $y_1-y_2$ is an element of this space and $N_P(z)$ with $z \in \textup{relint}(F)$ is contained in the corresponding orthogonal space. Note that \eqref{Pf:samefacessamelength} also holds if $F$ is a vertex in which case we have $\textup{relint}(F) = F$. Then $y_1, y_2 \in F$ implies $y_1 = y_2$ and the statement follows immediately. Now recall \eqref{eq:System}:
\begin{align*}
 q_{j+1}-q_j \in N_T(p_j)\, , \hspace*{1cm}&  q'_{j+1}-q'_j \in N_T(p'_j)\, , \\ p_{j+1} - p_j \in -N_K(q_{j+1})\, , \hspace*{1cm}& p'_{j+1} - p'_j \in -N_K(q'_{j+1}).
\end{align*}
The following calculation completes the proof:
{\allowdisplaybreaks \begin{align*}
\ell_T(q') & = \sum\limits_{j=1}^m \langle  q'_{j+1}-q'_j , p'_j \rangle\\
& = \sum\limits_{j=1}^m \langle  q_{j+1}-q_j, p'_j \rangle + \sum\limits_{j=1}^m \langle q_j - q_j', p'_j \rangle - \sum\limits_{j=1}^m \langle q_{j+1} - q'_{j+1}, p'_j \rangle\\
& =  \sum\limits_{j=1}^m \langle q_{j+1} - q_{j}, p'_j \rangle + \sum\limits_{j=1}^m \langle q_j - q'_j, p'_j \rangle - \sum\limits_{j=1}^m \langle q_{j} - q'_{j}, p'_{j-1} \rangle\\
& =  \sum\limits_{j=1}^m \langle q_{j+1} - q_{j}, p'_j \rangle + \sum\limits_{j=1}^m \langle q_j - q'_j, p'_j - p'_{j-1} \rangle\\
& = \sum\limits_{j=1}^m \langle q_{j+1} - q_{j}, p'_j \rangle\\
& = \sum\limits_{j=1}^m \langle q_{j+1} - q_{j}, p_j \rangle + \sum\limits_{j=1}^m \langle q_{j+1} - q_{j}, p'_j - p_j \rangle \\
& = \sum\limits_{j=1}^m \langle q_{j+1} - q_{j}, p_j \rangle \\
&= \ell_T(q).
\end{align*}
Note that we used \eqref{Pf:samefacessamelength} to show that
\begin{align*}
\langle q_j - q'_j, p'_{j-1} - p'_{j} \rangle = 0 = \langle q_{j+1} - q_{j}, p'_j - p_j \rangle.
\end{align*}}%
\epf

For the remainder of this section, we fix $n=2$. We proceed with the case $m=2$. In other words, for each choice of faces $F_1,F_2$ of $K$ and for each choice of faces $G_1,G_2$ of $T$, we compute closed polygonal curves
\beqq q = (q_1,q_1) \; \text{ and } \; p = (p_1,p_2)\eeqq
fulfilling \eqref{eq:System} such that $q_j \in F_j$ and $p_j \in G_j$ for $j\in \{1,2\}$. Note that Proposition \ref{Prop:samefacessamelength} suggests that we ask for
\beqq q_j \in \textup{relint}(F_j) \; \text{ and } \; p_j \in \textup{relint}(G_j).\eeqq
Instead, for the sake of simplicity, we allow $q_j \in \partial F_j$ and replace $N_K(q_j)$ with $N_K(z)$ for some $z \in \textup{relint}(F_j)$ (this applies analogously to $p_j \in G_j$). The proof of Proposition \ref{Prop:samefacessamelength} extends directly to this case. 

After considering all choice of $F_1,F_2,G_1,G_2$, we compare the $\ell_T$-length of all found closed polygonal curves. Before starting the calculation, it is beneficial to check whether
\beqq N_K(q_1) \cap -N_K(q_2) \; \text{ and } \; N_T(p_1) \cap -N_T(p_2)\eeqq
are nonempty. The reason for this is, that the existence of a pair $(q, p)$ of closed polygonal curves satisfying \eqref{eq:System} implies:
\begin{align*}
-N_K(q_2) \ni p_2 - p_1 = -(p_1 - p_2) \in N_K(q_1),\\
N_T(p_1) \ni q_2 - q_1 = -(q_1 - q_2) \in -N_T(p_2).
\end{align*}
Note that the normal cones $N_K(q_j),N_T(p_j)$ only depend on the faces $F_j,G_j$. So, in the following, we can assume that these intersections are indeed nonempty. The goal is now to calculate a pair of suitable polygonal curves $(q, p)$ if possible. For this, it is helpful to distinguish whether the faces $F_j,G_j$ are facets (i.e., edges) or vertices. We consider the following cases:
\begin{itemize}
\item[1)] $F_1,F_2,G_1,G_2$ are vertices.
\item[2)] $F_1,F_2$ are vertices and among $G_1,G_2$ there is at least one facet.
\item[3)] Among $F_1,F_2$ as well as among $G_1,G_2$ there is at least one facet.
\end{itemize}
All remaining cases can be covered by switching the roles of $K$ and $T$. 

The first case is easy. If all chosen faces are vertices, the resulting closed polygonal curves are unique and \eqref{eq:System} can be checked directly.

We start the second case by assuming that both $G_1$ and $G_2$ are facets. Therefore, $N_T(p_1)$ and $N_T(p_2)$ are one-dimensional cones. Let 
\begin{align*}
w_j \in N_T(p_j)\setminus \{0\}\, \textup{ for }\,j \in \{1,2\}.
\end{align*}
We can ensure that
\beqq q_2-q_1 \in N_T(p_1)\eeqq
holds by checking whether $w_1$ is a positive multiple of $q_2-q_1$. If this is the case,
\beqq q_1-q_2 \in N_T(p_2)\eeqq
follows directly since we assume that
\beqq N_T(p_1) \cap -N_T(p_2)\eeqq
is nonempty. Alternatively, we can check whether $w_2$ is a positive multiple of $q_1-q_2$ and get
\beqq q_2-q_1 \in N_T(p_1)\eeqq
for free. It remains to solve the following problem:
\begin{align*}
\textup{Find } & p_1,p_2 \textup{ such that:}\\
& p_1 \in G_1\, , \ p_2 \in G_2,\\
& p_2 - p_1 \in -N_K(q_2),\\
& p_1 - p_2 \in -N_K(q_1).
\end{align*}
The constraints can be expressed with linear equations and inequalities. For this, recall the definition of the outer normal cone of a convex set $C$ at $z \in \partial C$:
\begin{align*}
N_C(z) = \{v \colon \langle v, y-z \rangle \leq 0\, ,\ \forall \ y \in C\}.
\end{align*}
If $C$ is a convex polytope, it is sufficient to demand
\beqq \langle v, y-z \rangle \leq 0\eeqq
for every vertex $y$ of $C$. Therefore, membership of $N_C(z)$ can be modeled by finitely many linear inequalities. Altogether, finding suitable points $p_1,p_2$, can be realized by using linear programming techniques. The same approach can be used if either $G_1$ or $G_2$ is a vertex. In this case, the linear program remains unchanged except for the fact that one of the two variable vectors is replaced by a constant vector.

In the third case, we start with the assumption that $F_1,F_2,G_1,G_2$ are facets. Then, all relevant normal cones are one-dimensional and we let 
\begin{align*}
u_j \in N_K(q_j)\setminus \{0\} \; \text{ and } \; w_j \in N_T(p_j)\setminus \{0\} \; \textup{ for } \; j\in \{1,2\}.
\end{align*}
Now, we solve the following problem:
\begin{align*}
\textup{Find } & q_1,q_2,p_1,p_2,\alpha_1,\alpha_2 \textup{ such that:}\\
& q_1 \in F_1\, , \ q_2 \in F_2\, , \ p_1 \in G_1\, , \ p_2 \in G_2,\\
& \alpha_1, \alpha_2 \geq 0,\\
& q_2 - q_1 = \alpha_1 w_1,\\
& p_2 - p_1 = -\alpha_2 u_2.
\end{align*}
Similar to the previous case, this problem is a linear program. Note that the last two constraints suffice to imply \eqref{eq:System} since
\beqq N_T(p_1) \cap -N_T(p_2) \; \text{ and } \; N_K(q_1) \cap -N_K(q_2)\eeqq
are nonempty. If not all chosen faces are facets, for instance if $G_1$ is a vertex, the linear program has to be changed in two ways. First, much like in the second case, the corresponding variable vector, here $p_1$, is replaced by a constant vector. Second, if $G_1$ is a vertex, then the normal cone $N_T(p_1)$ is no longer one-dimensional and the definition of $w_1$ does not make sense any more. However, in this case, $G_2$ is a facet and we replace the constraint
\beqq q_2 - q_1 = \alpha_1 w_1 \; \text{ with } \; q_1 - q_2 = \alpha_1 w_2.\eeqq
We apply this reasoning also when $F_1$ or $F_2$ is not a facet.

There may be multiple ways to choose $(p,q)$ for given faces $F_1,F_2,G_1,G_2$. If this is the case, our algorithm chooses (if possible) $p$ such that $N_K(q_j)$ is one-dimensional for $j\in \{1,2\}$ (or equivalently such that $q_1,q_2$ are not vertices of $K$). This is achieved in the following way. If $F_1$ or $F_2$ is a vertex, the resulting closed polygonal curve $q$ always contains a vertex of $K$. So, we assume both $F_1$ and $F_2$ are facets. If $q_1$ and $q_2$ are smooth points (i.e., lie in the interior of $F_1$ and $F_2$), then there is nothing to do. Otherwise, we denote
\beqq N_K(q_1) = \mathbb{R}_+ u \; \text{ and } \; N_K(q_2) = \mathbb{R}_+ (-u)\eeqq
for some vector $u\in \mathbb{R}^2$. Let $v\neq 0$ be a vector orthogonal to $u$. Moving $q_j$ along the facet $F_j$ can only be done in at most two directions: $v$ or $-v$. If we can move both $q_1$ and $q_2$ in the same direction, we simply translate the closed polygonal curve $q$. If $q_1$ and $q_2$ can only be moved in opposite directions, it is necessary to check whether the normal cones $N_T(p_1)$ and $N_T(p_2)$ allow such movement. If not, it is not possible to find suitable points
\beqq q_1 \in \textup{relint}(F_1) \; \text{ and } \; q_2 \in \textup{relint}(F_2).\eeqq

\begin{figure}[h!]
\centering
\def\svgwidth{300pt}
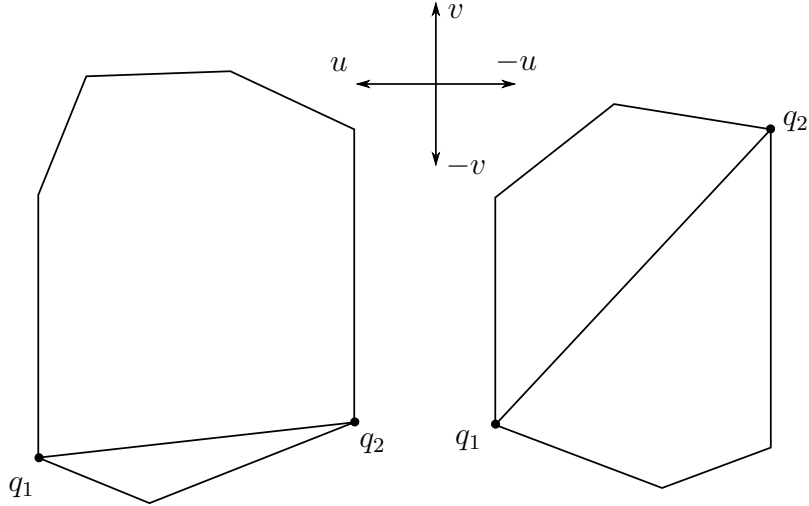
\caption[Illustration of an argument within the preparation of the implementation]{Two closed polygonal curves such that $q_1,q_2$ are vertices of $K$. On the left, we can translate $q_1,q_2$ upwards (in direction $v$). On the right, we need to move $q_1$ upwards and $q_2$ downwards. Whether this is possible depends on $N_T(p_1)$ and $N_T(p_2)$.}
\label{img:Picture 7}
\end{figure}

This concludes the algorithm for the case $m = 2$. Next, we will discuss the algorithm for $m = 3$. We start as described in the previous subsection and choose facets $F_1,F_2,F_3$ of $K$. For each $j \in \{1,2,3\}$, we let $n_{K,j}$ be the outer unit normal vector of $K$ at some point in the relative interior of $F_j$. If possible, we construct a triangle $\Delta$ by only using negative multiples of these three vectors. Here a triangle is the convex hull of three affinely independent points. This can easily be done by solving a system of linear equations. The task is now to find all $(\Delta,T)$-inbodies which we define by:
\bdefi[$(\Delta,T)$-inbody]\label{Def:Inbody}
Let $\Delta\subseteq \mathbb{R}^2$ be a triangle and $T\subseteq \mathbb{R}^2$ be a nonempty convex set. A $(\Delta,T)$-inbody is a set $S \subseteq \mathbb{R}^2$ which fulfills:
\begin{itemize}
\item[(i)]
\beqq S = \lambda \Delta + u\eeqq
for some $\lambda > 0$ and $u \in \mathbb{R}^2$.
\item[(ii)] All three vertices of $S$ are contained in $\partial T$. 
\item[(iii)] If
\beqq \{v_1,v_2,v_3\} = V(S),\eeqq
then there is no hyperplane $H$ through the origin, such that
\beqq N_T(v_1) \cup N_T(v_2) \cup N_T(v_3)\eeqq
is contained in one of the two closed halfspaces defined by $H$.
\end{itemize}
\edefi
\noindent Here, $V(S)$ denotes the set of vertices of $S$. The idea behind such a $(\Delta,T)$-inbody is to find the closed polygonal curve $p$. More precisely, we will choose $p$ as the closed polygonal curve having $v_1,v_2,v_3$ as vertices. (i) ensures that the pair $(p,q)$ fulfills the second line of \eqref{eq:System}. Later on, in this algorithm, we need to construct another triangle from outer normal vectors at the points $v_1,v_2,v_3$. Property (iii) ensures that this is possible. Finding all $(\Delta,T)$-inbodies is simple due to the following proposition. 

\bprop[Lemma 4.3.8 in \cite{Krupp2021}] \label{NiceInbody}
Let $T \subset \mathbb{R}^2$ be a convex polytope and $\Delta \subset \mathbb{R}^2$ be a triangle. If there is a $(\Delta,T)$-inbody, then
\beqq S^\ast = \lambda ^\ast \Delta + u^\ast\eeqq
is the only $(\Delta,T)$-inbody, where $\lambda ^\ast, u^\ast$ is a solution of 
\begin{align}\label{eq:bigInbody}
\max\limits &\hspace*{5mm} \lambda \notag\\
\textup{such that} & \hspace*{5mm} \lambda \geq 0, \ u\in \mathbb{R}^2,\\
& \hspace*{5mm} \lambda \Delta + u \subseteq T.\notag
\end{align}
\eprop

\begin{figure}[h]
\centering
\begin{tikzpicture}[scale = 1.5]
\draw (-2,0) -- (5,0);
\draw (0,0) -- (1.5,2.2) -- (3,0) -- (0,0);
\draw (0.4,1.75) -- (1.4,3.216666667) -- (2.4,1.75) -- (0.4,1.75);
\draw[dashed, ->, > = stealth] (1.4,3.216666667) -- (1.0,1.466666667);
\draw[dashed, ->, > = stealth] (0.4,1.75) -- (0,0);
\node at (1.4,3.4) {\small $v_1$};
\node at (1.5,2.3) {\small $w_1$};
\node at (0.2,1.8) {\small $v_2$};
\node at (2.6,1.8) {\small $v_3$};
\node at (-0.1,-0.2) {\small $w_2$};
\node at (3.1,-0.2) {\small $w_3$};
\node at (1.5,0.75) {\small $S^\ast$};
\node at (2.1,2.6) {\small $S$};
\node at (4.5,1) {\small $H_1^+$};
\node at (5.35,0) {\small $H_1$};
\end{tikzpicture}
\caption{Depiction of the $(\Delta,T)$-inbody $S$ with vertices $v_1,v_2,v_3$ and the triangle $S^\ast$ with vertices $w_1,w_2,w_3$. The dashed arrows indicate the location of the line segment $[v_1,v_2]$ after shifting it by $w_2-v_2$.}
\label{img:Picture 9}
\end{figure}
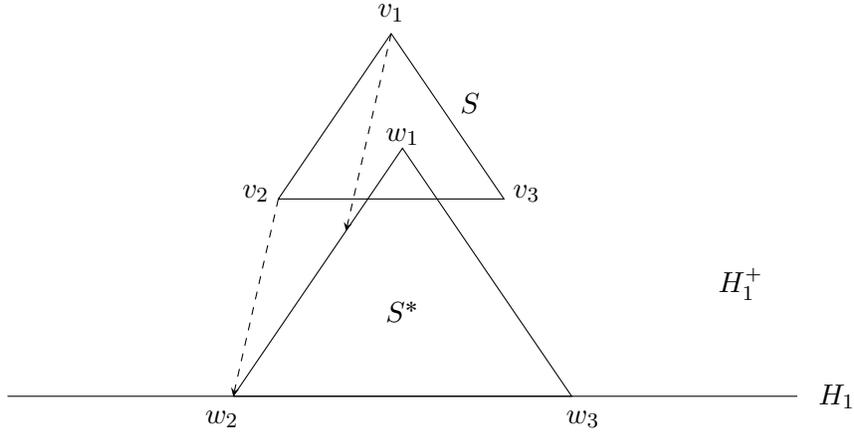

\bpf
Let
\beqq S = \textup{conv} \{v_1,v_2,v_3\}\eeqq
\beqq S^\ast = \textup{conv} \{w_1,w_2,w_3\}\eeqq
as in the claim. Here, we choose the names for the vertices such that there are $\mu > 0$ and $x \in \mathbb{R}^2$ with
\beqq w_j = \mu v_j + x\eeqq
for $j \in \{1,2,3\}$. Note that such a naming is possible since $S$ is a scaled translate of $S^\ast$. We start by letting $H_1,H_2,H_3$ be three lines defined by 
\begin{align*}
w_2,w_3 \in H_1\, , \ w_1,w_3 \in H_2 \, \textup{ and } \, w_1,w_2 \in H_3.
\end{align*}
Each of these lines is the affine hull of a facet of $S^\ast$
. Furthermore, each line $H_j$ devides the plane $\mathbb{R}^2$ in two halfspaces. We denote the halfspace which contains $S^\ast$ by $H_j^+$. If $S = S^\ast$ there is nothing to show. So, we assume $S\neq S^\ast$. Because $S$ is a smaller (or equal size) version of $S^\ast$, it is contained in $H_j^+$ for some $j \in \{1,2,3\}$. Without loss of generality we assume $S\subset H_1^+$, as the other cases can be treated similarly. This situation is depicted in Figure \ref{img:Picture 9}.

Consider the triangle with vertices $v_1,w_2,w_3$. Neither $v_2$ nor $v_3$ can be contained in the interior of this triangle. The reason for this is that by property (ii) of Definition \ref{Def:Inbody} $v_2$ and $v_3$ are boundary points of $T$ but
\begin{align*}
\textup{conv} \{v_1,w_2,w_3\} \subseteq T.
\end{align*}
We will now show that property (iii) is violated for $S$. This contradiction implies $S=S^\ast$ and finishes the proof. More precisely, we claim that $N_T(v_j)$ is contained in the halfspace
\begin{align*}
I = \{x \in \mathbb{R}^2 \colon \langle w_2 - v_2, x \rangle \leq 0\}.
\end{align*}
for every $j\in \{1,2,3\}$. By definition of the outer normal cone it follows immediately that $N_T(v_2) \subset I$. Because $S$ is a smaller (or equal size) version of $S^\ast$ we have $w_2-w_3 = \alpha (v_2 - v_3)$ for some $\alpha \geq 1$. For any $x \in N_T(v_3)$ this implies:
\begin{align*}
\langle w_2 - v_2, x \rangle &= \langle w_2 - v_2 + w_3 -w_3 + v_3-v_3, x \rangle\\
&= \langle (w_2-w_3) - (v_2-v_3), x \rangle + \langle w_3 - v_3, x \rangle\\
&\leq \langle w_2-w_3, x \rangle - \langle v_2-v_3, x \rangle\\
&= (\alpha - 1) \langle v_2-v_3, x \rangle\\
&\leq 0.
\end{align*}
Thus,
\beqq N_T(v_3)\subseteq I.\eeqq
Next, we observe that if we shift $S$ by $w_2-v_2$, then the face $[v_1,v_2]$ of $S$ is contained in $[w_1,w_2]$ (see Figure \ref{img:Picture 9}). So,
\beqq v_1 + (w_2-v_2)\eeqq
is contained in $S^\ast \subseteq T$. Now, for any $y \in N_T(v_1)$ we get:
\begin{align*}
0 \geq \langle v_1 + (w_2-v_2) - v_1, y \rangle = \langle w_2 - v_2,y\rangle 
\end{align*}
As desired this yields
\beqq N_T(v_1)\subseteq I.\eeqq
\epf

\begin{figure}[h!]
\centering
\def\svgwidth{250pt}
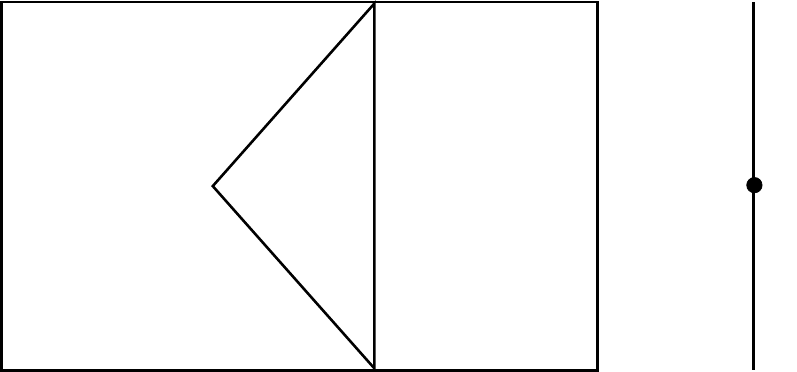
\caption[A situation without $(\Delta,T)$-inbody]{A situation, where there is no $(\Delta,T)$-inbody. The gray area is a scaled translate of $\Delta$. Also shown is a hyperplane $H$ through the origin $O$.}
\label{img:Picture12}
\end{figure}

We point out that there is not always a $(\Delta,T)$-inbody. For example, if \eqref{eq:bigInbody} has multiple optimal solutions, the proof shows that there is no $(\Delta,T)$-inbody. An example for this situation is depicted in Figure \ref{img:Picture12}.

As we can see, there are multiple optimal solutions for \eqref{eq:bigInbody}, since we can shift the gray area to the left and right. The only way to have all vertices of this area on $\partial T$ is to shift it to the left. Then all the corresponding normal vectors are contained in the halfspace on the left of $H$. 

With Proposition \ref{NiceInbody}, we can reduce the search of $(\Delta,T)$-inbodies to a simple maximization problem which we can formulate as a linear problem (LP). It is clear that this problem has an optimal solution as long as $T$ is compact. After we found a solution $\lambda^\ast,u^\ast$, we check whether
\beqq S^\ast=\lambda^\ast \Delta + u^\ast\eeqq
fulfills properties (i)-(iii) in Definition \ref{Def:Inbody}. As pointed out before, it suffices to consider any optimal solution. If $S^\ast$ does not meet properties (i)-(iii), then there is no $(\Delta,T)$-inbody and we proceed with the next choice of $F_1,F_2,F_3$. Otherwise, we take a unit normal vector from $N_T(v)$ for each vertex $v$ of $S$. We construct another triangle as before, using only positive multiples of these normal vectors. It is notable, that $v$ may be a vertex of $T$. In this case, $N_T(v)$ is not one-dimensional and the choice of the corresponding normal vector is not unique. As mentioned earlier, one way of handling this case is to slightly perturb the vertices of $T$. As follows from Proposition \ref{Prop:samefacessamelength}, it is sufficient to find one vector in $N_T(v)$ such that the remaining steps of the algorithm are carried out successfully. So, another way is to sample $N_T(v)$, i.e., only consider finitely many unit normal vectors. The remainder of the algorithm is straight forward and only uses strategies which have been discussed before.

Regarding efficiency, we point out that the algorithm for finding closed polygonal curves with 2 vertices takes
\beqq \mathcal{O}(\vert V(K) \vert^2 \cdot \vert V(T) \vert^2)\eeqq
iterations before it terminates. This is clear since the number of facets of a two dimensional convex polytope equals the number of its vertices. For each of the two convex polytopes $K$ and $T$, the algorithm considers at most one choice for $F_1,F_2,G_1,G_2$ per iteration. In each iteration, we search for the points $q_1,q_2,p_1,p_2$. In the worst case (i.e., if $F_1,F_2,G_1,G_2$ are facets), we solve an LP with $10$ variables and
\beqq 2(\vert V(K) \vert + \vert V(T) \vert + 3)\eeqq
constraints. In order to solve LPs, we use the \texttt{conelp} solver of CVXOPT. This solver relies on a primal-dual path-following method. It is well known that linear problems can be solved in polynomial time (cf.\;\cite{Schrijver1998}).

For finding closed polygonal curves with 3 vertices, the algorithm takes
\beqq \mathcal{O}(\vert V(K) \vert^3)\eeqq
iterations to consider every choice of faces $F_1,F_2,F_3$. In each iteration, we solve the maximization problem stated in Proposition \ref{NiceInbody}. This is an LP with $3$ variables and
\beqq 3\vert V(T) \vert + 1\eeqq
constraints. The remainder of the loop for $F_1,F_2,F_3$ can be realized with running time $\mathcal{O}(1)$. Finally, we note that the calculations for each choice of faces are independent of each other. Therefore, we use parallel computing to speed up the calculations.

In the following Table \ref{table:timesminkbill}, we examine the running time of the algorithm outlined in the Sections \ref{Subsec:GeneralConstruction} and \ref{Subsec:Implementation}. To do so, we let $K$ and $T$ be two-dimensional convex polytopes and consider three different cases. First, we regard the case where both $K$ and $T$ have the same number of vertices. In the second case, $K$ will have a small number of vertices and, in the third case, we chose $T$ to have few vertices. Each time, the convex polytopes $K$ and $T$ have been chosen randomly in the following way. We take a number of normally distributed points and compute their convex hull. Since many of these points will be close to the origin and are unnecessary, we scale each of these normally distributed points to have a random length in $[1,3]$ before we calculate the convex hull. As the number of points grows, the convex hull resembles a ball of radius $3$ due to the normal distribution. So, to accelerate this process, we reduced the interval $[1,3]$ for convex polytopes with many vertices ($\geq 30$).

We compare the running time for finding an $\ell_T$-minimizing closed polygonal curve with 2 vertices to the running time for finding an $\ell_T$-minimizing closed polygonal curve with 3 vertices. As we can see the running time for 2 vertices is approximately symmetric in $\left\lvert V(K) \right\rvert$ and $\left\lvert V(T) \right\rvert$. In contrast to this, the running time for 3 vertices mainly depends on $\left\lvert V(K) \right\rvert$.

All calculations have been done on a Dell Latitude E6530 laptop with Intel Core i7-3520M processor, 2.9 GHz (capable of running four threads). The algorithm and a detailed description on how to choose the input is available on the website \href{www.github.com/S-Krupp/EHZ-capacity-of-polytopes}{www.github.com/S-Krupp/EHZ-capacity-of-polytopes}.

\begin{table}
\centering
\begin{tabular}[h]{c|c|c|c}
$\left\lvert V(K) \right\rvert$ & $\left\lvert V(T) \right\rvert$ & time for 2 vert. in s. & time for 3 vert. in s.\\
\hline
5  & 5  & 1.28237 & 0.22158\\
10 & 10 & 7.85068 & 2.54632\\
15 & 15 & 27.39097 & 11.71706\\
20 & 20 & 54.75007 & 25.03621\\
25 & 25 & 82.30637 & 60.11260\\
30 & 30 & 125.05111 & 110.61238\\
35 & 35 & 170.22497 & 181.65273\\
40 & 40 & 259.88731 & 302.30844\\
45 & 45 & 266.73415 & 385.03827\\
50 & 50 & 361.56254 & 609.04153\\
55 & 55 & 451.56054 & 786.54793\\
\hline
5 & 10 & 3.02675 & 0.22164\\
5 & 15 & 5.57299 & 0.22637\\
5 & 20 & 11.40925 & 0.22114\\
5 & 25 & 16.91015 & 0.23931\\
5 & 30 & 19.89903 & 0.21554\\
5 & 35 & 23.96365 & 0.39383\\ 
5 & 40 & 29.05107 & 0.21106\\
5 & 45 & 32.18348 & 0.54072\\
5 & 50 & 36.41029 & 0.38885\\
5 & 55 & 49.02657 & 0.57020\\
5 & 65 & 59.80655 & 0.79811\\
5 & 75 & 67.33951 & 0.72834\\
\hline
10 & 5 & 3.39280 & 1.19793\\
15 & 5 & 5.74532 & 4.32675\\
20 & 5 & 10.50168 & 11.31948\\
25 & 5 & 14.59203 & 24.99738\\
30 & 5 & 17.76183 & 45.60183\\
35 & 5 & 20.62535 & 90.62127\\
40 & 5 & 23.89690 & 137.52914\\
45 & 5 & 25.73543 & 170.43779\\
50 & 5 & 30.23246 & 266.66650\\
55 & 5 & 33.68478 & 345.84228\\
65 & 5 & 41.49229 & 558.22820\\
75 & 5 & 51.92742 & 937.36931
\end{tabular}
\vspace*{3mm}
\caption[Running times of the algorithm computing the lengths of shortest Minkowski billiard trajectories]{Running times for the calculations of an $\ell_T$-minimizing closed $(K,T)$-Minkowski billiard trajectory with 2 (resp. 3) vertices as described in Section \ref{Sec:Construction}. All numbers are given in seconds.}\label{table:timesminkbill}
\end{table}

\section{A note on Minkowski billiard trajectories on obtuse triangles}\label{Sec:ObtuseTriangle}

It is an open problem for already a long time whether obtuse triangles $\Delta_{>\frac{\pi}{2}}\subset\R^2$ possess closed regular $(\Delta_{>\frac{\pi}{2}},B_1^2)$-Minkowski billiard trajectories (cf.\;\cite{HalbeisenHungerbuhler2000}), i.e., closed classical Euclidean billiard trajectories. The strongest result so far is the existence of a closed classical Euclidean billiard trajectory on triangles with angles not greater than $100^\circ$ (cf.\;\cite{Schwarz2009}).

Obviously, there cannot exist closed regular $(\Delta_{>\frac{\pi}{2}},B_1^2)$-Minkowski billiard trajectories with two bouncing points (cf.\;\cite[Proposition 2.6]{KruppRudolf2020}). Using our algorithm described in Section \ref{Sec:Construction}, we can reason that there cannot be closed regular $(\Delta_{>\frac{\pi}{2}},B_1^2)$-Minkowski billiard trajectories with three bouncing points neither: searching for closed regular $(\Delta_{>\frac{\pi}{2}},B_1^2)$-Minkowski billiard trajectories with three bouncing points means, among other aspects, searching for closed polygonal curves (one for every choice of order of the edges of $\Delta_{>\frac{\pi}{2}}$) with three vertices consisting of line segments given by negative multiples of the outer unit normal vectors at the edges of $\Delta_{>\frac{\pi}{2}}$ that have all three vertices on the sphere $S^1$. But since $\Delta_{>\frac{\pi}{2}}$ is obtuse, for geometrical reasons neither of these two closed polygonal curves with vertices (which are meant to be the closed dual billiard trajectories on $B_1^2$) on $S^1$ are in $F(B_1^2)$. With Proposition \ref{Prop:notranslation}, this implies that there is no closed regular $(\Delta_{>\frac{\pi}{2}},B_1^2)$-Minkowski billiard trajectory.

However, instead of solving the original problem, we can use our algorithm from Section \ref{Sec:Construction} in order to determine the family $\mathcal{T}$ of all convex bodies $T\subset\R^2$ admitting the existence of a closed regular $(\Delta_{>\frac{\pi}{2}},T)$-Minkowski billiard trajectory with three bouncing points: Let $D(\alpha)$ be the rotation matrix in $\R^2$ rotating counter clockwise by angle $\alpha$. Then, $\mathcal{T}$ is the set of all convex bodies $T\subset\R^2$ for which there are either $\lambda^+>0$ and $\xi^+\in\R^2$ with
\beqq (\lambda^+D(\pi/2)\Delta_{>\frac{\pi}{2}}+\xi^+)\cap \partial T = \{p_1,p_2,p_3\}\in F(T)\eeqq
and maximally spanning unit normal vectors $n_T(p_1),n_T(p_2),n_T(p_3)$ in the outer normal cones at $p_1,p_2,p_3$, or there are $\lambda^- >0$ and $\xi^-\in\R^2$ such that the same holds for
\beqq \lambda^-D(-\pi/2)\Delta_{>\frac{\pi}{2}}+\xi^-.\eeqq

\section*{Acknowledgement}
This research is supported by the SFB/TRR 191 'Symplectic Structures in Geometry, Algebra and Dynamics', funded by the \underline{German Research Foundation}, and was carried out under the supervision of Alberto Abbondandolo (Ruhr-Universit\"at Bochum) and Frank Vallentin (Universit\"at zu K\"oln). The authors are thankful to the supervisors' support.


\medskip

\section*{Stefan Krupp, Universit\"at zu K\"oln, Mathematisches Institut, Weyertal 86-90, D-50931 K\"oln, Germany.}
\center{E-mail address: krupp@math.uni-koeln.de}

\section*{Daniel Rudolf, Ruhr-Universit\"at Bochum, Fakult\"at f\"ur Mathematik, Universit\"atsstrasse 150, D-44801 Bochum, Germany.}
\center{E-mail address: daniel.rudolf@ruhr-uni-bochum.de}

\end{document}

%% file: images/Stossregelmink.pdf_tex
\begingroup%
  \makeatletter%
  \providecommand\color[2][]{%
    \errmessage{(Inkscape) Color is used for the text in Inkscape, but the package 'color.sty' is not loaded}%
    \renewcommand\color[2][]{}%
  }%
  \providecommand\transparent[1]{%
    \errmessage{(Inkscape) Transparency is used (non-zero) for the text in Inkscape, but the package 'transparent.sty' is not loaded}%
    \renewcommand\transparent[1]{}%
  }%
  \providecommand\rotatebox[2]{#2}%
  \newcommand*\fsize{\dimexpr\f@size pt\relax}%
  \newcommand*\lineheight[1]{\fontsize{\fsize}{#1\fsize}\selectfont}%
  \ifx\svgwidth\undefined%
    \setlength{\unitlength}{285.81452862bp}%
    \ifx\svgscale\undefined%
      \relax%
    \else%
      \setlength{\unitlength}{\unitlength * \real{\svgscale}}%
    \fi%
  \else%
    \setlength{\unitlength}{\svgwidth}%
  \fi%
  \global\let\svgwidth\undefined%
  \global\let\svgscale\undefined%
  \makeatother%
  \begin{picture}(1,0.69875334)%
    \lineheight{1}%
    \setlength\tabcolsep{0pt}%
    \put(0,0){\includegraphics[width=\unitlength,page=1]{Stossregelmink.pdf}}%
    \put(0.86349328,0.42473384){\color[rgb]{0,0,0}\makebox(0,0)[lt]{\lineheight{1.25}\smash{\begin{tabular}[t]{l}$H_j$\end{tabular}}}}%
    \put(0.71786167,0.0082515){\color[rgb]{0,0,0}\makebox(0,0)[lt]{\lineheight{1.25}\smash{\begin{tabular}[t]{l}$q_{j+1}$\end{tabular}}}}%
    \put(-0.00290425,0.09233721){\color[rgb]{0,0,0}\makebox(0,0)[lt]{\lineheight{1.25}\smash{\begin{tabular}[t]{l}$q_{j-1}$\end{tabular}}}}%
    \put(0.56097889,0.60368355){\color[rgb]{0,0,0}\makebox(0,0)[lt]{\lineheight{1.25}\smash{\begin{tabular}[t]{l}$q_j$\end{tabular}}}}%
    \put(0.0728596,0.42675179){\color[rgb]{0,0,0}\makebox(0,0)[lt]{\lineheight{1.25}\smash{\begin{tabular}[t]{l}$K$\end{tabular}}}}%
    \put(0,0){\includegraphics[width=\unitlength,page=2]{Stossregelmink.pdf}}%
    \put(0.38372694,0.64717047){\color[rgb]{0,0,0}\makebox(0,0)[lt]{\lineheight{1.25}\smash{\begin{tabular}[t]{l}$\widebar{q}_j$\end{tabular}}}}%
    \put(0,0){\includegraphics[width=\unitlength,page=3]{Stossregelmink.pdf}}%
  \end{picture}%
\endgroup%

%% file: images/ReflectionRule1.pdf_tex
\begingroup%
  \makeatletter%
  \providecommand\color[2][]{%
    \errmessage{(Inkscape) Color is used for the text in Inkscape, but the package 'color.sty' is not loaded}%
    \renewcommand\color[2][]{}%
  }%
  \providecommand\transparent[1]{%
    \errmessage{(Inkscape) Transparency is used (non-zero) for the text in Inkscape, but the package 'transparent.sty' is not loaded}%
    \renewcommand\transparent[1]{}%
  }%
  \providecommand\rotatebox[2]{#2}%
  \newcommand*\fsize{\dimexpr\f@size pt\relax}%
  \newcommand*\lineheight[1]{\fontsize{\fsize}{#1\fsize}\selectfont}%
  \ifx\svgwidth\undefined%
    \setlength{\unitlength}{400.87343144bp}%
    \ifx\svgscale\undefined%
      \relax%
    \else%
      \setlength{\unitlength}{\unitlength * \real{\svgscale}}%
    \fi%
  \else%
    \setlength{\unitlength}{\svgwidth}%
  \fi%
  \global\let\svgwidth\undefined%
  \global\let\svgscale\undefined%
  \makeatother%
  \begin{picture}(1,0.66776735)%
    \lineheight{1}%
    \setlength\tabcolsep{0pt}%
    \put(0,0){\includegraphics[width=\unitlength,page=1]{ReflectionRule1.pdf}}%
    \put(0.25458959,0.48286674){\color[rgb]{0,0,0}\makebox(0,0)[lt]{\lineheight{1.25}\smash{\begin{tabular}[t]{l}$T^\circ$\end{tabular}}}}%
    \put(0,0){\includegraphics[width=\unitlength,page=2]{ReflectionRule1.pdf}}%
    \put(0.0148858,0.09328872){\color[rgb]{0,0,0}\makebox(0,0)[lt]{\lineheight{1.25}\smash{\begin{tabular}[t]{l}$q_{j-1}$\end{tabular}}}}%
    \put(0.4468371,0.05656639){\color[rgb]{0,0,0}\makebox(0,0)[lt]{\lineheight{1.25}\smash{\begin{tabular}[t]{l}$q_{j+1}$\end{tabular}}}}%
    \put(0.31748623,0.34010566){\color[rgb]{0,0,0}\makebox(0,0)[lt]{\lineheight{1.25}\smash{\begin{tabular}[t]{l}$q_j$\end{tabular}}}}%
    \put(0.01173623,0.37419702){\color[rgb]{0,0,0}\makebox(0,0)[lt]{\lineheight{1.25}\smash{\begin{tabular}[t]{l}$H_j$\end{tabular}}}}%
    \put(0.57114943,0.5924815){\color[rgb]{0,0,0}\makebox(0,0)[lt]{\lineheight{1.25}\smash{\begin{tabular}[t]{l}$\nabla \mu_{T^\circ}(q_j-q_{j-1})$\end{tabular}}}}%
    \put(0.18425321,0.08180308){\color[rgb]{0,0,0}\makebox(0,0)[lt]{\lineheight{1.25}\smash{\begin{tabular}[t]{l}$\nabla\mu_{T^\circ}(q_{j+1}-q_j)$\end{tabular}}}}%
    \put(0.07126836,0.23623355){\color[rgb]{0,0,0}\makebox(0,0)[lt]{\lineheight{1.25}\smash{\begin{tabular}[t]{l}$K$\end{tabular}}}}%
    \put(0,0){\includegraphics[width=\unitlength,page=3]{ReflectionRule1.pdf}}%
    \put(0.24571383,0.40803105){\color[rgb]{0,0,0}\makebox(0,0)[lt]{\lineheight{1.25}\smash{\begin{tabular}[t]{l}$n_{H_j}$\end{tabular}}}}%
    \put(0,0){\includegraphics[width=\unitlength,page=4]{ReflectionRule1.pdf}}%
    \put(0.93073398,0.14994206){\color[rgb]{0,0,0}\makebox(0,0)[lt]{\lineheight{1.25}\smash{\begin{tabular}[t]{l}$n_{H_j}$\end{tabular}}}}%
    \put(0.76370728,0.00294157){\color[rgb]{0,0,0}\makebox(0,0)[lt]{\lineheight{1.25}\smash{\begin{tabular}[t]{l}$\nabla\mu_{T^\circ}(q_{j+1}-q_j)$\end{tabular}}}}%
    \put(0.77380325,0.26420151){\color[rgb]{0,0,0}\makebox(0,0)[lt]{\lineheight{1.25}\smash{\begin{tabular}[t]{l}$\nabla \mu_{T^\circ}(q_j-q_{j-1})$\end{tabular}}}}%
  \end{picture}%
\endgroup%

%% file: images/ReflectionRule2.pdf_tex
\begingroup%
  \makeatletter%
  \providecommand\color[2][]{%
    \errmessage{(Inkscape) Color is used for the text in Inkscape, but the package 'color.sty' is not loaded}%
    \renewcommand\color[2][]{}%
  }%
  \providecommand\transparent[1]{%
    \errmessage{(Inkscape) Transparency is used (non-zero) for the text in Inkscape, but the package 'transparent.sty' is not loaded}%
    \renewcommand\transparent[1]{}%
  }%
  \providecommand\rotatebox[2]{#2}%
  \newcommand*\fsize{\dimexpr\f@size pt\relax}%
  \newcommand*\lineheight[1]{\fontsize{\fsize}{#1\fsize}\selectfont}%
  \ifx\svgwidth\undefined%
    \setlength{\unitlength}{453.06517667bp}%
    \ifx\svgscale\undefined%
      \relax%
    \else%
      \setlength{\unitlength}{\unitlength * \real{\svgscale}}%
    \fi%
  \else%
    \setlength{\unitlength}{\svgwidth}%
  \fi%
  \global\let\svgwidth\undefined%
  \global\let\svgscale\undefined%
  \makeatother%
  \begin{picture}(1,0.48893334)%
    \lineheight{1}%
    \setlength\tabcolsep{0pt}%
    \put(0,0){\includegraphics[width=\unitlength,page=1]{ReflectionRule2.pdf}}%
    \put(-0.00091607,0.10212351){\color[rgb]{0,0,0}\makebox(0,0)[lt]{\lineheight{1.25}\smash{\begin{tabular}[t]{l}$q_{j-1}$\end{tabular}}}}%
    \put(0.48050791,0.15010533){\color[rgb]{0,0,0}\makebox(0,0)[lt]{\lineheight{1.25}\smash{\begin{tabular}[t]{l}$q_{j+1}$\end{tabular}}}}%
    \put(0.15356825,0.31747748){\color[rgb]{0,0,0}\makebox(0,0)[lt]{\lineheight{1.25}\smash{\begin{tabular}[t]{l}$q_j$\end{tabular}}}}%
    \put(0.01959378,0.35523095){\color[rgb]{0,0,0}\makebox(0,0)[lt]{\lineheight{1.25}\smash{\begin{tabular}[t]{l}$N_K(q_j)$\end{tabular}}}}%
    \put(0.33632386,0.29106843){\color[rgb]{0,0,0}\makebox(0,0)[lt]{\lineheight{1.25}\smash{\begin{tabular}[t]{l}$K$\end{tabular}}}}%
    \put(0,0){\includegraphics[width=\unitlength,page=2]{ReflectionRule2.pdf}}%
    \put(0.80604683,0.37370725){\color[rgb]{0,0,0}\makebox(0,0)[lt]{\lineheight{1.25}\smash{\begin{tabular}[t]{l}$p_{j-1}$\end{tabular}}}}%
    \put(0.928456,0.10134673){\color[rgb]{0,0,0}\makebox(0,0)[lt]{\lineheight{1.25}\smash{\begin{tabular}[t]{l}$p_j$\end{tabular}}}}%
    \put(0.58800534,0.12735873){\color[rgb]{0,0,0}\makebox(0,0)[lt]{\lineheight{1.25}\smash{\begin{tabular}[t]{l}$p_{j+1}$\end{tabular}}}}%
    \put(0.62396306,0.33009896){\color[rgb]{0,0,0}\makebox(0,0)[lt]{\lineheight{1.25}\smash{\begin{tabular}[t]{l}$T$\end{tabular}}}}%
    \put(0,0){\includegraphics[width=\unitlength,page=3]{ReflectionRule2.pdf}}%
  \end{picture}%
\endgroup%

%% file: images/upperlowerboundary.pdf_tex
\begingroup%
  \makeatletter%
  \providecommand\color[2][]{%
    \errmessage{(Inkscape) Color is used for the text in Inkscape, but the package 'color.sty' is not loaded}%
    \renewcommand\color[2][]{}%
  }%
  \providecommand\transparent[1]{%
    \errmessage{(Inkscape) Transparency is used (non-zero) for the text in Inkscape, but the package 'transparent.sty' is not loaded}%
    \renewcommand\transparent[1]{}%
  }%
  \providecommand\rotatebox[2]{#2}%
  \newcommand*\fsize{\dimexpr\f@size pt\relax}%
  \newcommand*\lineheight[1]{\fontsize{\fsize}{#1\fsize}\selectfont}%
  \ifx\svgwidth\undefined%
    \setlength{\unitlength}{355.31280466bp}%
    \ifx\svgscale\undefined%
      \relax%
    \else%
      \setlength{\unitlength}{\unitlength * \real{\svgscale}}%
    \fi%
  \else%
    \setlength{\unitlength}{\svgwidth}%
  \fi%
  \global\let\svgwidth\undefined%
  \global\let\svgscale\undefined%
  \makeatother%
  \begin{picture}(1,0.60610578)%
    \lineheight{1}%
    \setlength\tabcolsep{0pt}%
    \put(0,0){\includegraphics[width=\unitlength,page=1]{upperlowerboundary.pdf}}%
    \put(0.04741876,0.50493735){\color[rgb]{0,0,0}\makebox(0,0)[lt]{\lineheight{1.25}\smash{\begin{tabular}[t]{l}$n_{H_j}$\end{tabular}}}}%
    \put(0,0){\includegraphics[width=\unitlength,page=2]{upperlowerboundary.pdf}}%
    \put(0.26724518,0.08367878){\color[rgb]{0,0,0}\makebox(0,0)[lt]{\lineheight{1.25}\smash{\begin{tabular}[t]{l}$\partial T_{H_j^-}$\end{tabular}}}}%
    \put(0.87222903,0.52108515){\color[rgb]{0,0,0}\makebox(0,0)[lt]{\lineheight{1.25}\smash{\begin{tabular}[t]{l}$\partial T_{H_j^+}$\end{tabular}}}}%
    \put(0,0){\includegraphics[width=\unitlength,page=3]{upperlowerboundary.pdf}}%
    \put(-0.00233619,0.35506942){\color[rgb]{0,0,0}\makebox(0,0)[lt]{\lineheight{1.25}\smash{\begin{tabular}[t]{l}$H_j$\end{tabular}}}}%
    \put(0.17956447,0.3241762){\color[rgb]{0,0,0}\makebox(0,0)[lt]{\lineheight{1.25}\smash{\begin{tabular}[t]{l}$H_j^-$\end{tabular}}}}%
    \put(0.17654902,0.383279){\color[rgb]{0,0,0}\makebox(0,0)[lt]{\lineheight{1.25}\smash{\begin{tabular}[t]{l}$H_j^+$\end{tabular}}}}%
  \end{picture}%
\endgroup%

%% file: images/bouncingruleinqj.pdf_tex
\begingroup%
  \makeatletter%
  \providecommand\color[2][]{%
    \errmessage{(Inkscape) Color is used for the text in Inkscape, but the package 'color.sty' is not loaded}%
    \renewcommand\color[2][]{}%
  }%
  \providecommand\transparent[1]{%
    \errmessage{(Inkscape) Transparency is used (non-zero) for the text in Inkscape, but the package 'transparent.sty' is not loaded}%
    \renewcommand\transparent[1]{}%
  }%
  \providecommand\rotatebox[2]{#2}%
  \newcommand*\fsize{\dimexpr\f@size pt\relax}%
  \newcommand*\lineheight[1]{\fontsize{\fsize}{#1\fsize}\selectfont}%
  \ifx\svgwidth\undefined%
    \setlength{\unitlength}{337.84782477bp}%
    \ifx\svgscale\undefined%
      \relax%
    \else%
      \setlength{\unitlength}{\unitlength * \real{\svgscale}}%
    \fi%
  \else%
    \setlength{\unitlength}{\svgwidth}%
  \fi%
  \global\let\svgwidth\undefined%
  \global\let\svgscale\undefined%
  \makeatother%
  \begin{picture}(1,0.56822496)%
    \lineheight{1}%
    \setlength\tabcolsep{0pt}%
    \put(0,0){\includegraphics[width=\unitlength,page=1]{bouncingruleinqj.pdf}}%
    \put(-0.00122848,0.14207749){\color[rgb]{0,0,0}\makebox(0,0)[lt]{\lineheight{1.25}\smash{\begin{tabular}[t]{l}$q_{j-1}$\end{tabular}}}}%
    \put(0.59667114,0.13675091){\color[rgb]{0,0,0}\makebox(0,0)[lt]{\lineheight{1.25}\smash{\begin{tabular}[t]{l}$q_{j+1}$\end{tabular}}}}%
    \put(0.33326036,0.45458154){\color[rgb]{0,0,0}\makebox(0,0)[lt]{\lineheight{1.25}\smash{\begin{tabular}[t]{l}$q_{j}$\end{tabular}}}}%
    \put(0.26677902,0.52436946){\color[rgb]{0,0,0}\makebox(0,0)[lt]{\lineheight{1.25}\smash{\begin{tabular}[t]{l}$n_{H_j}$\end{tabular}}}}%
    \put(0.13143784,0.46355145){\color[rgb]{0,0,0}\makebox(0,0)[lt]{\lineheight{1.25}\smash{\begin{tabular}[t]{l}$H_j$\end{tabular}}}}%
    \put(0.29073304,0.2229803){\color[rgb]{0,0,0}\makebox(0,0)[lt]{\lineheight{1.25}\smash{\begin{tabular}[t]{l}$K$\end{tabular}}}}%
    \put(0.76647774,0.30307331){\color[rgb]{0,0,0}\makebox(0,0)[lt]{\lineheight{1.25}\smash{\begin{tabular}[t]{l}$T$\end{tabular}}}}%
    \put(0.87477847,0.08060414){\color[rgb]{0,0,0}\makebox(0,0)[lt]{\lineheight{1.25}\smash{\begin{tabular}[t]{l}$p_{j}$\end{tabular}}}}%
    \put(0.92552005,0.45926154){\color[rgb]{0,0,0}\makebox(0,0)[lt]{\lineheight{1.25}\smash{\begin{tabular}[t]{l}$p_{j-1}$\end{tabular}}}}%
  \end{picture}%
\endgroup%

%% file: images/hyperplanechoice.pdf_tex
\begingroup%
  \makeatletter%
  \providecommand\color[2][]{%
    \errmessage{(Inkscape) Color is used for the text in Inkscape, but the package 'color.sty' is not loaded}%
    \renewcommand\color[2][]{}%
  }%
  \providecommand\transparent[1]{%
    \errmessage{(Inkscape) Transparency is used (non-zero) for the text in Inkscape, but the package 'transparent.sty' is not loaded}%
    \renewcommand\transparent[1]{}%
  }%
  \providecommand\rotatebox[2]{#2}%
  \newcommand*\fsize{\dimexpr\f@size pt\relax}%
  \newcommand*\lineheight[1]{\fontsize{\fsize}{#1\fsize}\selectfont}%
  \ifx\svgwidth\undefined%
    \setlength{\unitlength}{340.67058538bp}%
    \ifx\svgscale\undefined%
      \relax%
    \else%
      \setlength{\unitlength}{\unitlength * \real{\svgscale}}%
    \fi%
  \else%
    \setlength{\unitlength}{\svgwidth}%
  \fi%
  \global\let\svgwidth\undefined%
  \global\let\svgscale\undefined%
  \makeatother%
  \begin{picture}(1,0.78120034)%
    \lineheight{1}%
    \setlength\tabcolsep{0pt}%
    \put(0,0){\includegraphics[width=\unitlength,page=1]{hyperplanechoice.pdf}}%
    \put(0.37579882,0.71772376){\color[rgb]{0,0,0}\makebox(0,0)[lt]{\lineheight{1.25}\smash{\begin{tabular}[t]{l}$\pi_U(q_1=q_{i_1})$\end{tabular}}}}%
    \put(0.05373587,0.46452061){\color[rgb]{0,0,0}\makebox(0,0)[lt]{\lineheight{1.25}\smash{\begin{tabular}[t]{l}$\pi_U(q_2=q_{i_2})$\end{tabular}}}}%
    \put(0.52636465,0.25475872){\color[rgb]{0,0,0}\makebox(0,0)[lt]{\lineheight{1.25}\smash{\begin{tabular}[t]{l}$\pi_U(q_3=q_{i_3})$\end{tabular}}}}%
    \put(0.53327821,0.50619429){\color[rgb]{0,0,0}\makebox(0,0)[lt]{\lineheight{1.25}\smash{\begin{tabular}[t]{l}$\pi_U(q_4)$\end{tabular}}}}%
    \put(-0.0012183,0.62381034){\color[rgb]{0,0,0}\makebox(0,0)[lt]{\lineheight{1.25}\smash{\begin{tabular}[t]{l}$H_1\cap U=\widetilde{H}_1$\end{tabular}}}}%
    \put(0.13022475,0.75530098){\color[rgb]{0,0,0}\makebox(0,0)[lt]{\lineheight{1.25}\smash{\begin{tabular}[t]{l}$H_2\cap U=\widetilde{H}_2$\end{tabular}}}}%
    \put(0.09309093,0.29120379){\color[rgb]{0,0,0}\makebox(0,0)[lt]{\lineheight{1.25}\smash{\begin{tabular}[t]{l}$H_3\cap U$\end{tabular}}}}%
    \put(0.75488156,0.77004931){\color[rgb]{0,0,0}\makebox(0,0)[lt]{\lineheight{1.25}\smash{\begin{tabular}[t]{l}$H_4\cap U$\end{tabular}}}}%
    \put(0.38718079,0.11937277){\color[rgb]{0,0,0}\makebox(0,0)[lt]{\lineheight{1.25}\smash{\begin{tabular}[t]{l}$\widetilde{H}_{3}$\end{tabular}}}}%
    \put(0.3680655,0.543614){\color[rgb]{0,0,0}\makebox(0,0)[lt]{\lineheight{1.25}\smash{\begin{tabular}[t]{l}$\pi_U(K)$\end{tabular}}}}%
  \end{picture}%
\endgroup%

%% file: images/Counterexample.pdf_tex
\begingroup%
  \makeatletter%
  \providecommand\color[2][]{%
    \errmessage{(Inkscape) Color is used for the text in Inkscape, but the package 'color.sty' is not loaded}%
    \renewcommand\color[2][]{}%
  }%
  \providecommand\transparent[1]{%
    \errmessage{(Inkscape) Transparency is used (non-zero) for the text in Inkscape, but the package 'transparent.sty' is not loaded}%
    \renewcommand\transparent[1]{}%
  }%
  \providecommand\rotatebox[2]{#2}%
  \newcommand*\fsize{\dimexpr\f@size pt\relax}%
  \newcommand*\lineheight[1]{\fontsize{\fsize}{#1\fsize}\selectfont}%
  \ifx\svgwidth\undefined%
    \setlength{\unitlength}{340.66867711bp}%
    \ifx\svgscale\undefined%
      \relax%
    \else%
      \setlength{\unitlength}{\unitlength * \real{\svgscale}}%
    \fi%
  \else%
    \setlength{\unitlength}{\svgwidth}%
  \fi%
  \global\let\svgwidth\undefined%
  \global\let\svgscale\undefined%
  \makeatother%
  \begin{picture}(1,0.51545115)%
    \lineheight{1}%
    \setlength\tabcolsep{0pt}%
    \put(0,0){\includegraphics[width=\unitlength,page=1]{Counterexample.pdf}}%
    \put(0.24114849,0.14752172){\color[rgb]{0,0,0}\makebox(0,0)[lt]{\lineheight{1.25}\smash{\begin{tabular}[t]{l}$q_1$\end{tabular}}}}%
    \put(0.32763802,0.24344651){\color[rgb]{0,0,0}\makebox(0,0)[lt]{\lineheight{1.25}\smash{\begin{tabular}[t]{l}$q_2$\end{tabular}}}}%
    \put(0.10865057,0.24344654){\color[rgb]{0,0,0}\makebox(0,0)[lt]{\lineheight{1.25}\smash{\begin{tabular}[t]{l}$q_3$\end{tabular}}}}%
    \put(0.72894158,0.41477409){\color[rgb]{0,0,0}\makebox(0,0)[lt]{\lineheight{1.25}\smash{\begin{tabular}[t]{l}$\widetilde{p}_1$\end{tabular}}}}%
    \put(0.5266613,0.23262954){\color[rgb]{0,0,0}\makebox(0,0)[lt]{\lineheight{1.25}\smash{\begin{tabular}[t]{l}$\widetilde{p}_2$\end{tabular}}}}%
    \put(0.72160702,0.05094562){\color[rgb]{0,0,0}\makebox(0,0)[lt]{\lineheight{1.25}\smash{\begin{tabular}[t]{l}$\widetilde{p}_3$\end{tabular}}}}%
    \put(0,0){\includegraphics[width=\unitlength,page=2]{Counterexample.pdf}}%
    \put(0.27038426,0.34157499){\color[rgb]{0,0,0}\makebox(0,0)[lt]{\lineheight{1.25}\smash{\begin{tabular}[t]{l}$H_3$\end{tabular}}}}%
    \put(0.16105893,0.34208628){\color[rgb]{0,0,0}\makebox(0,0)[lt]{\lineheight{1.25}\smash{\begin{tabular}[t]{l}$H_2$\end{tabular}}}}%
    \put(0.42957286,0.17550738){\color[rgb]{0,0,0}\makebox(0,0)[lt]{\lineheight{1.25}\smash{\begin{tabular}[t]{l}$H_1$\end{tabular}}}}%
    \put(0,0){\includegraphics[width=\unitlength,page=3]{Counterexample.pdf}}%
    \put(0.89497916,0.38706391){\color[rgb]{0,0,0}\makebox(0,0)[lt]{\lineheight{1.25}\smash{\begin{tabular}[t]{l}$p_1$\end{tabular}}}}%
    \put(0.55059346,0.05761722){\color[rgb]{0,0,0}\makebox(0,0)[lt]{\lineheight{1.25}\smash{\begin{tabular}[t]{l}$p_2$\end{tabular}}}}%
    \put(0,0){\includegraphics[width=\unitlength,page=4]{Counterexample.pdf}}%
    \put(0.11893188,0.182634){\color[rgb]{0,0,0}\makebox(0,0)[lt]{\lineheight{1.25}\smash{\begin{tabular}[t]{l}$K$\end{tabular}}}}%
    \put(0.83050518,0.10636592){\color[rgb]{0,0,0}\makebox(0,0)[lt]{\lineheight{1.25}\smash{\begin{tabular}[t]{l}$T$\end{tabular}}}}%
  \end{picture}%
\endgroup%

%% file: images/cs1.pdf_tex
\begingroup%
  \makeatletter%
  \providecommand\color[2][]{%
    \errmessage{(Inkscape) Color is used for the text in Inkscape, but the package 'color.sty' is not loaded}%
    \renewcommand\color[2][]{}%
  }%
  \providecommand\transparent[1]{%
    \errmessage{(Inkscape) Transparency is used (non-zero) for the text in Inkscape, but the package 'transparent.sty' is not loaded}%
    \renewcommand\transparent[1]{}%
  }%
  \providecommand\rotatebox[2]{#2}%
  \newcommand*\fsize{\dimexpr\f@size pt\relax}%
  \newcommand*\lineheight[1]{\fontsize{\fsize}{#1\fsize}\selectfont}%
  \ifx\svgwidth\undefined%
    \setlength{\unitlength}{240.37560696bp}%
    \ifx\svgscale\undefined%
      \relax%
    \else%
      \setlength{\unitlength}{\unitlength * \real{\svgscale}}%
    \fi%
  \else%
    \setlength{\unitlength}{\svgwidth}%
  \fi%
  \global\let\svgwidth\undefined%
  \global\let\svgscale\undefined%
  \makeatother%
  \begin{picture}(1,0.5687003)%
    \lineheight{1}%
    \setlength\tabcolsep{0pt}%
    \put(0,0){\includegraphics[width=\unitlength,page=1]{cs1.pdf}}%
    \put(0.14251541,0.11596467){\color[rgb]{0,0,0}\makebox(0,0)[lt]{\lineheight{1.25}\smash{\begin{tabular}[t]{l}$q_1$\end{tabular}}}}%
    \put(0.22111901,0.39004237){\color[rgb]{0,0,0}\makebox(0,0)[lt]{\lineheight{1.25}\smash{\begin{tabular}[t]{l}$q_2$\end{tabular}}}}%
    \put(0.14009935,0.55289658){\color[rgb]{0,0,0}\makebox(0,0)[lt]{\lineheight{1.25}\smash{\begin{tabular}[t]{l}$n_K(q_2)$\end{tabular}}}}%
    \put(0.0188244,0.45814818){\color[rgb]{0,0,0}\makebox(0,0)[lt]{\lineheight{1.25}\smash{\begin{tabular}[t]{l}$\widetilde{n}_K(q_2)$\end{tabular}}}}%
    \put(0.28799958,0.0525861){\color[rgb]{0,0,0}\makebox(0,0)[lt]{\lineheight{1.25}\smash{\begin{tabular}[t]{l}$\widetilde{n}_K(q_1)$\end{tabular}}}}%
    \put(0.07594533,0.02220561){\color[rgb]{0,0,0}\makebox(0,0)[lt]{\lineheight{1.25}\smash{\begin{tabular}[t]{l}$n_K(q_1)$\end{tabular}}}}%
    \put(0,0){\includegraphics[width=\unitlength,page=2]{cs1.pdf}}%
    \put(0.71587091,0.1060833){\color[rgb]{0,0,0}\makebox(0,0)[lt]{\lineheight{1.25}\smash{\begin{tabular}[t]{l}$p_1$\end{tabular}}}}%
    \put(0.71769388,0.43457233){\color[rgb]{0,0,0}\makebox(0,0)[lt]{\lineheight{1.25}\smash{\begin{tabular}[t]{l}$p_2$\end{tabular}}}}%
    \put(0.85409628,0.10218377){\color[rgb]{0,0,0}\makebox(0,0)[lt]{\lineheight{1.25}\smash{\begin{tabular}[t]{l}$\widetilde{p}_1$\end{tabular}}}}%
    \put(0.51301304,0.43068859){\color[rgb]{0,0,0}\makebox(0,0)[lt]{\lineheight{1.25}\smash{\begin{tabular}[t]{l}$\widetilde{p}_2$\end{tabular}}}}%
    \put(0.85468909,0.03854986){\color[rgb]{0,0,0}\makebox(0,0)[lt]{\lineheight{1.25}\smash{\begin{tabular}[t]{l}$n_T(\widetilde{p}_1)$\end{tabular}}}}%
    \put(0.60285139,0.03230962){\color[rgb]{0,0,0}\makebox(0,0)[lt]{\lineheight{1.25}\smash{\begin{tabular}[t]{l}$n_T(p_1)$\end{tabular}}}}%
    \put(0.72030158,0.50255585){\color[rgb]{0,0,0}\makebox(0,0)[lt]{\lineheight{1.25}\smash{\begin{tabular}[t]{l}$n_T(p_2)$\end{tabular}}}}%
    \put(0.45865712,0.5014415){\color[rgb]{0,0,0}\makebox(0,0)[lt]{\lineheight{1.25}\smash{\begin{tabular}[t]{l}$n_T(\widetilde{p}_2)$\end{tabular}}}}%
    \put(0.7749032,0.35212166){\color[rgb]{0,0,0}\makebox(0,0)[lt]{\lineheight{1.25}\smash{\begin{tabular}[t]{l}$T$\end{tabular}}}}%
    \put(0.12257596,0.31423455){\color[rgb]{0,0,0}\makebox(0,0)[lt]{\lineheight{1.25}\smash{\begin{tabular}[t]{l}$K$\end{tabular}}}}%
    \put(0,0){\includegraphics[width=\unitlength,page=3]{cs1.pdf}}%
    \put(0.36215649,0.42678159){\color[rgb]{0,0,0}\makebox(0,0)[lt]{\lineheight{1.25}\smash{\begin{tabular}[t]{l}$H_2$\end{tabular}}}}%
    \put(0,0){\includegraphics[width=\unitlength,page=4]{cs1.pdf}}%
    \put(0.32872664,0.52818536){\color[rgb]{0,0,0}\makebox(0,0)[lt]{\lineheight{1.25}\smash{\begin{tabular}[t]{l}$\widetilde{H}_2$\end{tabular}}}}%
    \put(0,0){\includegraphics[width=\unitlength,page=5]{cs1.pdf}}%
  \end{picture}%
\endgroup%

%% file: images/exampleB2.pdf_tex
\begingroup%
  \makeatletter%
  \providecommand\color[2][]{%
    \errmessage{(Inkscape) Color is used for the text in Inkscape, but the package 'color.sty' is not loaded}%
    \renewcommand\color[2][]{}%
  }%
  \providecommand\transparent[1]{%
    \errmessage{(Inkscape) Transparency is used (non-zero) for the text in Inkscape, but the package 'transparent.sty' is not loaded}%
    \renewcommand\transparent[1]{}%
  }%
  \providecommand\rotatebox[2]{#2}%
  \newcommand*\fsize{\dimexpr\f@size pt\relax}%
  \newcommand*\lineheight[1]{\fontsize{\fsize}{#1\fsize}\selectfont}%
  \ifx\svgwidth\undefined%
    \setlength{\unitlength}{265.50305183bp}%
    \ifx\svgscale\undefined%
      \relax%
    \else%
      \setlength{\unitlength}{\unitlength * \real{\svgscale}}%
    \fi%
  \else%
    \setlength{\unitlength}{\svgwidth}%
  \fi%
  \global\let\svgwidth\undefined%
  \global\let\svgscale\undefined%
  \makeatother%
  \begin{picture}(1,0.65176269)%
    \lineheight{1}%
    \setlength\tabcolsep{0pt}%
    \put(0,0){\includegraphics[width=\unitlength,page=1]{exampleB2.pdf}}%
    \put(0.06142437,0.11532112){\color[rgb]{0,0,0}\makebox(0,0)[lt]{\lineheight{1.25}\smash{\begin{tabular}[t]{l}$q_1$\end{tabular}}}}%
    \put(-0.00003512,0.02490859){\color[rgb]{0,0,0}\makebox(0,0)[lt]{\lineheight{1.25}\smash{\begin{tabular}[t]{l}$n_K(q_1)$\end{tabular}}}}%
    \put(0.05983322,0.53630251){\color[rgb]{0,0,0}\makebox(0,0)[lt]{\lineheight{1.25}\smash{\begin{tabular}[t]{l}$q_2$\end{tabular}}}}%
    \put(-0.00156322,0.60994974){\color[rgb]{0,0,0}\makebox(0,0)[lt]{\lineheight{1.25}\smash{\begin{tabular}[t]{l}$n_K(q_2)$\end{tabular}}}}%
    \put(0.27310759,0.32242281){\color[rgb]{0,0,0}\makebox(0,0)[lt]{\lineheight{1.25}\smash{\begin{tabular}[t]{l}$q_3$\end{tabular}}}}%
    \put(0,0){\includegraphics[width=\unitlength,page=2]{exampleB2.pdf}}%
    \put(0.39114503,0.34764446){\color[rgb]{0,0,0}\makebox(0,0)[lt]{\lineheight{1.25}\smash{\begin{tabular}[t]{l}$n_K(q_3)$\end{tabular}}}}%
    \put(0.29769503,0.22088678){\color[rgb]{0,0,0}\makebox(0,0)[lt]{\lineheight{1.25}\smash{\begin{tabular}[t]{l}$N_K(q_3)$\end{tabular}}}}%
    \put(0,0){\includegraphics[width=\unitlength,page=3]{exampleB2.pdf}}%
    \put(0.78499867,0.44734448){\color[rgb]{0,0,0}\makebox(0,0)[lt]{\lineheight{1.25}\smash{\begin{tabular}[t]{l}$p_1$\end{tabular}}}}%
    \put(0.77293082,0.52762059){\color[rgb]{0,0,0}\makebox(0,0)[lt]{\lineheight{1.25}\smash{\begin{tabular}[t]{l}$n_T(p_1)$\end{tabular}}}}%
    \put(0.77915498,0.1754309){\color[rgb]{0,0,0}\makebox(0,0)[lt]{\lineheight{1.25}\smash{\begin{tabular}[t]{l}$p_2=p_3$\end{tabular}}}}%
    \put(0.64295781,0.07383798){\color[rgb]{0,0,0}\makebox(0,0)[lt]{\lineheight{1.25}\smash{\begin{tabular}[t]{l}$n_T(p_3)$\end{tabular}}}}%
    \put(0,0){\includegraphics[width=\unitlength,page=4]{exampleB2.pdf}}%
    \put(0.80377802,0.07631715){\color[rgb]{0,0,0}\makebox(0,0)[lt]{\lineheight{1.25}\smash{\begin{tabular}[t]{l}$n_T(p_2)$\end{tabular}}}}%
    \put(0.13993722,0.38396363){\color[rgb]{0,0,0}\makebox(0,0)[lt]{\lineheight{1.25}\smash{\begin{tabular}[t]{l}$K$\end{tabular}}}}%
    \put(0.70894687,0.33873676){\color[rgb]{0,0,0}\makebox(0,0)[lt]{\lineheight{1.25}\smash{\begin{tabular}[t]{l}$T$\end{tabular}}}}%
  \end{picture}%
\endgroup%

%% file: images/normalvectorspanning.pdf_tex
\begingroup%
  \makeatletter%
  \providecommand\color[2][]{%
    \errmessage{(Inkscape) Color is used for the text in Inkscape, but the package 'color.sty' is not loaded}%
    \renewcommand\color[2][]{}%
  }%
  \providecommand\transparent[1]{%
    \errmessage{(Inkscape) Transparency is used (non-zero) for the text in Inkscape, but the package 'transparent.sty' is not loaded}%
    \renewcommand\transparent[1]{}%
  }%
  \providecommand\rotatebox[2]{#2}%
  \newcommand*\fsize{\dimexpr\f@size pt\relax}%
  \newcommand*\lineheight[1]{\fontsize{\fsize}{#1\fsize}\selectfont}%
  \ifx\svgwidth\undefined%
    \setlength{\unitlength}{313.14881727bp}%
    \ifx\svgscale\undefined%
      \relax%
    \else%
      \setlength{\unitlength}{\unitlength * \real{\svgscale}}%
    \fi%
  \else%
    \setlength{\unitlength}{\svgwidth}%
  \fi%
  \global\let\svgwidth\undefined%
  \global\let\svgscale\undefined%
  \makeatother%
  \begin{picture}(1,0.60580175)%
    \lineheight{1}%
    \setlength\tabcolsep{0pt}%
    \put(0,0){\includegraphics[width=\unitlength,page=1]{normalvectorspanning.pdf}}%
    \put(0.158073,0.13116923){\color[rgb]{0,0,0}\makebox(0,0)[lt]{\lineheight{1.25}\smash{\begin{tabular}[t]{l}$q_1$\end{tabular}}}}%
    \put(0.15637946,0.46394985){\color[rgb]{0,0,0}\makebox(0,0)[lt]{\lineheight{1.25}\smash{\begin{tabular}[t]{l}$q_2$\end{tabular}}}}%
    \put(0.3038409,0.27342659){\color[rgb]{0,0,0}\makebox(0,0)[lt]{\lineheight{1.25}\smash{\begin{tabular}[t]{l}$q_3$\end{tabular}}}}%
    \put(0.51951948,0.46078318){\color[rgb]{0,0,0}\makebox(0,0)[lt]{\lineheight{1.25}\smash{\begin{tabular}[t]{l}$p_1$\end{tabular}}}}%
    \put(0.86590219,0.46165631){\color[rgb]{0,0,0}\makebox(0,0)[lt]{\lineheight{1.25}\smash{\begin{tabular}[t]{l}$p_1'$\end{tabular}}}}%
    \put(0.87475799,0.29433119){\color[rgb]{0,0,0}\makebox(0,0)[lt]{\lineheight{1.25}\smash{\begin{tabular}[t]{l}$p_2$\end{tabular}}}}%
    \put(0.50950054,0.29807993){\color[rgb]{0,0,0}\makebox(0,0)[lt]{\lineheight{1.25}\smash{\begin{tabular}[t]{l}$p_3$\end{tabular}}}}%
    \put(0.68665152,0.19853124){\color[rgb]{0,0,0}\makebox(0,0)[lt]{\lineheight{1.25}\smash{\begin{tabular}[t]{l}$T$\end{tabular}}}}%
    \put(0.02392167,0.39589931){\color[rgb]{0,0,0}\makebox(0,0)[lt]{\lineheight{1.25}\smash{\begin{tabular}[t]{l}$K$\end{tabular}}}}%
    \put(0.16244068,0.55951996){\color[rgb]{0,0,0}\makebox(0,0)[lt]{\lineheight{1.25}\smash{\begin{tabular}[t]{l}$n_K(q_2)$\end{tabular}}}}%
    \put(0.16061337,0.02871007){\color[rgb]{0,0,0}\makebox(0,0)[lt]{\lineheight{1.25}\smash{\begin{tabular}[t]{l}$n_K(q_1)$\end{tabular}}}}%
    \put(0.37412271,0.31826114){\color[rgb]{0,0,0}\makebox(0,0)[lt]{\lineheight{1.25}\smash{\begin{tabular}[t]{l}$n_K(q_3)$\end{tabular}}}}%
    \put(0,0){\includegraphics[width=\unitlength,page=2]{normalvectorspanning.pdf}}%
  \end{picture}%
\endgroup%

%% file: images/notranslation.pdf_tex
\begingroup%
  \makeatletter%
  \providecommand\color[2][]{%
    \errmessage{(Inkscape) Color is used for the text in Inkscape, but the package 'color.sty' is not loaded}%
    \renewcommand\color[2][]{}%
  }%
  \providecommand\transparent[1]{%
    \errmessage{(Inkscape) Transparency is used (non-zero) for the text in Inkscape, but the package 'transparent.sty' is not loaded}%
    \renewcommand\transparent[1]{}%
  }%
  \providecommand\rotatebox[2]{#2}%
  \newcommand*\fsize{\dimexpr\f@size pt\relax}%
  \newcommand*\lineheight[1]{\fontsize{\fsize}{#1\fsize}\selectfont}%
  \ifx\svgwidth\undefined%
    \setlength{\unitlength}{305.21601367bp}%
    \ifx\svgscale\undefined%
      \relax%
    \else%
      \setlength{\unitlength}{\unitlength * \real{\svgscale}}%
    \fi%
  \else%
    \setlength{\unitlength}{\svgwidth}%
  \fi%
  \global\let\svgwidth\undefined%
  \global\let\svgscale\undefined%
  \makeatother%
  \begin{picture}(1,0.37453812)%
    \lineheight{1}%
    \setlength\tabcolsep{0pt}%
    \put(0,0){\includegraphics[width=\unitlength,page=1]{notranslation.pdf}}%
    \put(0.3819319,0.07554272){\color[rgb]{0,0,0}\makebox(0,0)[lt]{\lineheight{1.25}\smash{\begin{tabular}[t]{l}$q_1$\end{tabular}}}}%
    \put(0.18506163,0.20547138){\color[rgb]{0,0,0}\makebox(0,0)[lt]{\lineheight{1.25}\smash{\begin{tabular}[t]{l}$q_2$\end{tabular}}}}%
    \put(0.54350372,0.20264434){\color[rgb]{0,0,0}\makebox(0,0)[lt]{\lineheight{1.25}\smash{\begin{tabular}[t]{l}$q_3$\end{tabular}}}}%
    \put(0.42310912,0.24040648){\color[rgb]{0,0,0}\makebox(0,0)[lt]{\lineheight{1.25}\smash{\begin{tabular}[t]{l}$K$\end{tabular}}}}%
    \put(0.61504847,0.14902225){\color[rgb]{0,0,0}\makebox(0,0)[lt]{\lineheight{1.25}\smash{\begin{tabular}[t]{l}$n_K(q_3)$\end{tabular}}}}%
    \put(0.38885914,0.01819771){\color[rgb]{0,0,0}\makebox(0,0)[lt]{\lineheight{1.25}\smash{\begin{tabular}[t]{l}$n_K(q_1)$\end{tabular}}}}%
    \put(0.06433391,0.15952759){\color[rgb]{0,0,0}\makebox(0,0)[lt]{\lineheight{1.25}\smash{\begin{tabular}[t]{l}$n_K(q_2)$\end{tabular}}}}%
    \put(0,0){\includegraphics[width=\unitlength,page=2]{notranslation.pdf}}%
    \put(0.8663025,0.19769209){\color[rgb]{0,0,0}\makebox(0,0)[lt]{\lineheight{1.25}\smash{\begin{tabular}[t]{l}$T$\end{tabular}}}}%
    \put(0.90105351,0.33793762){\color[rgb]{0,0,0}\makebox(0,0)[lt]{\lineheight{1.25}\smash{\begin{tabular}[t]{l}$p_2$\end{tabular}}}}%
    \put(0.90149233,0.02451275){\color[rgb]{0,0,0}\makebox(0,0)[lt]{\lineheight{1.25}\smash{\begin{tabular}[t]{l}$p_2'$\end{tabular}}}}%
    \put(0.80548788,0.10212642){\color[rgb]{0,0,0}\makebox(0,0)[lt]{\lineheight{1.25}\smash{\begin{tabular}[t]{l}$p_3$\end{tabular}}}}%
    \put(0.80808953,0.26332036){\color[rgb]{0,0,0}\makebox(0,0)[lt]{\lineheight{1.25}\smash{\begin{tabular}[t]{l}$p_1$\end{tabular}}}}%
    \put(0.76080799,0.33421431){\color[rgb]{0,0,0}\makebox(0,0)[lt]{\lineheight{1.25}\smash{\begin{tabular}[t]{l}$n_T(p_1)$\end{tabular}}}}%
    \put(0.76763405,0.0301419){\color[rgb]{0,0,0}\makebox(0,0)[lt]{\lineheight{1.25}\smash{\begin{tabular}[t]{l}$n_T(p_3)$\end{tabular}}}}%
    \put(0.9490935,0.06628471){\color[rgb]{0,0,0}\makebox(0,0)[lt]{\lineheight{1.25}\smash{\begin{tabular}[t]{l}$n_T(p_2')$\end{tabular}}}}%
    \put(0.95504193,0.30132487){\color[rgb]{0,0,0}\makebox(0,0)[lt]{\lineheight{1.25}\smash{\begin{tabular}[t]{l}$n_T(p_2)$\end{tabular}}}}%
    \put(0,0){\includegraphics[width=\unitlength,page=3]{notranslation.pdf}}%
  \end{picture}%
\endgroup%

%% file: images/ExampleF.pdf_tex
\begingroup%
  \makeatletter%
  \providecommand\color[2][]{%
    \errmessage{(Inkscape) Color is used for the text in Inkscape, but the package 'color.sty' is not loaded}%
    \renewcommand\color[2][]{}%
  }%
  \providecommand\transparent[1]{%
    \errmessage{(Inkscape) Transparency is used (non-zero) for the text in Inkscape, but the package 'transparent.sty' is not loaded}%
    \renewcommand\transparent[1]{}%
  }%
  \providecommand\rotatebox[2]{#2}%
  \newcommand*\fsize{\dimexpr\f@size pt\relax}%
  \newcommand*\lineheight[1]{\fontsize{\fsize}{#1\fsize}\selectfont}%
  \ifx\svgwidth\undefined%
    \setlength{\unitlength}{182.58999632bp}%
    \ifx\svgscale\undefined%
      \relax%
    \else%
      \setlength{\unitlength}{\unitlength * \real{\svgscale}}%
    \fi%
  \else%
    \setlength{\unitlength}{\svgwidth}%
  \fi%
  \global\let\svgwidth\undefined%
  \global\let\svgscale\undefined%
  \makeatother%
  \begin{picture}(1,0.6813886)%
    \lineheight{1}%
    \setlength\tabcolsep{0pt}%
    \put(0,0){\includegraphics[width=\unitlength,page=1]{ExampleF.pdf}}%
    \put(0.37685752,0.61912842){\color[rgb]{0,0,0}\makebox(0,0)[lt]{\lineheight{1.25}\smash{\begin{tabular}[t]{l}$N_K(q_2)$\end{tabular}}}}%
    \put(0.37772228,0.04255755){\color[rgb]{0,0,0}\makebox(0,0)[lt]{\lineheight{1.25}\smash{\begin{tabular}[t]{l}$N_K(q_1)$\end{tabular}}}}%
    \put(0.05553357,0.43521364){\color[rgb]{0,0,0}\makebox(0,0)[lt]{\lineheight{1.25}\smash{\begin{tabular}[t]{l}$K$\end{tabular}}}}%
    \put(0.86982294,0.50607627){\color[rgb]{0,0,0}\makebox(0,0)[lt]{\lineheight{1.25}\smash{\begin{tabular}[t]{l}$p_1$\end{tabular}}}}%
    \put(0.87672785,0.16496595){\color[rgb]{0,0,0}\makebox(0,0)[lt]{\lineheight{1.25}\smash{\begin{tabular}[t]{l}$p_2=p_3$\end{tabular}}}}%
    \put(0.54715916,0.6417805){\color[rgb]{0,0,0}\makebox(0,0)[lt]{\lineheight{1.25}\smash{\begin{tabular}[t]{l}$U$\end{tabular}}}}%
    \put(0,0){\includegraphics[width=\unitlength,page=2]{ExampleF.pdf}}%
    \put(0.20477651,0.6060939){\color[rgb]{0,0,0}\makebox(0,0)[lt]{\lineheight{1.25}\smash{\begin{tabular}[t]{l}$n_K(q_2)$\end{tabular}}}}%
    \put(0,0){\includegraphics[width=\unitlength,page=3]{ExampleF.pdf}}%
    \put(0.7678735,0.33098193){\color[rgb]{0,0,0}\makebox(0,0)[lt]{\lineheight{1.25}\smash{\begin{tabular}[t]{l}$T$\end{tabular}}}}%
    \put(0.84969687,0.65831782){\color[rgb]{0,0,0}\makebox(0,0)[lt]{\lineheight{1.25}\smash{\begin{tabular}[t]{l}$n_T(p_1)$\end{tabular}}}}%
    \put(0.90446712,0.09163184){\color[rgb]{0,0,0}\makebox(0,0)[lt]{\lineheight{1.25}\smash{\begin{tabular}[t]{l}$n_T(p_2)$\end{tabular}}}}%
    \put(0,0){\includegraphics[width=\unitlength,page=4]{ExampleF.pdf}}%
    \put(0.48320253,0.24370256){\color[rgb]{0,0,0}\makebox(0,0)[lt]{\lineheight{1.25}\smash{\begin{tabular}[t]{l}$N_K(q_3)$\end{tabular}}}}%
    \put(0,0){\includegraphics[width=\unitlength,page=5]{ExampleF.pdf}}%
    \put(0.55362132,0.35449679){\color[rgb]{0,0,0}\makebox(0,0)[lt]{\lineheight{1.25}\smash{\begin{tabular}[t]{l}$n_K(q_3)$\end{tabular}}}}%
    \put(0,0){\includegraphics[width=\unitlength,page=6]{ExampleF.pdf}}%
    \put(0.27341995,0.15774781){\color[rgb]{0,0,0}\makebox(0,0)[lt]{\lineheight{1.25}\smash{\begin{tabular}[t]{l}$q_1$\end{tabular}}}}%
    \put(0.2671011,0.51735379){\color[rgb]{0,0,0}\makebox(0,0)[lt]{\lineheight{1.25}\smash{\begin{tabular}[t]{l}$q_2$\end{tabular}}}}%
    \put(0.39107681,0.37171228){\color[rgb]{0,0,0}\makebox(0,0)[lt]{\lineheight{1.25}\smash{\begin{tabular}[t]{l}$q_3$\end{tabular}}}}%
    \put(0,0){\includegraphics[width=\unitlength,page=7]{ExampleF.pdf}}%
    \put(0.6793616,0.09263322){\color[rgb]{0,0,0}\makebox(0,0)[lt]{\lineheight{1.25}\smash{\begin{tabular}[t]{l}$n_T(p_3)$\end{tabular}}}}%
    \put(0.73569069,0.00387491){\color[rgb]{0,0,0}\makebox(0,0)[lt]{\lineheight{1.25}\smash{\begin{tabular}[t]{l}$N_T(p_2)=N_T(p_3)$\end{tabular}}}}%
    \put(0,0){\includegraphics[width=\unitlength,page=8]{ExampleF.pdf}}%
    \put(0.68382267,0.56008392){\color[rgb]{0,0,0}\makebox(0,0)[lt]{\lineheight{1.25}\smash{\begin{tabular}[t]{l}$N_T(p_1)$\end{tabular}}}}%
    \put(0.20822299,0.06572785){\color[rgb]{0,0,0}\makebox(0,0)[lt]{\lineheight{1.25}\smash{\begin{tabular}[t]{l}$n_K(q_1)$\end{tabular}}}}%
  \end{picture}%
\endgroup%

%% file: images/ExampleG.pdf_tex
\begingroup%
  \makeatletter%
  \providecommand\color[2][]{%
    \errmessage{(Inkscape) Color is used for the text in Inkscape, but the package 'color.sty' is not loaded}%
    \renewcommand\color[2][]{}%
  }%
  \providecommand\transparent[1]{%
    \errmessage{(Inkscape) Transparency is used (non-zero) for the text in Inkscape, but the package 'transparent.sty' is not loaded}%
    \renewcommand\transparent[1]{}%
  }%
  \providecommand\rotatebox[2]{#2}%
  \newcommand*\fsize{\dimexpr\f@size pt\relax}%
  \newcommand*\lineheight[1]{\fontsize{\fsize}{#1\fsize}\selectfont}%
  \ifx\svgwidth\undefined%
    \setlength{\unitlength}{239.13214714bp}%
    \ifx\svgscale\undefined%
      \relax%
    \else%
      \setlength{\unitlength}{\unitlength * \real{\svgscale}}%
    \fi%
  \else%
    \setlength{\unitlength}{\svgwidth}%
  \fi%
  \global\let\svgwidth\undefined%
  \global\let\svgscale\undefined%
  \makeatother%
  \begin{picture}(1,0.59456071)%
    \lineheight{1}%
    \setlength\tabcolsep{0pt}%
    \put(0,0){\includegraphics[width=\unitlength,page=1]{ExampleG.pdf}}%
    \put(0.29492104,0.30519352){\color[rgb]{0,0,0}\makebox(0,0)[lt]{\lineheight{1.25}\smash{\begin{tabular}[t]{l}$q_2^{1/2}$\end{tabular}}}}%
    \put(0.05873219,0.29486586){\color[rgb]{0,0,0}\makebox(0,0)[lt]{\lineheight{1.25}\smash{\begin{tabular}[t]{l}$q_1^{1/2}$\end{tabular}}}}%
    \put(0.2646339,0.22615266){\color[rgb]{0,0,0}\makebox(0,0)[lt]{\lineheight{1.25}\smash{\begin{tabular}[t]{l}$q_2^{3/4}$\end{tabular}}}}%
    \put(0.08749826,0.21987997){\color[rgb]{0,0,0}\makebox(0,0)[lt]{\lineheight{1.25}\smash{\begin{tabular}[t]{l}$q_1^{3/4}$\end{tabular}}}}%
    \put(0.85521892,0.16604735){\color[rgb]{0,0,0}\makebox(0,0)[lt]{\lineheight{1.25}\smash{\begin{tabular}[t]{l}$p_1$\end{tabular}}}}%
    \put(0.85469166,0.41215945){\color[rgb]{0,0,0}\makebox(0,0)[lt]{\lineheight{1.25}\smash{\begin{tabular}[t]{l}$p_1'$\end{tabular}}}}%
    \put(0.50643216,0.27081871){\color[rgb]{0,0,0}\makebox(0,0)[lt]{\lineheight{1.25}\smash{\begin{tabular}[t]{l}$p_2$\end{tabular}}}}%
    \put(0.58481386,0.52253083){\color[rgb]{0,0,0}\makebox(0,0)[lt]{\lineheight{1.25}\smash{\begin{tabular}[t]{l}$T$\end{tabular}}}}%
    \put(0.24796826,0.39133701){\color[rgb]{0,0,0}\makebox(0,0)[lt]{\lineheight{1.25}\smash{\begin{tabular}[t]{l}$K$\end{tabular}}}}%
    \put(0,0){\includegraphics[width=\unitlength,page=2]{ExampleG.pdf}}%
  \end{picture}%
\endgroup%

%% file: images/Algorithm_grey.pdf_tex
\begingroup%
  \makeatletter%
  \providecommand\color[2][]{%
    \errmessage{(Inkscape) Color is used for the text in Inkscape, but the package 'color.sty' is not loaded}%
    \renewcommand\color[2][]{}%
  }%
  \providecommand\transparent[1]{%
    \errmessage{(Inkscape) Transparency is used (non-zero) for the text in Inkscape, but the package 'transparent.sty' is not loaded}%
    \renewcommand\transparent[1]{}%
  }%
  \providecommand\rotatebox[2]{#2}%
  \newcommand*\fsize{\dimexpr\f@size pt\relax}%
  \newcommand*\lineheight[1]{\fontsize{\fsize}{#1\fsize}\selectfont}%
  \ifx\svgwidth\undefined%
    \setlength{\unitlength}{504.2383772bp}%
    \ifx\svgscale\undefined%
      \relax%
    \else%
      \setlength{\unitlength}{\unitlength * \real{\svgscale}}%
    \fi%
  \else%
    \setlength{\unitlength}{\svgwidth}%
  \fi%
  \global\let\svgwidth\undefined%
  \global\let\svgscale\undefined%
  \makeatother%
  \begin{picture}(1,0.62446618)%
    \lineheight{1}%
    \setlength\tabcolsep{0pt}%
    \put(0,0){\includegraphics[width=\unitlength,page=1]{Algorithm_grey.pdf}}%
    \put(0.45537434,0.22805218){\color[rgb]{0,0,0}\makebox(0,0)[lt]{\lineheight{1.25}\smash{\begin{tabular}[t]{l}(a)\end{tabular}}}}%
    \put(0.50420537,0.09207658){\color[rgb]{0,0,0}\makebox(0,0)[lt]{\lineheight{1.25}\smash{\begin{tabular}[t]{l}(b)\end{tabular}}}}%
    \put(0.74986296,0.3294704){\color[rgb]{0,0,0}\makebox(0,0)[lt]{\lineheight{1.25}\smash{\begin{tabular}[t]{l}(c)\end{tabular}}}}%
    \put(0.51021538,0.5638593){\color[rgb]{0,0,0}\makebox(0,0)[lt]{\lineheight{1.25}\smash{\begin{tabular}[t]{l}(d)\end{tabular}}}}%
    \put(0.21300873,0.24856237){\color[rgb]{0,0,0}\makebox(0,0)[lt]{\lineheight{1.25}\smash{\begin{tabular}[t]{l}(e)/(f)\end{tabular}}}}%
    \put(0.4118021,0.38881888){\color[rgb]{0,0,0}\makebox(0,0)[lt]{\lineheight{1.25}\smash{\begin{tabular}[t]{l}$K$\end{tabular}}}}%
    \put(0.84902745,0.39708261){\color[rgb]{0,0,0}\makebox(0,0)[lt]{\lineheight{1.25}\smash{\begin{tabular}[t]{l}$T$\end{tabular}}}}%
    \put(0,0){\includegraphics[width=\unitlength,page=2]{Algorithm_grey.pdf}}%
  \end{picture}%
\endgroup%

%% file: images/Picture7.pdf_tex
\begingroup%
  \makeatletter%
  \providecommand\color[2][]{%
    \errmessage{(Inkscape) Color is used for the text in Inkscape, but the package 'color.sty' is not loaded}%
    \renewcommand\color[2][]{}%
  }%
  \providecommand\transparent[1]{%
    \errmessage{(Inkscape) Transparency is used (non-zero) for the text in Inkscape, but the package 'transparent.sty' is not loaded}%
    \renewcommand\transparent[1]{}%
  }%
  \providecommand\rotatebox[2]{#2}%
  \newcommand*\fsize{\dimexpr\f@size pt\relax}%
  \newcommand*\lineheight[1]{\fontsize{\fsize}{#1\fsize}\selectfont}%
  \ifx\svgwidth\undefined%
    \setlength{\unitlength}{336.1947399bp}%
    \ifx\svgscale\undefined%
      \relax%
    \else%
      \setlength{\unitlength}{\unitlength * \real{\svgscale}}%
    \fi%
  \else%
    \setlength{\unitlength}{\svgwidth}%
  \fi%
  \global\let\svgwidth\undefined%
  \global\let\svgscale\undefined%
  \makeatother%
  \begin{picture}(1,0.63722559)%
    \lineheight{1}%
    \setlength\tabcolsep{0pt}%
    \put(0,0){\includegraphics[width=\unitlength,page=1]{Picture7.pdf}}%
    \put(-0.00271594,0.01605985){\color[rgb]{0,0,0}\makebox(0,0)[lt]{\lineheight{1.25}\smash{\begin{tabular}[t]{l}$q_1$\end{tabular}}}}%
    \put(0.97357402,0.47865866){\color[rgb]{0,0,0}\makebox(0,0)[lt]{\lineheight{1.25}\smash{\begin{tabular}[t]{l}$q_2$\end{tabular}}}}%
    \put(0,0){\includegraphics[width=\unitlength,page=2]{Picture7.pdf}}%
    \put(0.40270132,0.5465898){\color[rgb]{0,0,0}\makebox(0,0)[lt]{\lineheight{1.25}\smash{\begin{tabular}[t]{l}$u$\end{tabular}}}}%
    \put(0.61115515,0.54770354){\color[rgb]{0,0,0}\makebox(0,0)[lt]{\lineheight{1.25}\smash{\begin{tabular}[t]{l}$-u$\end{tabular}}}}%
    \put(0.54906271,0.41621728){\color[rgb]{0,0,0}\makebox(0,0)[lt]{\lineheight{1.25}\smash{\begin{tabular}[t]{l}$-v$\end{tabular}}}}%
    \put(0.55133649,0.61285817){\color[rgb]{0,0,0}\makebox(0,0)[lt]{\lineheight{1.25}\smash{\begin{tabular}[t]{l}$v$\end{tabular}}}}%
    \put(0,0){\includegraphics[width=\unitlength,page=3]{Picture7.pdf}}%
    \put(0.5605224,0.07523444){\color[rgb]{0,0,0}\makebox(0,0)[lt]{\lineheight{1.25}\smash{\begin{tabular}[t]{l}$q_1$\end{tabular}}}}%
    \put(0.43983941,0.07150658){\color[rgb]{0,0,0}\makebox(0,0)[lt]{\lineheight{1.25}\smash{\begin{tabular}[t]{l}$q_2$\end{tabular}}}}%
  \end{picture}%
\endgroup%

%% file: images/Picture12.pdf_tex
\begingroup%
  \makeatletter%
  \providecommand\color[2][]{%
    \errmessage{(Inkscape) Color is used for the text in Inkscape, but the package 'color.sty' is not loaded}%
    \renewcommand\color[2][]{}%
  }%
  \providecommand\transparent[1]{%
    \errmessage{(Inkscape) Transparency is used (non-zero) for the text in Inkscape, but the package 'transparent.sty' is not loaded}%
    \renewcommand\transparent[1]{}%
  }%
  \providecommand\rotatebox[2]{#2}%
  \newcommand*\fsize{\dimexpr\f@size pt\relax}%
  \newcommand*\lineheight[1]{\fontsize{\fsize}{#1\fsize}\selectfont}%
  \ifx\svgwidth\undefined%
    \setlength{\unitlength}{229.25727774bp}%
    \ifx\svgscale\undefined%
      \relax%
    \else%
      \setlength{\unitlength}{\unitlength * \real{\svgscale}}%
    \fi%
  \else%
    \setlength{\unitlength}{\svgwidth}%
  \fi%
  \global\let\svgwidth\undefined%
  \global\let\svgscale\undefined%
  \makeatother%
  \begin{picture}(1,0.46714737)%
    \lineheight{1}%
    \setlength\tabcolsep{0pt}%
    \put(0,0){\includegraphics[width=\unitlength,page=1]{Picture12.pdf}}%
    \put(0.96904793,0.23542905){\color[rgb]{0,0,0}\makebox(0,0)[lt]{\lineheight{1.25}\smash{\begin{tabular}[t]{l}$O$\end{tabular}}}}%
    \put(0.97458888,0.03446946){\color[rgb]{0,0,0}\makebox(0,0)[lt]{\lineheight{1.25}\smash{\begin{tabular}[t]{l}$H$\end{tabular}}}}%
    \put(0.57413923,0.33824548){\color[rgb]{0,0,0}\makebox(0,0)[lt]{\lineheight{1.25}\smash{\begin{tabular}[t]{l}$T$\end{tabular}}}}%
    \put(0,0){\includegraphics[width=\unitlength,page=2]{Picture12.pdf}}%
  \end{picture}%
\endgroup%